
\pageno=0

%
%


\font\Bigtenrm=cmr10 scaled\magstep {5.5}
\font\Bigtenbf=cmbx10 scaled\magstep {2}

\vskip 6cm
\centerline{\Bigtenrm Extensions of Normed Algebras}
\vskip 3cm
\centerline{\Bigtenrm by Thomas Dawson, M.Math.}
\vskip 45mm
\centerline{\Bigtenrm Thesis submitted to}\vskip 4mm
\centerline{\Bigtenrm The University of Nottingham} \vskip 4mm
\centerline{\Bigtenrm for the degree of}        \vskip 4mm
\centerline{\Bigtenrm Doctor of Philosophy,}   \vskip 4mm
\centerline{\Bigtenrm June 2003}   \vskip 3cm
\centerline{(Minor corrections and additions have been made. The page layout
has also been changed.)}
\vfill

\eject

\pageno=1

\font\bigtenrm=cmr10 scaled\magstep {2.5}
\font\Bigtenbf=cmbx10 scaled\magstep {3.5}

\def\leaderfill{\leaders\hbox to 1em{\hss.\hss}\hfill}



\def\C{{\bf C}}	
\def\N{{\bf N}}	
\def\Z{{\bf Z}}	
\def\R{{\bf R}}	



\def\prin{\bigl(\alpha (x)\bigr)}

\def\sqr#1#2{{\vcenter{\hrule height.#2pt
	\hbox{\vrule width.#2pt height#1pt \kern #1pt
		\vrule width.#2pt}
	\hrule height.#2pt}}}
\def\square{\mathchoice\sqr77\sqr77\sqr77\sqr77}
\def\udot{\cdot\hskip -7pt\cup}		

\font\bigtenrm=cmr10 scaled\magstep {2.5}
\font\bigtenbf=cmbx10 scaled\magstep {2.5}
\font\medtenbf=cmbx10 scaled\magstep {1.7}

\def\V{{\cal V}}
\def\No{\N_0}
\def\eop{\hfill $\square$\vskip10pt}
\def\thm{\vskip 10pt\noindent THEOREM }
\def\dfn{\vskip 10pt\noindent DEFINITION }
\def\ex{\vskip 10pt\noindent EXAMPLE }
\def\cor{\vskip 10pt\noindent COROLLARY }
\def\lem{\vskip 10pt\noindent LEMMA }
\def\prop{\vskip 10pt\noindent PROPOSITION }
\def\pf{\par\noindent {\it Proof.} }
\def\ref#1{[{\bf #1}]}

\def\npar{\vskip 10pt\noindent}

\def\na{{\bf na}}
\def\Ba{{\bf Ba}}
\def\ua{{\bf ua}}

\def\bGA{\overline{G(A)}}
\def\xb{\bar x}
\def\pri{{^\prime}}
\def\htt{\char'136}
\def\qs#1#2{q\s_{#1}(#2)}
\def\T#1{\Theta(#1)}
\def\int{^\circ}
\def\ov#1{\overline{#1}}
\def\c#1{\widetilde{#1}}
\def\snc#1{\overline{#1\htt}}
\def\S#1{\check S\left(#1\right)}
\def\inv{{^{-1}}}

\def\Ac{\tilde A}
\def\aA{A^\a}
\def\ph{\varphi}

\def\k{\kappa}
\def\l{\lambda}
\def\e{\varepsilon}
\def\ths#1{\theta_{#1}}
\def\th{\theta}
\def\d{\delta}
\def\g{\gamma}
\def\G{\Gamma}
\def\r{\rho}
\def\s{\sigma}
\def\p{\pi}
\def\t{\tau}
\def\m{\mu}
\def\u{\upsilon}
\def\f{\phi}
\def\o{\omega}
\def\ps{\pi^*}
\def\pss#1#2{\ps_{#1,#2}}
\def\Ps#1#2{\pi_{#1,#2}}
\def\tss#1#2{\ths{#1,#2}}
\def\Om#1{\Omega(#1)}
\def\inc{\iota}

\def\J{{\cal J}}
\def\U{{\cal U}}
\def\ao{{\alpha_0}}
\def\oo{{\o_0}}
\def\oi{{\o_1}}
\def\ui{{\u_1}}
\def\fs{^{<\oo}}

\def\ker{{\rm ker}\,}	
\def\im{{\rm im}\,}	
\def\OA{\Omega(A)}	
\def\ha{\hat a}		
\def\hA{\hat A}		
\def\eh#1{\varepsilon_{#1}}	
\def\id#1{{\rm id}_{#1}}

\def\a{\alpha}
\def\b{\beta}
\def\O{\Omega}		
\def\Oa{\O_\alpha}	
\def\Ou{\O_\u}
\def\Ot{\O_\t}
\def\Aa{A_{\alpha}}
\def\prin{(\a(x))}
\def\aA{A^\a}
\def\norm#1{\left\Vert#1\right\Vert}
\def\ab#1{\left\vert#1\right\vert}
\def\map#1#2{\colon#1\to#2}
\def\mapto#1#2#3#4{\colon#1\to#2;\,#3\mapsto#4}
\def\st{\,\colon\;}
\def\set#1{\left\{#1\right\} }
\def\AU{A_\U}
\def\UA{A^\U}
\def\OU{\O_\U}
\def\aiu{{\a\in\U}}
\def\Om#1{\O(#1)}
\def\ol{{(\o,\l)}}
\def\kl{{(\k,\l)}}

\def\Nb#1{{\bf N}\left(#1\right)}
\def\Nbb#1#2{{\bf N}_{#1}\left(#2\right)}
\def\Ho#1{H^1(#1,\Z)}
\def\Hn#1#2{H^{#1}(#2,\Z)}
\def\fib#1{\pi\inv\left(#1\right)}
\def\ri{{\rm i}}
\def\rest#1{\vert_{#1}}


\def\mapright#1{\smash{
	\mathop{\longrightarrow}\limits^{#1}}}
\def\mapleft#1{\smash{
	\mathop{\longleftarrow}\limits_{#1}}}
\def\mapdown#1{\Big\downarrow
	\rlap{$\vcenter{\hbox{$\scriptstyle#1$}}$}}
\def\mapup#1{\Big\uparrow
	\rlap{$\vcenter{\hbox{$\scriptstyle#1$}}$}}
\def\mapne#1{\nearrow${$\scriptstyle#1$}$}
\def\mapnw#1{ \nwarrow${$\scriptstyle{#1}$}$  }
\def\mapse#1{\searrow${$\scriptstyle#1$}$}
\def\mapsw#1{ \swarrow${$\scriptstyle{#1}$}$  }

\def\dl#1{{\lim_\rightharpoondown}\mathstrut_{#1}}	
\def\il#1{{\lim_\leftharpoondown}\mathstrut_{#1}}		


\noindent {\Bigtenbf Contents}

\vskip 16pt\noindent
{\bf Contents\hfill 1}

\vskip 9pt\noindent
{\bf Abstract\hfill 2}

\vskip 9pt
\noindent
{\bf Chapter 1: Constructions and Elementary Properties and Applications of Algebraic Extensions
of Normed Algebras\hfill 3}\vskip 4pt
\line{1. Introduction\leaderfill 3}
\line{2. Construction of Extensions\leaderfill 4}
\line{3. Comparison of the Extensions\leaderfill 12}
\line{4. Maps back from the Extensions\leaderfill 17}
\line{5. The \v Silov and Choquet Boundaries of the Extensions\leaderfill 23}

\vskip 9pt
\noindent
{\bf Chapter 2: Further Properties and Applications of Algebraic Extensions\hfill 29}\vskip 4pt
\line{1. Integrally Closed Rings\leaderfill 29}
\line{2.  Regularity and Localness \leaderfill 30}

\vskip 9pt
\noindent
{\bf Chapter 3: Algebraic Extensions and Topology\hfill 35}\vskip 4pt
\line{1. Dimension\leaderfill 35}
\line{2. \v Cech Cohomology\leaderfill 36}

\vskip 9pt
\noindent
{\bf Chapter 4: The Invertible Group of an Extension of a Normed Algebra\hfill 41}\vskip 4pt
\line{1. Introductory Results and Examples\leaderfill 41}
\line{2. Algebraic Extensions\leaderfill 44}
\line{3. Context: Uniform Algebras\leaderfill 48}

\vskip 9pt
\noindent
{\bf Chapter 5: Uniform Algebras with Connected, Dense Invertible Groups\hfill 52}\vskip 4pt
\line{1. Uniform Algebras on Compact Plane-Sets\leaderfill 52}
\line{2. Systems of Cole Extensions\leaderfill 53}
\line{3. Logarithmic Extensions\leaderfill 54}
\line{4. Conclusion\leaderfill 65}


\vskip 9pt
\noindent
{\bf Appendix 1:} Gelfand Theory\hfill {\bf 67}

\vskip 9pt
\noindent
{\bf Appendix 2:} Direct and Inverse Limits\hfill {\bf 70}

\vskip 9pt
\noindent
{\bf Appendix 3:} Algebra\hfill {\bf 73}

\vskip 9pt
\noindent
{\bf References \hfill 75}

\vfill\eject


\noindent {\bigtenbf Abstract}

\vskip 25pt
\noindent We review and analyse techniques from the literature for extending a normed algebra, $A$,
to a normed algebra, $B$, so that $B$ has interesting or desirable properties which
$A$ may lack. For example, $B$ might include roots of monic polynomials over $A$.

These techniques have been important historically for constructing examples in the
theory of Banach algebras. In Section 2.2.2, we construct a new example in this way. Elsewhere
in Chapter 2, we contribute to the related programme of determining which properties of an
algebra are shared by certain extensions of it.

Similarly, we consider the relations between $\OA$ and $\Om B$, the
topological spaces of closed, maximal ideals of $A$ and $B$
respectively. For example, it is shown that if $\Om B$ has trivial
first \v Cech-cohomology group, then so has $\OA$
if $B$ is one of the extensions of $A$ constructed in Chapter 1.

The invertible group of a normed algebra is studied in Chapter 4; it is shown that if
a Banach algebra, $A$,
has dense invertible group then so has
every integral extension of $A$. The context for this work is also explained: some
new results characterising trivial uniform algebras by means of approximation by invertible elements
are given. We show how these results partially answer a famous, open problem of Gelfand.

Results in Chapter 4 lead to the conjecture that a uniform algebra is trivial if the group
of exponentials of its elements is dense in the algebra. We investigate
this conjecture in Chapter 5. In the search for a counterexample, we
construct and establish some properties of `logarithmic extensions' of a regular uniform algebra.

\vfill
\npar
{\bf Disclaimer.} Except where otherwise stated, all the results in this thesis
are the author's.

Some parts of Chapters 1 and 4 have appeared in $\ref{DawS}$
and $\ref{DF}$ respectively.
\npar
{\bf Acknowledgement.} The author is very grateful to Joel Feinstein for useful suggestions,
guidance with the literature, and
comments on this work.

The author also thanks the EPSRC for providing financial support for
this work.

\vfill\eject



\noindent {\bigtenbf Chapter 1}
\vskip 10pt
\noindent {\bigtenrm Constructions and Elementary Properties and Applications of Algebraic Extensions
of Normed Algebras}
\vskip 25pt


\def\xma#1{\xi_{\a_{#1}}^{m_{#1}}}
\def\xmaj{\xma j}
\def\tma#1{t_{\a_{#1}}^{m_{#1}}}

\def\smaj#1{s_{\a_{#1}}^{m_{#1}}}

\def\naj{n(\a_j)}


\vskip 10pt

\noindent In this chapter we synthesise the methods from the literature
for extending normed algebras so as to include roots of sets of monic polynomials over the
original normed algebra. In Sections 3 to 5 we establish some elementary properties
of the extensions.

Throughout the thesis, `Proposition 1.2.3.4' will be abbreviated to `Proposition 4'
within Section 1.2.3, and so on.


\vskip 15pt
\noindent{\medtenbf 1.1. Introduction}

\vskip 15pt
\noindent{\bf 1.1.1. Context and Notation}\vskip 10pt

\noindent We assume that the reader is familiar with the
basic results of functional analysis and in particular the general theory of normed algebras.
A good reference for the latter subject is Chapter 2 of $\ref{Dal}$.

In this work, $A$ will, except where otherwise stated, denote a commutative,
complex, normed, associative algebra. This means that $A$ is a non-empty set together
with  maps
$$\eqalign{&A\times A\to A\;;\;(a_1,a_2)\mapsto a_1+a_2\qquad\qquad\hbox{`addition',}\cr
&A\times A\to A\;;\;(a_1,a_2)\mapsto a_1\cdot a_2\qquad\qquad\hbox{`multiplication',}\cr
&\C\times A\to A\;;\;(\l,a)\mapsto \l a\qquad\qquad\hbox{`scalar multiplication',}\cr
&A\to [0,+\infty)\;;\;a\mapsto \norm a_A\qquad\qquad\hbox{`the norm',}\cr}$$
such that $A$ is an algebra over $\C$ and a normed space over $\C$. Axioms for the basic
algebraic concepts are given in Appendix 3.
Multiplication is often written as juxtaposition
and the subscript on the norm will often be omitted.
Algebra norms are assumed to satisfy the following:
$$\norm{a_1a_2}\le\norm{a_1}\norm{a_2}\qquad\qquad(a_1,a_2\in A).$$

If $A$ is an algebra, we shall assume it is unital: it has an identity, $1_A$, such that
$$\norm {1_A}=1.$$ When no confusion is likely to result we shall omit the subscript $A$.

A {\it Banach algebra} is a normed
algebra which is {\it complete}: every Cauchy sequence in $A$ converges.
The completion (see Appendix 1) of a normed algebra, $A$, is written $\c A$.

In this thesis we shall use the term `algebraic extension' to mean an
extension constructed by one of the methods in Section 1.2.
Historically (see $\ref{DawS}$) the study of algebraic extensions
in functional analysis appears to have been mainly
motivated by the search for examples in the theory of uniform algebras. We shall define, $\ua$,
the category of uniform
algebras in Section 1.1.3.
The categories of normed and Banach algebras will sometimes be
abbreviated to $\na$ and $\Ba$ respectively. No knowledge of category
theory is necessary to understand this thesis.

A mapping $\th\map AB$ between
algebras is called {\it unital} if $\th(1_A)=1_B$. We shall assume that the homomorphisms
in $\na$ are continuous and unital.

For future reference we note that
$\th$ induces a unital homomorphism, $\th_\#$, between the associated algebras
of polynomials
$$\th_\#\mapto{A[x]}{B[x]}{\sum_{j=0}^m a_mx^m}{\sum_{j=0}^m \th(a_m)x^m}.$$
Furthermore, if $\a(x)\in A[x]$ then we shall often write $\th(\a)(x)$ for
$\th_\#(\a(x))$.

Another example of one map inducing a second will frequently be used during this thesis.
We describe the situation now. Suppose that $\th\map{E_1}{E_2}$ is a map and
that $F_1$ and $F_2$ are sets. Suppose that $E_j^*$ is a non-empty subset of $F_j^{E_j}$
($j=1,2$). Suppose further that for each $f\in E_2^*$ we have $f\circ\th\in E_1^*$.
Then $\th$ induces an adjoint map
$$\th^*\mapto{E_2^*}{E_1^*}f{f\circ\th}.$$
Strictly speaking, the notation for this adjoint should refer to $E_1^*$ and $E_2^*$.
We shall always make clear the domain and codomain of adjoints wherever they occur. Note
that the usual situation of the adjoint map in functional analysis (that is,
in which $E_1$ and
$E_2$ are normed spaces and $E_j^*$ is the set of bounded, linear functionals on
$E_j$ ($j=1,2$)) will not arise in this thesis.

Let $E$ be a set. Then $E\fs$ will denote the set of finite subsets of $E$.
If $E$ is a subset of a ring then $(E)$ will denote the ideal generated by $E$ in that ring.
If $E$ is a subset of a topological space then $E^\circ$ and $\overline E$ denote the
interior and closure respectively of $E$ with respect to the topology. It is
convenient to mention
now that
neighbourhoods need not be open in this thesis. For a
normed space, $A$,  $a\in A$, and $r>0$ we
set
$$\eqalign{B_A[a,r]&:=\set{b\in A\st\norm{b-a}\le r}\qquad\hbox{and}\cr
B_A(a,r)&:=\set{b\in A\st\norm{b-a}< r}.\cr}$$
If the subscript $A$ is omitted, the normed space should be taken to be $\C$.

\vskip 15pt
\noindent{\bf 1.1.2. Gelfand Theory}\vskip 10pt

\noindent 
We relegate a summary of Gelfand theory
to Appendix 1; the reader may wish to glance at this now to check
the conventions and notation which will be used. In particular, $\OA$ will stand
for the space of continuous characters of the normed algebra $A$.
When $\O$ appears on its own it will refer to $A$.

\vskip 15pt
\noindent{\bf 1.1.3. Uniform Algebras}\vskip 10pt

\noindent There is a voluminous literature on uniform algebras.
Introductions can be found in $\ref{Lei}$, $\ref{Bro}$, and
$\ref{Sto}$.

If $X$ is a topological space we
write $C(X)$ for the algebra of all continuous functions $X\to\C$.
Recall that a set of functions $S\subseteq\C^X$  is said to {\it
separate the points of $X$} if for each pair of distinct points $x_1,
x_2\in X$ there exists $f\in S$ such that $f(x_1)\ne f(x_2)$.

\dfn 1.1.3.1. A {\it uniform algebra}, $A$, is a subalgebra of $C(X)$ for some
compact, Hausdorff space $X$ such that $A$ is closed with respect to the supremum
norm, separates the points of $X$, and contains the constant functions. We speak
of `the uniform algebra $(A,X)$'. The notation $\varepsilon$ will be
used for the map sending $\k\in X$ to the character $\eh\k\in\O$ given by
$\eh\k(f)=f(\k)\ (f\in A)$. The uniform algebra is {\it natural} if
$\eh{}$ is surjective,
and $A$ is called {\it trivial} if $A=C(X)$.

\npar  It is standard that (see $\ref{Lei}$, p. 28), all trivial uniform algebras
are natural.

An important question in the theory of uniform algebras is which properties of $(A,X)$
force $A$ to be trivial. For example it is sufficient that $A$ be self-adjoint, by
the Stone-Weierstrass theorem. In $\ref{Col}$ an example is given of a non-trivial
uniform algebra, $(B,X)$, which is natural and such that
every point of $X$ is a `peak-point' (see Section 1.5.1). It had previously been conjectured
that no such algebra existed. The example was constructed using algebraic extensions.


\def\ts{\th_\#}
\def\xma#1{\xi_{\a_{#1}}^{m_{#1}}}
\def\xmaj{\xma j}
\def\tma#1{t_{\a_{#1}}^{m_{#1}}}

\def\smaj#1{s_{\a_{#1}}^{m_{#1}}}

\def\naj{n(\a_j)}

\vskip 15pt
\noindent {\medtenbf 1.2. Construction of Extensions}

\vskip 15pt
\noindent {\bf 1.2.1. Arens-Hoffman Extensions}\vskip 10pt

\noindent Let $\a(x)=a_0+\cdots+a_{n-1}x^{n-1}+x^n$
be a monic polynomial over the algebra $A$. The basic construction arising from
$A$ and $\a(x)$ is the `Arens-Hoffman extension', $\Aa$. This
was introduced in $\ref{Are}$ and will be described after the following definition.

\dfn 1.2.1.1.
Let $A$ be a normed algebra. An {\it extension}
of $A$ is a commutative, unital normed algebra, $B$, together
with an isometric monomorphism $\th\map AB$.

\npar The {\it Arens-Hoffman extension} of $A$ with respect to $\a(x)$ is
the algebra $\Aa:=A[x]/(\a(x))$ under a certain norm; the embedding is
given by the map $\nu\,\colon\,a\mapsto(\a(x))+a$.

To simplify notation, we shall let $\xb$ denote the coset of $x$ and
often omit the indeterminate when using a polynomial as an index.

It is a purely algebraic fact ($\ref{Are}$) that each element of $\Aa$ has a unique
representative of degree less than $n$, the degree of $\a(x)$. Arens
and Hoffman proved that, provided the positive
number $t$ satisfies the inequality $t^n\ge\sum_{k=0}^{n-1}\norm{a_k}t^k$,
then
$$\norm{ \sum_{k=0}^{n-1}b_k\xb^k}=
\sum_{k=0}^{n-1}\norm{b_k}t^k\qquad(b_0,\ldots,b_{n-1}\in A)$$
defines an algebra norm on $\Aa$. We shall refer to $t$ as the (Arens-Hoffman)
norm-parameter of $\Aa$. From now on we shall usually identify $A$ with its canonical
image in $\Aa$.

\npar Our first proposition shows that Arens-Hoffman extensions satisfy a
certain universal property which is very useful when investigating algebraic
extensions. It is not explicitly stated anywhere in the literature (excepting
$\ref{4DIS}$) and seems to be taken as obvious.

\prop 1.2.1.2. Let $A$ be a normed algebra and
let $\th\map {A}{B^{(2)}}$ be a homomorphism of
normed algebras. Let $\a_1(x)=a_0+\cdots+a_{n-1}x^{n-1}+x^n\in A[x]$
and $B^{(1)}=A_{\a_1}$ with norm-parameter $t$.
Let $y\in B^{(2)}$ be a root
of the polynomial $\a_2(x):=\th(\a_1)(x)$ (see Section 1.1.1 for this
notation). Then there is a unique homomorphism
$\f\map{B^{(1)}}{B^{(2)}}$ such that $\phi(\xb)=y$ and the following diagram
is commutative
$$\matrix{B^{(1)}&\mapright{\phi}&B^{(2)}\cr
	\mapup{\nu}&\mapne{\theta}&\cr
	A&&\cr}$$
The map $\f$ is continuous if and only if $\th$ is continuous.

\pf This is elementary; we provide a proof for completeness.

The map $\th$ induces a homomorphism
$\ts$ as above. Let $\phi^0$ be the map $\ts$ followed by evaluation
at $y$. Thus $\phi^0(\b(x))=\th(\b)(y)$.
By the assumption on $y$, the principal ideal, $(\a_1(x))$, in $A[x]$ generated
by $\a_1(x)$ is contained in the kernel of $\phi^0$. This map therefore
induces a homomorphism
$$\phi\mapto{\Aa}{B^{(2)}}{\b(\xb)}{\th(\b)(y)}.$$
Clearly $\theta=\phi\circ\nu$. It is also clear that $\phi$ is unique.
If $\th$ is continuous then so is $\phi$ for if
$\b(x)=\sum_{j=0}^{n-1}b_jx^j\in A[x]$ then we see
$$\eqalign{\norm{\phi(\b(\xb))} &= \norm{\sum_{j=0}^{n-1}\th(b_j)y^j}\cr
&\le\norm\th\sum_{j=0}^{n-1}\norm{b_j}\norm{y}^j\cr
&\le M\sum_{j=0}^{n-1}\norm{b_j}t^j\cr
&= M\norm{\b(\xb)},\cr}$$
where $M=\norm\th\max\left((\norm y/t)^j\st j=0,\ldots,n-1\right)$.

Conversely it is easy to see that $\th$ must be bounded if $\phi$ is.
\eop

\vskip 15pt
\noindent {\bf 1.2.2. Extensions by Infinite Sets of Polynomials}\vskip 10pt
\noindent
We now introduce the methods for extending $A$ by an infinite set, $\U$, of monic
polynomials over $A$. For each $\a\in\U$, $n(\a)$ will stand for the degree of $\a(x)$.
In this section we shall make extensive use of transfinite
methods and set theory. The reader is referred to $\ref{Hal}$ for this
subject.

\vskip 10pt\noindent{\bf Standard Extensions}
\npar
The basic construction which generalises the Arens-Hoffman
extension is the standard extension.
Standard extensions solve the problem of adjoining roots of elements of $\U$ in the
category of algebras; Lindberg showed how to extend the norms. The constructions of Narmania
and Cole perform the same duty in the categories $\Ba$ and $\ua$ respectively.

In particular, Cole's extensions are uniform algebras: it is well-known that
an Arens-Hoffman extension of a uniform algebra need not be a uniform algebra. To see this
one can use the criterion for semisimplicity given in $\ref{Are}$ to construct Arens-Hoffman
extensions of a uniform algebra which are not semisimple; it is obvious that all
uniform algebras are semisimple.

Standard extensions are defined in the following theorem of Lindberg.

\thm 1.2.2.1 ($\ref{LinIE}$). Let $A$ be a normed algebra and $\U$ be a set of monic
polynomials over $A$. Let $\le$ be a well-ordering on $\U$ with least element
$\ao$. Then there exists a normed algebra, $B_\U$, with a family of
subalgebras, $(B_\a)_{\a\in\U}$, such that $B_\U=\bigcup_{\a\in\U}B_\a$ and
\item{(i)} for all $\a,\b\in\U$, $B_\a\subseteq B_\b$ if $\a\le\b$, and,
\item{(ii)} for all $\b\in\U$, there is an isometric isomorphism, $\psi_\b$, from
the Arens-Hoffman extension of $B_{<\b}$ by $\b(x)$ onto $B_\b$, where
$$B_{<\b}:=\cases{ \bigcup_{\a<\b} B_\a & if $\ao<\b$, and\cr
A & if $\ao=\b$.\cr}$$
The map $\psi_\b$ satisfies $\psi_\b(b)=b$ for all $b\in B_{<\b}$.
\pf See $\ref{LinIE}$.\eop

\npar We shall refer to $\xi_\b=\psi_\b(\xb)$ as
the {\it standard root} of $\b(x)$ in $B_\U$. The norm of $\xi_\b$
is equal to $t_\b$, the norm-parameter chosen for $(B_{<\b})_\b$.
Unless otherwise stated, when we discuss standard extensions
the notation will be as it is here.

Again the following two results have not
been explicitly  given in the literature but,
apart from the continuity condition in Proposition 3, are probably regarded as
obvious by those working in the field. We state the results in full as we shall
make use of them later.

\lem 1.2.2.2. Let $B_\U$ be a standard extension of the normed algebra $A$
with respect to the well-ordered set of monic polynomials, $\U$. For each $b\in B_\U$
there exist $S=\set{\a_1,\ldots,\a_N}\in\U\fs$ and
$b_m\in A$ ($m\in L:=\times_{j=1}^N\set{0,1,\ldots,n(\a_j)-1}$) such that
$$\eqalign{b&=\sum_{m\in L} b_m \xi_{\a_1}^{m_1}\cdots\xi_{\a_N}^{m_N}\qquad\hbox{and}\cr
\norm b&=\sum_{m\in L} \norm{b_m}
t_{\a_1}^{m_1}\cdots t_{\a_N}^{m_N}.\cr}$$

\pf
Let, with $L$ depending on $\a_1,\ldots,\a_N$ as above,
$$\eqalign{ \J=\Biggl\{ \a\in\U \st &\hbox{ for all }b\in B_\a\hbox{ there exist }
\a_1,\ldots,\a_N\le\a\hbox{ and }
b_m\in A\ (m\in L)\cr
&\quad\hbox{such that }b=\sum_{m\in L}b_m\xma 1\cdots\xma N\hbox{ and }
\norm b
=\sum_{m\in L} \norm{b_m}\tma 1\cdots \tma N\Biggr\}.\cr}$$

Evidently $\ao\in\J$.
Let $\b>\ao$ and assume that $[\ao,\b)\subseteq\J$. Let $b\in B_\b$. By the
construction of standard extensions,
there exist unique $c_0,\ldots,c_{n(\b)-1}\in B_{<\b}$ such that
$$b=\sum_{j=0}^{n(\b)-1} c_j\xi_\b^j.$$
Since $\b>\ao$, $B_{<\b}=\cup_{\a<\b}B_\a$ and by hypothesis there exist
$\a_1,\ldots,\a_N<\b$  such that for each $j\in\set{0,\ldots,n(\b)-1}$,
$c_j$ has a representation of the form
$$c_j=\sum_{m\in L}b_{m,j}\xma 1\cdots\xma N\qquad(b_{m,j}\in A),$$
where $L=\times_{j=1}^N\set{0,1,\ldots,n(\a_j)-1}$ and
$\norm{c_j}=\sum_{m\in L}\norm{b_{m,j}}\tma1\cdots\tma{N}$. Set $\a_{N+1}=\b$ and
$L\pri=\times_{j=1}^{N+1}\set{0,1,\ldots,n(\a_j)-1}$. Thus (admitting
a slight abuse in the notation for the brackets in the subscripts)
$$\eqalign{b&=\sum_{m\in L\pri}b_m\xma 1\cdots\xma {N+1}\qquad\hbox{and}\cr
\norm b&=\sum_{m_{N+1}=0}^{n(\a_{N+1})-1} \norm{c_{m_{N+1}}}\tma{N+1}\qquad\hbox{(by definition of the standard norm)}\cr
&=\sum_{m_{N+1}=0}^{n(\a_{N+1})-1} \sum_{m\in L}\norm{b_{m,m_{N+1}}}\tma1\cdots\tma{N+1}\qquad\hbox{(by inductive hypothesis)}\cr
&=\sum_{m\in L\pri}\norm{b_m}\tma1\cdots\tma{N+1}.\cr}$$

Therefore $\b\in\J$ and by the transfinite induction theorem $\J=\U$ as required.\eop

\prop 1.2.2.3.
Let $A$ be a normed algebra and $\U$ be a non-empty set of monic polynomials over
$A$. Let $B^{(1)}=B_\U$ be a standard extension
of $A$ with respect to $\U$ and
$\th\map{A}{B^{(2)}}$ be a unital homomorphism of normed algebras.
Let $\xi_\a$ be the standard root of ${\a\in\U}$, with associated
norm parameter $t_\a$, and suppose
$(\eta_\a)_{\a\in\U}\subseteq B^{(2)}$ is such that $\th(\a)(\eta_\a)=0$
for all $\a\in\U$.
Then there is a unique, unital homomorphism $\phi\map{B^{(1)}}{B^{(2)}}$ such that
for all $\a\in\U$ $\phi(\xi_\a)=\eta_\a$ and the
following diagram is commutative
 $$\matrix{ B^{(1)} &\mapright{\phi} &B^{(2)}	\cr
	\mapup{\subseteq} &\mapne{\theta}&	\cr
	A &&	\cr}$$
The map $\phi$ is continuous if $\th$ is continuous and
$$\sum_\aiu \left(n(\a)-1\right)\log^+\left({\norm{\eta_\a}\over{t_\a}}\right)<+\infty.$$
(Here $\log^+$ denotes the positive part of the logarithm, $\max(\log ,0)$.) In this
case $$\norm\th\le\sup\left(\prod_{\a\in U}
\left(\ {\norm{\eta_\a}\over{t_\a}}\right)^{m_\a}
\st0\le m_\a<n(\a),\ \a\in U\in\U\fs\ \right).$$

\pf The existence of the homomorphism follows from Proposition 1 and the
transfinite recursion theorem, but we include the details here
for completeness. We use the transfinite recursion theorem informally.

Let $\b\in\U$ and suppose that $(\f_\a\map{B_\a}{B^{(2)}})_{\a<\b}$ is a family
of homomorphisms satisfying
\item{(i)} for all $\a_1,\a_2<\b$, $\f_{\a_2}\rest{B_{\a_1}}=\f_{\a_1}$ if $\a_1\le\a_2$,
\item{(ii)} for all $\a<\b$ and $a\in A$, $\f_\a(a)=\th(a)$, and $\f_\a(\xi_\a)=\eta_\a$.

Let $\psi_\b\map{B_{<\b}[x]/(\b(x))}{B_\b}$ be the canonical, isometric isomorphism. Let
$$\f_{<\b}\;:\;B_{<\b}\to B^{(2)}\;;\;b\mapsto\cases{\th(b) &$b\in A$,\cr
\f_\a(b)& if $b\in B_\a, \a<\b$.\cr}$$
By hypothesis, $\f_{<\b}$ is well-defined.
By Proposition 1.2.1.2, there is a unique homomorphism $\ths\b\map{(B_{<\b})_\b}{B^{(2)}}$ such that
$\ths\b\circ\nu_{<\b}=\f_{<\b}$ and $\ths\b(\xb)=\eta_\b$. Let the canonical embedding
of $B_{<\b}$ into the Arens-Hoffman extension $(B_{<\b})_\b$ be denoted by $\nu_{<\b}$.
We then have the following commutative diagram
$$\matrix{ B_\b& \mapleft{\psi_\b} &(B_{<\b})_\b &\mapright{\th_\b} &B^{(2)}\cr
&&\mapup{\nu_{<\b}}&\mapne{\f_{<\b}}&\cr
&&B_{<\b}&&\cr}$$

Set $\f_\b=\th_\b\circ\psi_\b^{-1}$. Then

\item{(i)} for all $\a<\b$, $\f_\b\rest{B_\a}=\f_\a$ since for $b\in B_\a$ we have
$\f_\b(b)=\th_\b(b)=\f_{<\b}(b)=\f_\a(b)$;
\item{(ii)} for all $a\in A$, $\f_\b(a)=\th_\b(a)=\f_{<\b}(a)= a$.

Finally, $\f_\b(\xi_\b)=\f_\b(\psi_\b(\xb))=\th_\b(\xb)=\eta_\b$.

Hence the new family of homomorphisms, $(\f_\a)_{\a\le\b}$, is consistent with
conditions (i) and (ii) and by the transfinite
recursion theorem there exists a family $(\f_\a)_{\a\in\U}$ which is also
consistent with them. The required homomorphism is defined by setting
$\phi(b)=\f_\a(b)$ where $\a\in\U$ is such that $b\in B_\a$.

It is obvious that $\f$ is unique.

Finally we check the statement about continuity. In view of Lemma 2, it is
plainly sufficient that $\th$ is continuous and that
$$\left(\prod_{\a\in U}
\left(\ {\norm{\eta_\a}\over{t_\a}}\right)^{m_\a}
\st0\le m_\a<n(\a),\ \a\in U\in\U\fs\ \right)\qquad\hbox{is bounded.}$$

An elementary exercise in summability shows that this condition is
equivalent to the one given in the statement.\eop

\npar Purely
algebraic standard extensions are defined in $\ref{LinIE}$ and
the main content of Proposition 3
is a statement about these.

Notice that if the condition on the norm-parameters is satisfied then, by standard
facts about summability (see, for
example, $\ref{Kel}$, p.\ 78), only countably many $\aiu$ are such that
$\log^+\left(\norm{\eta_\a}/{t_\a}\right)$ is non-zero. Equivalently,
$\set{ \aiu\st\norm{\eta_\a}>t_\a}$ must be countable.

Notice too that although the definition of a standard extension, $B_\U$, appears
to depend on the choice of well-ordering of $\U$, Proposition 3 shows that, provided
the same norm parameters are used for all $\a\in\U$, then a standard extension
with respect to any other well-ordering of $\U$ is isometrically isomorphic to
$B_\U$.

\vskip 10pt\noindent{\bf Narmania Extensions}
\npar
From now on in Section 1.2 we assume that $A$ is complete.

The {\it Narmania extension}, $A_\U$, of $A$ by $\U$ is equal to the
Banach-algebra direct limit of $(A_S$ : $S$ is a finite subset of $\U)$
where each $A_S$ is isometrically isomorphic to $A$ extended finitely many times
by the Arens-Hoffman construction. We shall give the full details of this below.
See Appendix 2 for a discussion of
direct limits of Banach algebras.

Let $t_\a\ (\a\in\U)$
be a valid choice of Arens-Hoffman norm-parameters for the Arens-Hoffman extensions
$\Aa\ (\a\in\U)$. It is important to insist that distinct elements $\a,\b\in\U$ are associated
with distinct indeterminates $x_\a,x_\b$. Thus $S=\set{\a_1,\ldots,\a_N}\in\U\fs$ is an abbreviation for
$\set{\a_1(x_{\a_1}),\ldots,\a_m(x_{\a_m})}$.

It is proved carefully in $\ref{Nar}$ that for $q=\sum_m q_m x_{\a_1}^{m_1} \cdots x_{\a_N}^{m_N}
\in A[x_{\a_1},\ldots,x_{\a_N}]$, the algebra of polynomials in $N$
commuting indeterminates over $A$, ($m$ is a multiindex in $\No^N$
where $\No=\set0\cup\N$) then $(S)+q$ has a unique
representative whose degree in $x_{\a_j}$ is less than than $n(\a_j)$, the degree
of $\a_j(x_{\a_j})$ ($j=1,\ldots,N$). For convenience we shall
call such representatives {\it minimal}. Then if $q$ is the minimal
representative
of $(S)+q$,
$$\norm{(S)+q}:=\sum_m \norm{q_m} t_{\a_1} ^{m_1} \cdots t_{\a_N} ^{m_N}$$
defines an algebra norm on $A_S$ (see $\ref{Nar}$). The index set, $\U\fs$ is a directed set,
directed by $\subseteq$. The connecting homomorphisms $\nu_{S,T}$ (for $S\subseteq T\in\U\fs$)
are the natural maps;
they are isometries. Thus
$A_\U$ is the completion of the normed direct limit,
$D:=\bigcup_{S\in\U\fs}A_S\bigm/\sim$,
where $\sim$ is an equivalence relation given by
$(S)+q\sim (T)+r$ if and only if $q-r\in (S\cup T)$ for $S, T\in\U\fs$. Furthermore (see Appendix 2),
the canonical map, $\nu_S$, which sends an element of $A_S$ to its equivalence
class in $D$, is an isometry. Note that $A_\emptyset$ is defined to be $A$.

The embedding $A\to\AU$ is given by the canonical map $\nu_\emptyset \; \colon\;A_\emptyset\to \AU$.

The structure of $A_\U$ will be clarified in Section 1.3.

\vskip 10pt\noindent{\bf Cole Extensions}
\npar
The small step of generalising Cole's method of adjoining square roots
to adjoining roots of arbitrary monic polynomials over uniform algebras
was taken explicitly in $\ref{4DIS}$. It was doubtless known to Cole, Lindberg, Gorin, and others
working in the field.

\prop 1.2.2.4 ($\ref{Col},\ref{4DIS}$). Let $\U$ be a set of monic polynomials over the
uniform algebra $(A,X)$. There exists a uniform algebra
$(A^\U,X^\U)$ and a continuous, open surjection $\pi\map{X^\U}X$
such that
\item{(i)} the adjoint map $\ps\mapto{C(X)}{C(X^\U)}g{g\circ \pi}$ induces an
isometric, unital monomorphism $A\to A^\U$,
\item{(ii)} for every $\a\in\U$ there exists $p_\a\in A^\U$
such that $\ps(\a)(p_\a)=0$, and
\item{(iii)} the uniform algebra $A^\U$ is generated by $\ps(A)\cup\set{p_\a\st \a\in\U}$.

\pf We let $X^\U$ be the subset of $X\times\C^\U$ consisting of the elements
$(\k,\l)$ such that for all $\a\in\U$
$$f_0^{(\alpha)}(\kappa)+\cdots+
f_{n(\alpha)-1}^{(\alpha)}(\kappa)\lambda_\a^{n(\alpha)-1}+\lambda_\a^{n(\alpha)}=0$$
where $\a(x)=f_0^{(\a)}+\cdots+f_{n(\a)-1}^{(\a)} x^{n(\a)-1}+x^{n(\a)}\in\U$. We
check that $X^\U$ is a compact, Hausdorff space in the relative product-topology.

The space is Hausdorff because subspaces and products of Hausdorff spaces are
again Hausdorff.
To show $X^\U$ is compact we require a lemma. Suppose $c_0,\ldots,c_{m-1}\in\C, \lambda\in\C,
\ab\lambda\ge 1$, and $c_0+\cdots+c_{m-1}\lambda^{m-1}+\lambda^m=0$.
By assumption, $-\lambda^m = c_0+\cdots+c_{m-1}\lambda^{m-1}$.
Hence $\ab\lambda^m =\left\vert c_0+\cdots+c_{m-1}\lambda^{m-1}\right\vert$,
from which we see that $\ab\lambda^m \le \ab {c_0}+\cdots+\ab {c_{m-1}} \ab\lambda^{m-1}
	\le \ab {c_0}\ab\lambda^{m-1}+\cdots+\ab {c_{m-1}}\ab\lambda^{m-1}$ since
$\ab\lambda\ge1$.
Therefore $\ab\lambda \le \ab {c_0}+\cdots+\ab{c_{m-1}}$.
Hence for all
${(\kappa, \lambda)}\in X^\U \hbox{ and }\a\in \U,\ \ab{\lambda_\a}\le 1\hbox{ or }\ab{\lambda_\a}\le
\left\vert\sum\nolimits_{k=0}^{n(\a)-1} f_k^{(\a)}(\kappa)\right\vert\le
\sum\nolimits_{k=0}^{n(\a)-1}\left\Vert f_k^{(\a)}\right\Vert$, and so
$$X^\U\subseteq X\times\prod_{\a\in\U}B_\C\Bigl[0,{{\rm max}\left(1,\,
\norm{f_0^{(\alpha)}}+\cdots+\norm{f_{n(\alpha)-1}^{(\alpha)}}\right)}\Bigr]$$
which is compact by Tychonoff's theorem.
Now by definition, $X^\U$ is the intersection of a family of zero sets of continuous
functions and so closed. It follows that $X^\U$ is compact.

Therefore the following functions are continuous:
$$\eqalign{\pi&\mapto{X^\U}X{(\k,\l)}\k,\qquad\hbox{and}\cr
p_\a&\mapto{X^\U}\C{(\k,\l)}\l_\a\qquad(\a\in\U).\cr}$$
The extension $A^\U$ is defined to be the closed subalgebra of $C(X^\U)$ generated
by $\ps(A)\cup\set{p_\a\st \a\in\U}$.

It is not hard to check that $A^\U$ is a uniform algebra on $X^\U$ with the required
properties.\eop

\npar We shall call the algebra $A^\U$ constructed above
the {\it Cole extension} of $A$ by $\U$. The embedding is
given by the restriction of the adjoint, $\ps$, to $A$. We shall use the same symbol
for both maps; this does not cause problems in practice.

We mention now that $\UA$ is non-trivial if $A$ is. Details
in the case when every element of $\U$ is of the form $x^2-f$ for some
$f\in A$ are in $\ref{Col}$. We shall prove a more general result in Section 1.4.

\vskip 15pt
\noindent {\bf 1.2.3. Systems of Algebraic Extensions}\vskip 10pt

\noindent We have just shown how to adjoin to $A$ roots of an arbitrary set of monic polynomials,
$\U$, over $A$. Let the embedding for the
extension be denoted by $\th=\tss A\U\map A{E(A,\U)}$. Thus in $\na$, $\th$ is inclusion,
in $\Ba$, $\th$ is given by the natural map into a set of equivalence classes (see Appendix 2), and
in $\ua$, $\th$ is given by an adjoint of the continuous surjection $\pi\map{X^\U}X$, for
some compact, Hausdorff space, $X$.
For future reference, we mention here that an extension will be called {\it simple} if
$\U$ is a singleton.

By repeating the process of forming algebraic extensions using transfinite recursion, one can
generate normed algebras which are integrally closed extensions of $A$.
The process, or a variant of it, has been
used frequently in the literature; see, for example, $\ref{Col}$ and $\ref{Nar}$.

We now describe some of the detail informally.
Let $\ui$ be a non-zero ordinal number. Let $\u<\ui$ and suppose that we have chosen
algebras $(A_\t)_{\t<\u}$ and homomorphisms $(\tss\s\t\map{A_\s}{A_\t})_{\s\le\t<\u}$ in
the appropriate category such that
\item{(i)} $A_0=A$ and $\tss00=\id{A_0}$;
\item{(ii)} $\left((A_\t)_{\t<\u}, (\tss\s\t)_{\s\le\t<\u}\right)$ is a direct system;
\item{(iii)} if $0<\t<\u$ and $\t=\s+1$ for some $\s<\t$ then there exists a set of
monic polynomials, $\U_\s$, over $A_\s$ such that $A_\t=E(A_\s,\U_\s)$ and
$\tss\s\t=\tss{A_\s}{\U_\s}$;
\item{(iv)} if $0<\t<\u$ and $\t$ is a limit ordinal then
$$\left((A_\t), (\tss\s\t)_{\s\le\t}\right)
=\dl{}\left((A_\s)_{\s<\t}, (\tss\r\s)_{\r\le\s<\r}\right).$$

There are two cases. First suppose that $\u=\t+1$ for some $\t<\u$. Then we can choose a set
of monic polynomials, $\U_\t$, over $A_\t$ and define
$A_\u:=E(A_\t,\U_\t)$ and
$$\tss\s\u:=\cases{\id{A_\u}&if $\s=\u$;\cr
	\tss{A_\t}{\U_\t}&if $\s=\t$;\cr
	\tss\t\u\circ\tss\s\t&if $\s<\t$.\cr} $$
It is not hard to check that $\left((A_\t)_{\t\le\u}, (\tss\s\t)_{\s\le\t\le\u}\right)$ is a direct system.
If $\r_1\le\r_2\le\r_3\le\u$ then for $\r_3<\u$, $\tss{\r_2}{\r_3}\circ\tss{\r_1}{\r_2}=\tss{\r_1}{\r_3}$
by hypothesis. Now suppose $\r_3=\u$: there are two more cases. First
suppose that $\r_2=\u$. Then $\tss{\r_2}{\r_3}=\id{A_\u}$ so
$\tss{\r_2}{\r_3}\circ\tss{\r_1}{\r_2}=\tss{\r_1}{\u}=\tss{\r_1}{\r_3}$. Otherwise $\r_2<\u$, in
which case we have
$$\eqalign{\tss{\r_2}{\r_3}\circ\tss{\r_1}{\r_2}&=\tss{\r_2}{\u}\circ\tss{\r_1}{\r_2}\cr
&=\tss\t\u\circ\tss{\r_2}{\t}\circ\tss{\r_1}{\r_2}\qquad\hbox{ (by definition of }\tss{\r_2}{\u})\cr
&=\tss\t\u\circ\tss{\r_1}{\t}\cr
&=\tss{\r_1}{\r_3}\qquad\qquad\hbox{ (by definition of }\tss{\r_1}{\u}).\cr}$$
So in all cases $\tss{\r_2}{\r_3}\circ\tss{\r_1}{\r_2}=\tss{\r_1}{\r_3}$. It is clear that
for all $\r\le\u$, $\tss\r\r=\id{A_\r}$.

Therefore (i)-(iv) are true with `$<\u$' replaced by `$\le\u$'.

Secondly suppose that $\u$ is a limit ordinal greater than $0$. We set
$$(A_\u,(\tss\t\u)_{\t\le\u}):=\dl{}\left((A_\t)_{\t<\u}, (\tss\s\t)_{\s\le\t<\u}\right),$$
with the direct limit taken in the appropriate category.
In this way we obtain a new direct system, $\left((A_\t)_{\t\le\u}, (\tss\s\t)_{\s\le\t\le\u}\right)$,
consistent with conditions (i)-(iv) above with `$<\u$' replaced by `$\le\u$'.
By an informal application of the transfinite recursion theorem there is a direct system
$\left((A_\t)_{\t\le\ui}, (\tss\s\t)_{\s\le\t\le\ui}\right)$ in which $A=A_0$ and for all $0<\t<\ui$ we have
\item{(i)} if $\t=\s+1$ for some $\s<\t$ then there exists a set of
monic polynomials, $\U_\s$, over $A_\s$ such that $A_\t=E(A_\s,\U_\s)$ and
$\tss\s\t=\tss{A_\s}{\U_\s}$;
\item{(ii)} if $\t$ is a limit ordinal then
$$\left((A_\t), (\tss\s\t)_{\s\le\t}\right)
=\dl{}\left((A_\s)_{\s<\t}, (\tss\r\s)_{\r\le\s<\t}\right).$$
\noindent
We are now justified in making the following definition.

\dfn 1.2.3.1. Let $\ui$ be an ordinal number.
Let ${\cal A}=\left((A_\t)_{\t\le\ui}, (\tss\s\t)_{\s\le\t\le\ui}\right)$ be a direct system
in $\na$, $\Ba$, or $\ua$. Then ${\cal A}$ is a
{\it system of algebraic extensions of $A$} provided the following three conditions hold:
\item{(i)} $A_0=A$;
\item{(ii)} if $0<\t\le\ui$ and $\t=\s+1$ for some $\s<\t$ then there exists a set of
monic polynomials, $\U_\s$, over $A_\s$ such that $A_\t=E(A_\s,\U_\s)$ and
$\tss\s\t=\tss{A_\s}{\U_\s}$;
\item{(iii)} if $0<\t\le\ui$ and $\t$ is a limit ordinal then
$$\left((A_\t), (\tss\s\t)_{\s\le\t}\right)
=\dl{}\left((A_\s)_{\s<\t}, (\tss\r\s)_{\r\le\s<\t}\right).$$

\noindent
We shall speak of, for example, a `system of Cole extensions', to mean
a system of algebraic extensions in $\ua$.

\vskip 15pt
\noindent {\bf 1.2.4. Integrally Closed Extensions}\vskip 10pt

\noindent It is well-known from the literature (see, for example,
$\ref{Col}$, $\ref{LinIE}$, and $\ref{Nar}$) that systems of algebraic extensions
can produce integrally closed algebras. We give details of the method for the
convenience of the reader.

The case of normed algebras is discussed extensively in $\ref{LinIE}$ so we
shall assume that $A$ is an algebra in $\Ba$ or $\ua$. Let $\oi$ denote the
first uncountable ordinal and
$\left((A_\t)_{\t\le\oi}, (\tss\s\t)_{\s\le\t\le\oi}\right)$ be a system of
algebraic extensions in the same category. Also suppose that for every $\t<\oi$, $\U_\t$
is taken to be the set of all monic polynomials over $A_\t$. The method relies
on the following well-known lemma. This does not appear in standard textbooks so
we state and prove it here.

\lem 1.2.4.1. Let $(X, d)$ be a metric space and $(F_\s)_{0\le\s<\omega_1}$ be a
family of closed sets in $X$ such that for all $\s\le\t<\omega_1$, $F_\s\subseteq F_\t$.
Then $F:=\bigcup_{\s<\omega_1}F_\s$ is closed.
\pf Let $(x_n)_{n\in\N}\subseteq F$ and $x_n\to x\in X\;(n\to+\infty)$.
Now for all $n\in\N$ there exists $\s_n<\omega_1$ such that $x_n\in F_{\s_n}$.
There also exists $\s<\omega_1$ such that for all $n\in\N$, $\s_n\le\s$ for if this were not true then
for all $\s<\omega_1$ there exists $n\in\N$ such that $\s<\s_n$.
Hence $[0, \omega_1)\subseteq\bigcup_{n\in\N}[0, \s_n)$.
But by the definition of $\omega_1$, $[0, \s_n)$ is countable for all $n\in\N$.
Therefore $\omega_1=[0, \omega_1)$ is countable since it is contained in a
countable union of countable sets. This is a contradiction.

Therefore $(x_n)\subseteq F_\s$ and since $F_\s$ is closed $x\in F_\s\subseteq F$.\eop

\npar Since $\oi$ is a limit ordinal, the algebra $D=\bigcup_{\t<\oi}\tss\t\oi(A_\t)$
is dense in $A_\oi$. But (see Section 1.2.2) $\tss\t\oi$ is an isometry
so $\tss\t\oi(A_\t)$ is closed in $A_\oi$ for each
$\t<\oi$. By the lemma, $D$ is closed
in $A_\oi$ and therefore equals $A_\oi$.
Let $\a(x)$ be a monic polynomial over $A_\oi$. Then $\a(x)$ is a polynomial
over $\tss\t\oi(A_\t)$ for some $\t<\oi$ and so there exist $a_0,\ldots,a_{n-1}\in A_\t$
such that
$$\a(x)=\sum_{k=0}^{n-1}\tss\t\oi(a_k)x^k+x^n.$$
By construction there exists $\xi_0\in A_{\t+1}$ such that
$$\sum_{k=0}^{n-1}\tss\t{\t+1}(a_k)\xi_0^k+\xi_0^n=0.$$
Set $\xi=\tss{\t+1}\oi(\xi_0)$. Then
$$\eqalign{\a(\xi)&=\sum_{k=0}^{n-1}\tss\t\oi(a_k)\tss{\t+1}\oi(\xi_0)^k+\tss{\t+1}\oi(\xi_0)^n\cr
&=\sum_{k=0}^{n-1}\tss{\t+1}\oi(\tss{\t}{\t+1}(a_k))\tss{\t+1}\oi(\xi_0^k)+\tss{\t+1}\oi(\xi_0^n)\cr
&=\tss{\t+1}\oi\left(\sum_{k=0}^{n-1}\tss\t{\t+1}(a_k)\xi_0^k+\xi_0^n\right)\cr
&=\tss{\t+1}\oi(0)=0.\cr}$$
Therefore $A_\oi$ is integrally closed.


\def\nb{{n(\b)-1}}
\def\Os{\O_\s}
\def\Ot{\O_\t}
\def\hai{\widehat{a_i}}

\vskip 15pt
\noindent {\medtenbf 1.3. Comparison of the Extensions}

\npar It is useful to be able to describe Narmania and Cole extensions in terms of standard extensions.
In order to do this, we identify the spaces of closed,
maximal ideals of the extensions. It is also appropriate to include a result about naturalness
and the uniform algebras in a system of Cole extensions.

\vskip 15pt
\noindent {\bf 1.3.1. The Spaces of Closed, Maximal Ideals of the Three Extensions}\vskip 10pt

\noindent It is not hard to use Arens' and Hoffman's description ($\ref{Are}$) of the
space of closed, maximal ideals of an Arens-Hoffman extension of $A$ to describe $\O(B_\U)$ where
$B_\U$ is a standard extension of $A$ with respect to a set, $\U$,
of monic polynomials over $A$. The next result was probably known to
Lindberg, but does not appear in the literature.

Let $\OU$ denote the following set
$$\set{(\o,\l)\in\O\times\C^\U\st\hbox{for all }\a\in\U\ \o(\a)(\l_\a)=0}.$$
It is clearly compact and Hausdorff (compare Proposition 1.2.2.4)
when given the relative product-topology.

\prop 1.3.1.1. The space $\O(B_\U)$ is homeomorphic to $\OU$.
\pf It is easy to check that the map
$$\O(B_\U)\to\OU\;;\;\o\mapsto(\o\rest A,(\o(\xi_\a))_\aiu)$$
is a homeomorphism. (To show surjectivity requires a routine application of
the transfinite recursion theorem. The full details of this, which are similar to
those in Proposition 1.2.2.3, are given in Theorem 5.6 of $\ref{4DIS}$.)\eop

\npar Note that if $A$ is a natural uniform algebra, then, taking $X=\O$ in
Proposition 1.2.2.4, $\OU=\O^\U$. It is mentioned
in $\ref{Col}$ (for the case of square roots) that
$\UA$ is natural if $A$ is. We shall prove a similar result for systems of Cole extensions.
First we need to establish a general result which is well-known in the literature. We
give details for the convenience of the reader.
The result is stated in $\ref{Kar}$, and also follows from standard results in
$\ref{Lei}$ (p. 212-213).

\prop 1.3.1.2. Let $\u$ be a non-zero, limit ordinal and
$$\left( (\O_\t)_{\t<\u},(\Ps\s\t)_{\s\le\t<\u}\right)$$ be an inverse system of
compact, Hausdorff spaces and continuous surjections. Suppose that
$$\left( (A_\t,\O_\t)_{\t<\u},(\pss\s\t\colon f\mapsto f\circ\Ps\s\t)_{\s\le\t<\u}\right)$$
is a direct system of natural uniform algebras.
Then the direct limit of this system, $A_\u$, is natural.
\pf
Let $(\O_\u,(\Ps\t\u)_{\t\le\u})$
be the inverse limit of the system of spaces and let
$$((A_\u,\O_\u),(\pss\t\u)_{\t\le\u})$$ be the direct limit of the system of uniform
algebras. From the discussion in Appendix 2, we see that $A_\u$ is the
uniform closure in $C(\O_\u)$ of $D=\bigcup_{\t<\u}\pss\t\u(A_\t)$.

Let $\ph\in\Om{A_\u}$. Then for each $\t<\u$, $\ph\circ\pss\t\u\in\Om{A_\t}$ so
there exists a unique $\o^\t\in\O_\t$ such that $\ph\circ\pss\t\u=\eh{\o^\t}$
where $\e$ is as in Section 1.1.3. If $\s\le\t<\u$ then
$$\eh{\o^\s}=\ph\circ\pss\s\u=\ph\circ\pss\t\u\circ\pss\s\t=\eh{\o^\t}\circ\pss\s\t$$
so for all $f\in A_\s$ we have
$$\eqalign{ f(\o^\s) &= \eh{\o^\t}(\pss\s\t(f))\cr
&=(f\circ\Ps\s\t)(\o^\t)\cr
&=f(\Ps\s\t(\o^\t)).\cr}$$
Therefore $\o^\s=\Ps\s\t(\o^\t)$ and so $\o^\u:=(\o^\t)_{\t<\u}\in\O_\u$. Moreover
$\ph=\eh{\o^\u}$ since both maps agree on $D$. Thus $\Om{A_\u}=\O_\u$ as required.\eop

\noindent We can now give the promised result on Cole extensions.

\prop 1.3.1.3. Suppose that
$\left( (A_\t)_{\t\le\u},(\pss\s\t\colon f\mapsto f\circ\Ps\s\t)_{\s\le\t\le\u}\right)$
is a system of Cole extensions of the natural uniform algebra, $(A,\O)$.
Then for all $\t\le\u$, $A_\t$ is
natural.

\pf Let $A_\t$ be a uniform algebra on $X_\t$ for all $\t\le\u$. Let
$$\J=\set{\t\le\u\st A_\t\hbox{ is natural }},$$ $\t\le\u$, and suppose that
$[0,\t)\subseteq\J$. By assumption $0\in\J$. Suppose now that $\t>0$. In
the case where $\t$ is a limit ordinal, it follows immediately from the proposition
above that $\t\in\J$. If $\t=\s+1$
for some $\s<\t$ then there exists a set, $\U_\s$, of monic polynomials over $A_\s$
such that $A_\t=A_\s^{\U_\s}$. Let $\o\in\Om{A_\t}$. Since $A_\s$ is natural
and $\o\circ\pss\s\t$ is a character on $A_\s$, there exists $\k\in X_\s$ such that
$\o\circ\pss\s\t=\eh\k$. Now let $a_0,\ldots,a_{n-1}\in A_\s$ be such that $\a(x)
=a_0+\cdots+a_{n-1}x^{n-1}+x^n\in\U_\s$.
By construction,
$$\pss\s\t(a_0)+\cdots+\pss\s\t(a_{n-1})p_\a^{n-1}+p_\a^n=0.$$
Since $\o$ is a homomorphism we have
$$\eh\k(a_0)+\cdots+\eh\k(a_{n-1})\l_\a^{n-1}+\l_\a^n=0$$
where $\l_\a:=\o(p_\a)$. Thus $y=(\k,(\l_\a)_{\a\in\U_\s})\in X_\t$.
Since the polynomials, $P$, in $\pss\s\t(A)\cup\set{p_\a\st\a\in\U_\s}$ are dense
in $A_\t$ and $\o\rest P=\eh y\rest P$ we have $\o=\eh y$ by continuity. Since $\o$
was arbitrary, $A_\t$ is natural and so $\t\in\J$.

By the transfinite induction theorem, $\J=[0,\u]$ as required.\eop

\vskip 15pt
\noindent {\bf 1.3.2. The Relationship between the Three Extensions}\vskip 10pt

\noindent
We can now show how the types of extensions we have been considering are related. Many
of the ideas behind the next proposition are due to Narmania. We take the step
of linking them to Cole and standard extensions.

\prop 1.3.2.1. Let $A$ be a commutative, unital Banach algebra and $\U$ a set of
monic polynomials over $A$. Then, provided the norm-parameters
for the standard extension, $B_\U$, and Narmania extension, $\AU$, are equal, we have
up to isometric isomorphism that $A_\U=\c{B_\U}$.
If $A$ is a uniform algebra then we have
$$A^\U=\overline{(A_\U)\htt}=\overline{(B_\U)\htt},$$
where the closures are taken with respect to the supremum norm and $\UA$ denotes the Cole extension.

\pf By Lemma A.1.5, if $B$ is a normed algebra then the
homeomorphism $\Om{\c B}\to\Om B$ induces an isometric isomorphism
$\snc B\to\snc{(\c B)}$. It is therefore sufficient to prove that $A_\U=
\c{B_\U}$ and that $A^\U=\snc{(B_\U)}$. The last equality follows very quickly from
Lemma 1.2.2.2 and the simplicity
of the definition of $A^\U$. (We may need to pass to $\hA$ if necessary first, in
order to obtain an isometrically isomorphic, {\it natural} uniform algebra.)
We shall therefore only prove the first identification.
Although what follows is routine,
we hope that it will help
to clarify the details of standard and Narmania extensions.

As before let $t_\a\ (\a\in\U)$
be a valid choice of Arens-Hoffman norm-parameters for the respective extensions $\Aa$.
We shall show that there is an isometric isomorphism between $B_\U$ and the
normed-algebra direct limit
$D:=\bigcup_{S\in\U\fs}\nu_S(A_S)$
(when defined by these parameters); the result then follows from the uniqueness
of completions (see Appendix 1).

For each $\a\in\U$ let $y_\a$ be the equivalence class $[(\set{\a(x_\a)})+x_\a]\in D$.
Since $y_\a$ is a root of $\nu_\emptyset(\a)(x)$ in $D$ there exists, by
the universal property of standard extensions, a (unique) homomorphism
$\f\map{B_\U}D$ such that $\f\vert_A=\nu_\emptyset$ and for all $\a\in\U$, $\f(\xi_\a)
=y_\a$. Here, $\xi_\a$ is the standard root of $\a(x)$,
the element of $B_\U$ associated with $\xb$ by
the isometric isomorphism $\psi_\a\;\colon\; B_{<\a}[x]/(\a(x))\to B_\a$ in the notation
of Theorem 1.2.2.1. Thus $\norm{\xi_\a}=t_\a$. Note that
$\psi_\a$ satisfies
$\psi_\a(a)=a$ for all $a\in B_{<\a}$.

It is clear that $\f$ is surjective; we now use the transfinite induction theorem, as
is customary for proving results about standard extensions, to show that $\f$ is
isometric.

Let $\J=\set{\b\in\U\st\f\vert_{B_\b}\hbox{ is isometric}}$.
It should be
clear to the reader that $\a_0\in\J$. Let $\b\in\U$ and suppose that
$[\ao,\b)\subseteq\J$. Let $b\in B_\b$.
Then, writing $n(\b)$ for the degree of $\b(x)$, there
exist unique $b_0,\ldots,b_{n(\b)-1}\in B_{<\b}$ such that
$b=\sum_{j=0}^\nb b_j\xi_\b^j$. Let $\a<\b$ with $b_0,\ldots,b_{n(\b)-1}\in B_\a$.

A routine exercise in the transfinite induction
theorem shows that for all $\g\in\U$,
$$\f(B_\g)\subseteq\bigcup_{S\in[0,\g]\fs}\nu_S(A_S).$$

We give the details here for the benefit of the reader.
Let $$\J\pri=\set{\g\in\U\st \f(B_\g)\subseteq
\cup_{S\in[0,\g]\fs} \nu_{S}(A_S)}$$ and suppose that $[\ao,\g)\subseteq\J\pri$.
Let $p$ be the degree of $\g(x)$ and $c\in B_\g$.
There exist unique $c_0,\ldots,c_{p-1}\in B_{<\g}$
such that $c=\sum_{j=0}^{p-1} c_j\xi_\g^j$. There is also some $\a<\g$ such that
$c_0,\ldots,c_{p-1}\in B_\a$ (unless $\g=\ao$ in which case $c_0,\ldots,c_{p-1}
\in A$). The inductive hypothesis implies that there exists a finite subset
$S\subseteq[\ao,\a]$ such that $\f(c_j)\in\nu_{S}(A_S)$ ($j=0,\ldots,p-1$).
Therefore $\f(c)\in\nu_{S\cup\set\g}(A_{S\cup\set\g })$ and so $\g\in\J\pri$.

By the transfinite induction theorem, $\J\pri=\U$

Since the algebras $\nu_{S}(A_S)$ are directed there exists $S\in[\ao,\a]\fs$ such that
$\f(b_j)\in\nu_{S}(A_S)\ j=0,\ldots,\nb$. Let $S=\set{\a_1,\ldots, \a_m}$
and $q_0,\ldots,q_\nb$$\in A[x_{\a_1},\ldots,x_{\a_m}]$
be the minimal representatives such that
$$\f(b_j)=\nu_S((S)+q_j)\qquad\qquad (j=0,\ldots,\nb).$$
In fact, let
$$q_j=\sum_{l^j\in L}q_{j, l^j}x_{\a_1}^{l_1^j}\cdots  x_{\a_m}^{l_m^j}$$
with $L:=\times_{k=1}^m\set{0,\ldots,n(\a_k)-1}$, and
$q_{j,l^j}\in A$ for each $l^j\in L$ and $j=0,\ldots,\nb$.
Set $T=S\;\udot\set\b=\set{\a_1,\ldots,\a_{m+1}}$ where $\a_{m+1}:=\b$. By considering the
degrees, we see that
$$\eqalign{q&=\sum_{j=0}^{n(\a_{m+1})-1} q_jx_{\a_{m+1}}^j\cr
&=\sum_{j=0}^{n(\a_{m+1})-1}\sum_{l^j\in L} q_{j,l^j}
x_{\a_{1}}^{l_1^j}\cdots x_{\a_{m}}^{l_m^j}x_{\a_{m+1}}^j\cr}$$
is the minimal representative for $(T)+q$ in $A_T$. Hence
$$\eqalign{\norm{\nu_T((T)+q)}&=\norm{(T)+q}_{A_T}\cr
&=\sum_{j=0}^{n(\a_{m+1})-1}\sum_{l^j\in L} \norm{q_{j,l^j}}
t_{\a_{1}}^{l_1^j}\cdots t_{\a_{m}}^{l_m^j}t_{\a_{m+1}}^j\cr
&=\sum_{j=0}^{n(\a_{m+1})-1} \norm{(S)+q_j}_{A_S}t_{\a_{m+1}}^j\cr
&=\sum_{j=0}^{\nb} \norm{\nu_S((S)+q_j)}t_{\b}^j\cr
&=\sum_{j=0}^{\nb} \norm{\f(b_j)}t_{\b}^j\qquad\hbox{ (from above)}\cr
&=\sum_{j=0}^{\nb} \norm{b_j}t_{\b}^j\qquad\hbox{ (by hypothesis).}\cr}$$
By construction of the standard norm, the quantity in the line above is equal to $\norm b$.

Now since $\f$ and the natural maps are homomorphisms,
$$\eqalign{\f(b)&= \sum_{j=0}^{\nb} \f(b_j)y_\b^j\cr
&= \sum_{j=0}^{\nb} \nu_S((S)+q_j)\nu_{\set\b}((\set\b)+x_\b^j)\cr
&= \sum_{j=0}^{\nb} \nu_T((T)+q_jx_\b^j)\cr
&= \nu_T((T)+q),\cr}$$
and so $\norm{\f(b)}=\norm b$ as required.

By the transfinite induction theorem, $\J=\U$.\eop

\npar We can also summarise this result with the commutative diagram (in
which the objects are subject to the same conditions as in Proposition 1.3.2.1
and the horizontal arrows have dense images):
$$\matrix{ B_\U &\mapright{\phi} &\AU &\mapright\G&\UA\cr
 &\mapnw{\inc} &\mapup\nu &\mapne{\ps}&\cr
 &&A&&\cr}$$
By results in Appendix 1, this yields the commutative diagram
(in which $\approx$ indicates a homeomorphism):
$$\matrix{ \Om{B_\U} &\mapleft{\approx} &\O(\AU) &\mapleft\approx&\O(\UA)=\OU\cr
 &\mapse{\inc^*} &\mapdown{\nu^*} &\mapsw{\pi^{**}}&\cr
 &&\O&&\cr}$$
Note that $\pi^{**}$ is identifiable with $\pi$. We shall, as has become traditional
in the literature on these extensions, identify all the spaces in the top row
with $\OU$ and denote all the downward maps by $\pi$.

\npar The following result will often be tacitly used. It is essentially due to
Lindberg, but we use the comments above to apply
it to Narmania and Cole extensions.

\prop 1.3.2.2. Let $((A_\s),(\tss\s\t)_{\s\le\t}\st\s,\t\le\u)$ be a system of
algebraic extensions for some ordinal $\u>0$ and let the adjoint maps between
the associated spaces of closed, maximal ideals be denoted $\Ps\s\t\mapto{\O_\t}{\O_\s}\o{\o\circ\tss\s\t}$
($\s\le\t\le\u$) where $\O_\t:=\Om{A_\t}\ (\t\le\u)$. Then for all $\s\le\t\le\u$,
$\Ps\s\t$ is a continuous, open surjection.
\pf The continuity and surjectivity of the maps $\Ps\s\t$ are well-known
in the literature on algebraic extensions and follow
easily from standard results on inverse limits of compact, Hausdorff spaces; see for
example $\ref{Lei}$. However we include a proof of these properties for
the reader's convenience.

By the comments above we can assume that we are dealing with a system of
standard extensions. As we have noted in Appendix 2, the direct limits can be realised as
unions in this case and so the maps $\Ps\s\t$ are just restrictions for all
$\s\le\t\le\u$.

Let $\J=\set{\t\le\u\st\hbox{ for all }\s<\t,\ \Ps\s\t\hbox{ is a continuous, open
surjection }}$, $\t\le\u$, and suppose $[0,\t)\subseteq\J$.

Obviously $0\in\J$.

If $\t=\s+1$ for some $\s<\t$ then $\Ps\s\t$ is a continuous, open surjection by
Theorem 3.5 of $\ref{LinIE}$. So $\t\in\J$.

Otherwise, $\t$ is a non-zero, limit ordinal. By assumption, the maps $\Ps\r\s$
are continuous, open surjections for $0\le\r\le\s<\t$.
We can identify the space $\Ot$ as the inverse limit
$$\O_\t=\il{} \left((\Os),(\Ps\r\s)\st\r\le\s<\t\right)\eqno (*)$$
because, as we shall show, there is a homeomorphism from $\O_\t$ onto
the realisation (see Appendix 2), $X_\t$, of the inverse limit ($*$) which is given by
$$X_\t:=\set{\k\in\times_{\s<\t}\Os\st\hbox{ for all }\r\le \s<\t,\quad \k_\r=\Ps \r\s(\k_\s)}.$$

We now digress to verify that the map
$$R\mapto{\O_\t}{X_\t}\o{(\o\rest{A_\s}})_{\s<\t}$$
is a homeomorphism.

The map $R$ is injective since $A_\t=\cup_{\s<\t}A_\s$ and surjective since
each element in $X_\t$ clearly defines a norm-one, character on $A_\t$.
Now let $(\o_t)$ be a net in $\O_\t$ which converges to $\o\in \O_\t$.
Then for all $\s<\t$ and $a\in A_\s$, $\o_t\rest{A_\s}(a)=\o_t(a)\to\o(a)=\o\rest{A_\s}(a)$.
So for all $\s<\t$, $\o_t\rest{A_\s}\to\o\rest{A_\s}$ , by definition of the
topology on $\O_\s$ (see Appendix 1). So by definition of the product topology,
$R(\o_t)\to R(\o)$. This shows that $R$ is continuous. Since $\O_\t$ is
compact and $X_\t$ is Hausdorff, it follows (by Proposition 1.6.8 in $\ref{Ped}$,
for example) that $R$ is a homemorphism, as required.

Under the identification ($*$), the canonical maps $\Ps\s\t$ ($\s<\t$)
are the coordinate maps $\Ot\to\Os\;;\;\k\mapsto\k_\s$.

Let $\r<\t$.
We must show that $\Ps\r\t$ is a continuous, open surjection. Continuity and surjectivity
follow from standard properties on inverse limits given in
Lemma 9 on p. 210 of $\ref{Lei}$. It remains to prove that
$\Ps\s\t$ is open. For this we reproduce the argument in the proof of Theorem 4.1 of
$\ref{LinIE}$.

Recall that $\Ot$ has the weak topology induced by the functions $\widehat{A_\t}$
so a basic open neighbourhood in $\Ot$ has the form
$$\eqalign{V_\t&=V(\k;a_1,\ldots,a_n;\e)\cr
&:=\set{\k\pri\in\Ot\st\ab{\hai(\k\pri)-\hai(\k)}<\e\qquad i=1,\ldots,n},\cr}$$
for some $\k\in\Ot$, $\e>0$, and $a_1,\ldots,a_n\in \bigcup_{\s<\t}A_\s=A_\t$.
Hence there exists $\s<\t$ such that $a_1,\ldots,a_n\in A_\s$. We show that
$\Ps\rho\t(V_\t)$ is open.

Consider $V_\s:=\Ps\s\t(V_\t)$. The element $y\in\Os$ belongs to $V_\s$ if and only if
there exists $\k\pri\in\Ot$ such that $y=\Ps\s\t(\k\pri)$ and
$$\ab{\k\pri(a_i)-\k(a_i)}=\ab{\Ps\s\t(\k\pri)(a_i)-\Ps\s\t(\k)(a_i)}<\e\qquad(i=1,\ldots,n).$$
Thus $V_\s=V(\Ps\s\t(\k);a_1,\ldots,a_n;\e)$.

If $\r\le\s$ then $\Ps\r\t(V_\t)=(\Ps\r\s\circ\Ps\s\t)(V_\t)$ is open by the inductive hypothesis.
Otherwise, $\s<\r$. Set $V_\r=\Ps\r\t(V_\t)$. Then $a_1,\ldots,a_n\in A_\s\subseteq A_\r$,
and similarly to above,
$V_\r=V(\Ps\r\t(\k);a_1,\ldots,a_n;\e)$
is an open neighbourhood of $\Ps\r\t(\k)$ and so $\Ps\r\t$ is open.
So $\t\in\J$ and by the transfinite induction theorem, $\J=\U$ as required.\eop


\def\ts{\th_\#}
\def\xma#1{\xi_{\a_{#1}}^{m_{#1}}}
\def\xmaj{\xma j}
\def\tma#1{t_{\a_{#1}}^{m_{#1}}}

\def\smaj#1{s_{\a_{#1}}^{m_{#1}}}
\def\saj#1{s_{\a}^{#1}}

\def\naj{n(\a_j)}
\def\nn{{n(\a_1)\cdots n(\a_N)}}

\def\na{n(\a)}
\def\Tss#1#2{T_{#1,#2}}
\def\tss#1#2{\th_{#1,#2}}
\def\Ss#1#2{S_{#1,#2}}
\def\dist#1#2{{\rm dist}\,(#1,#2)}
\def\CC{{\cal C}}
\def\AA{{\cal A}}
\def\Ts#1#2{T_{#1,#2}}
\def\Tes#1#2{T^e_{#1,#2}}

\vskip 15pt
\noindent {\medtenbf 1.4. Maps back from the Extensions}

\vskip 15pt
\noindent {\bf 1.4.1. Introduction}\vskip 10pt

\noindent Let $\th\map A{E(A,\U)}$ denote the embedding map of
one of the algebraic extensions introduced
in Section 1.2. We show that $\th$ has a special left-inverse, $T$, in the category of
normed spaces. The material in this section is due to Narmania ($\ref{Nar}$). However
some parts of the proofs are new.

Apart from having other special properties, $T$
will be a unital contraction and {\it $A$-linear} in the sense that
$$T(\th(a)b)=aT(b)\qquad\qquad (a\in A,\ b\in E(A,\U)).$$
(Recall that a linear map is called a {\it contraction} if its norm is
at most one.)

Applications of the operator $T$ will be given in the next section, Section 2.2, and Section 3.3.

\vskip 15pt
\noindent {\bf 1.4.2. The Construction}\vskip 10pt

\noindent For $\a(x)=a_0+\cdots+a_{n-1}x^{n-1}+x^n\in\U$
and $j\in\No$ let $\saj j$ denote the $j$th Newton-sum of $\a(x)$.
Newton-sums can be defined inductively as follows (see $\ref{Nar}$ or p. 140 of $\ref{JacI}$):
$$\cases{\saj 0=n;&\cr
\saj k+a_{n-1}\saj{k-1}+\cdots+a_{n-k+1}\saj1+ka_{n-k}=0&$\qquad(1\le k\le n);$\cr
\saj{n+k}+a_{n-1}\saj{n+k-1}+\cdots+a_0\saj k=0 &$\qquad (k\in\N)$.\cr }$$
In particular, $\saj j$ is a polynomial
in the coefficients of $\a(x)$. It has also been widely noted in the literature (see for example
$\ref{Nar}$ and $\ref{Kar}$) that for each $\o\in\O$,
$$\o(\saj j)=\sum_{k=1}^{n(\a)} \l_k^j$$
where $\l_1,\ldots,\l_{n(\a)}$ are the roots of $\o(\a)(x)$ listed according to their
multiplicities. (Recall that the homomorphism $\o$ can be extended to $\o_\#\map{A[x]}{\C[x]}$
as in Section 1.1.1.)

We shall use these facts to prove the following theorem of Narmania. Notation is as in
Section 1.2.3; $\U$ denotes a set of monic polynomials over $A$. We use
$\ab\cdot$ to denote the spectral radius (see Appendix 1).

\thm 1.4.2.1. Let $E(A,\U)$ be a standard, Narmania, or Cole extension of $A$. If $A$
is not a uniform algebra, suppose further that the polynomials and norm-parameters
satisfy
$$\hbox{for all }\aiu,\ j\in\set{0,\ldots,n(\a)-1}\qquad\norm{\saj j}\le t_\a^j\naj.\eqno(*)$$
Then there exists a unital, $A$-linear contraction $T\map{E(A,\U)}A$ such that
$T\circ\th=\id A$.

\pf We define the operator $T_\U\map{B_\U}A$ for the standard extension
and show that it is an $A$-linear left-inverse for $\th$. We then show that ($*$)
implies that $T_\U$ has norm one. From this and Proposition 1.3.2.1 we see that
$T_\U$ has a unique, norm-one extension to the Narmania extension, $A_\U$.
The extension has the required properties by continuity.

To obtain the operator $\UA\to A$ for the Cole extension, where applicable, it
is enough to verify that the map $T_\U$ satisfies
$$\ab{T_\U(b)}\le\ab b\qquad (b\in B_\U)$$
and so induces an $\hA$-linear map
$$\widehat{T_\U}\mapto{\widehat{B_\U}}\hA{\hat b}{T_\U(b)\htt}$$
of norm one. Again, Proposition 1.3.2.1 implies that $\widehat{T_\U}$ extends
by continuity to a map $T^\U\map\UA A$. It is trivial to verify that $T^\U$ has
the required properties.

{\it Construction.} Please refer to Section 1.2.2 for the notation associated
with standard extensions. We use the transfinite recursion theorem, as in Proposition
1.2.2.3, to construct $T_\U$.

Let $\g\in\U$ and assume that $(T_\b\map{B_\b}A)_{\b<\g}$ is a family
of linear maps such that the following four conditions hold:
\item{(i)} for all $\b<\g$, $T_\b\rest A=\id A$;
\item{(ii)} for all $\a\le\b<\g$, $T_\b\rest {B_\a}=T_\a$;
\item{(iii)} for all $\b<\g$, $a\in A$, and $b\in B_\b$, $T_\b(ab)=aT_\b(b)$;
\item{(iv)} for all $\b<\g$, $\norm{T_\b}\le 1$.

Recall that $B_{<\ao}=A$ and that if $\g>\ao$ then $B_{<\g}=\bigcup_{\b<\g}B_\b$.
If $\g=\ao$, the minimum element of $\U$, define
$$T_{<\g}\mapto{B_{<\g}}Aaa.$$
Otherwise, if $\g>\ao$, let
$$T_{<\g}\mapto{B_{<\g}}Ab{T_\b(b)}\qquad\qquad(\b\in\U\st b\in B_\b).$$
By hypothesis, $T_{<\g}$ is a well-defined, $A$-linear contraction. Let
$$T\mapto{(B_{<\g})_\g}A{\sum_{k=0}^{n(\g)-1}b_k\xb^k}
{\left({1\over{n(\g)}}\right)\sum_{k=0}^{n(\g)-1}T_{<\g}(b_k)s_\g^k}.$$
Set $T_\g=T\circ\psi_\g\inv\map{B_\g}A$ where $\psi_\g$ is the isometric
isomorphism $(B_{<\g})_\g\to B_\g$ of Theorem 1.2.2.1.

It is routine to check that (i)-(iii) now hold with `$\le\g$' in place
of `$<\g$'. To illustrate the calculations, we verify that $\norm{T_\g}\le1$.
As $\psi_\g\inv$ is isometric, it is enough to check that $\norm T\le 1$.

Let $b_0,\ldots,b_{n(\g)-1}\in B_{<\g}$. Then $b=\sum_{k=0}^{n(\g)-1}b_k\xb^k$
is the general element of $(B_{<\g})_\g$ and
$$\eqalign{\norm{T(b)}&={1\over{n(\g)}}\norm{
\sum_{k=0}^{n(\g)-1}T_{<\g}(b_k)s_\g^k }\cr
&\le\sum_{k=0}^{n(\g)-1}\norm{T_{<\g}(b_k)}(\norm{s_\g^k}\bigm/n(\g))\cr
&\le\sum_{k=0}^{n(\g)-1}\norm{b_k}t_\g^k\qquad\hbox{ by hypothesis and }(*)\cr
&=\norm b.\cr}$$
So $\norm T\le 1$ as required.

By the transfinite recursion theorem, there exists a family of linear contractions
$(T_\a\map{B_\a}A)_\aiu$ such that
\item{(i)} for all $\aiu$, $T_\a\rest A=\id A$;
\item{(ii)} for all $\a,\b\in\U$, $\a\le\b<\g$, $T_\b\rest {B_\a}=T_\a$;
\item{(iii)} for all $\aiu$, $a\in A$, and $b\in B_\a$, $T_\a(ab)=aT_\a(b)$.

Recall that $B_\U=\bigcup_{\aiu}B_\a$.
We obtain the required map by defining
$$T_\U\mapto{B_\U}Ab{T_\a(b)}\qquad\qquad(\a\in\U\st b\in B_\a).$$

Note that by omitting requirement (iv) we can always obtain an $A$-linear map
$T_\U\map{B_\U}A$ with $T_\U\rest A=\id A$. It remains to verify that
$$\ab{T_\U(b)}\le\ab b\qquad (b\in B_\U)\eqno(\dag)$$
The argument, as Narmania noted, does not depend on ($*$), so we do
not have to make assumption ($*$) in the setting of uniform algebras.

We use the representation
of elements of $B_\U$ obtained in Lemma 1.2.2.2 to verify (\dag). This expands
remarks in $\ref{Nar}$ to the effect that for each $b\in B_\U$,
$\widehat{T_\U(b)}$ sends each $\o\in\O$
to a weighted average of the values of $\hat b$ on the fibre $\pi\inv(\o)$. Here
$\pi$ is the usual map $\Om{B_\U}\to\O$.

By Lemma 1.2.2.2, for each $b\in B_\U$ there exists a set
$U=\set{\a_1,\ldots ,\a_N}\in\U\fs$ and $b_m\in A$
($m\in L:=\times_{j=1}^N\set{0,1,\ldots,n(\a_j)-1}$) such that
$$\eqalign{b&=\sum_{m\in L} b_m \xi_{\a_1}^{m_1}\cdots\xi_{\a_N}^{m_N}\qquad\hbox{ and}\qquad
\norm b=\sum_{m\in L} \norm{b_m}
t_{\a_1}^{m_1}\cdots t_{\a_N}^{m_N}.\cr}$$

Thus $b\in B_\b$ where $\b=\max(\a_1,\ldots,\a_N)$. Inductively we see that
$$T_\U(b)=  {1\over\nn}\sum_{m\in L} b_m \smaj1\cdots\smaj N.$$

Let $\o\in\O$ and $\l_{{j},0},\ldots,\l_{{j},\naj-1}$ be the roots of $\o(\a_j)(x)$
listed according to multiplicity ($j=1,\ldots, N$). Then
$$\eqalign{\o(T_\U b)&=\sum_{m\in L} \o(b_m)\prod_{j=1}^N{\o(\smaj j)\over \naj}\qquad
\hbox{where }L=\times^N_{k=1}\set{0,\ldots,n(\a_k)-1}\cr
&= \sum_{m\in L} \o(b_m)\prod_{j=1}^N\sum_{l_j=0}^{\naj-1}{\l_{j,l_j}^{m_j}\over \naj}\cr
&=\sum_{m\in L}\o(b_m)\sum_{l\in L}\prod_{j=1}^N {\l_{j,l_j}^{m_j}\over \naj}
\qquad\hbox{ (by an exercise in indices)}\cr
&= \sum_{m,l\in L}\o(b_m)\prod_{j=1}^N {\l_{j,l_j}^{m_j}\over \naj}.\cr}$$
For each $l\in L$, choose $y_l=(\o,\l_l)\in\OU$ such that
$$\l_{l,\a_j}=\l_{j,l_j}\qquad(j=1,\ldots, N).$$
Thus for each $m,l\in L$,
$$ \o(b_m)\prod_{j=1}^N {\l_{j,l_j}^{m_j}}=(\o,\l_l)\left(b_m\prod_{j=1}^N \xmaj\right),$$
whence
$$\eqalign{\o(T_\U b) &= \sum_{l,m\in L} (\o,\l_l)\left( b_m\prod_{j=1}^N{\xmaj\over\naj}\right)\cr
&=\sum_{l\in L} (\o,\l_l)\left( \sum_{m\in L} b_m\prod_{j=1}^N{\xmaj\over\naj}\right)\qquad\hbox{ and so}\cr
\ab{\o(T_\U b)} &\le\sum_{l\in L}{1\over\nn}\ab{(\o,\l_l)(b)}\cr
&\le\ab b.\cr}$$
This completes the proof.\eop

\npar We illustrate the effect of $T^\U$ with an example for future reference.
Let $(A,X)$ be a uniform algebra
and $\a(x)$  be a monic polynomial over $A$ of degree $n$. Let $(\aA,X^\a)$ denote the
Cole extension of $A$ with respect to $\a(x)$. Then, writing
$p\mapto{X^\a}\C{(\k,\l)}\l$, elements of the form
$$f=\sum_{k=0}^{n-1}\ps(f)p^k\qquad\qquad(f_0,\ldots,f_{n-1}\in A)$$
are dense in $\aA$ (by Proposition 1.3.2.1).
Let $s_k$ be the $k$th Newton-sum of $\a(x)$ for each $k\in\No$. If we identify
$X$ with a subset of $\O$ in the usual way, we have for each $\k\in X$,
writing $\l_1,\ldots,\l_n$ for the roots of $\k(a)(x)$ repeated by multiplicity,
$$\eqalign{T^{\set\a}(f)(\k)&=\left({1\over n}\sum_{k=0}^{n-1}f_ks_k\right)(\k)\qquad\hbox{ by construction}\cr
&={1\over n}\sum_{k=0}^{n-1}f_k(\k)\sum_{j=1}^n\l_j^k\cr
&={1\over n}\sum_{j=1}^n\sum_{k=0}^{n-1}f_k(\k)\l_j^k\cr
&={1\over n}\sum_{j=1}^nf(\k,\l_j).\cr}$$
It now follows from the continuity of $T^{\set\a}$ that $T^{\set\a}(g)(\k)$ is this
average of values of $g$ over the fibre $\pi\inv(\k)$ for every $g\in\aA$.
  
As Narmania noted, the
restriction ($*$) is not a hindrance to our main applications, as we can still generate an
integrally closed extension of any normed
algebra by using only polynomials for which ($*$) is true (using the procedure
described in Section 1.2.4). We provide a simple argument, which is omitted from
$\ref{Nar}$, to justify this.

Since we can always choose the Arens-Hoffman norm-parameters
to be at least one, we assume this and can simplify the condition to $\norm{\saj j}\le \naj$.
Let $\a(x)=a_0+\cdots+a_{n-1}x^{n-1}+x^n\in A[x]$ and $\m>0$. Let
$\a^\m(x)=\m^na_0+\cdots+(\m a_{n-1})x^{n-1}+x^n$. By the continuity of algebraic
operations in $A$, we have $s_{\a^\m}^j\to0$ as $\m\to0\ (j=1,\ldots, n-1)$. We
can replace $\a(x)$ by $\a^\m(x)$ for $\m$ so small that ($*$) holds
for $\a^\m(x)$. If $\a^\m(\eta)=0$ then $\a(\eta/\m)=0$.

\npar It is a routine matter to extend the operator defined above to systems of algebraic extensions.
We give details, which are omitted by Narmania, for the benefit of the reader.

\prop 1.4.2.2 ($\ref{Nar}$). Let $((A_\s),(\tss\s\t)_{\s\le\t}\st\s,\t\le\u)$ be a system
of algebraic extensions. If $A_0=A$ is not a uniform algebra, suppose also that
for all $\t\le\u,\ \a\in\U_\t,\ j\in\set{0,\ldots,\na-1}$ we have
$$\norm{\saj j}\le \naj,$$
and that all norm-parameters are at least one.
Then there exists a family of unital, linear contractions $(T_{\s,\t}\map{A_\t}{A_\s})
_{\s\le\t\le\u}$ such that for $\s\le\t\le\u$,
\item{(i)} $T_{\s,\t}$ is $A_\s$-linear, and
\item{(ii)} $T_{\s,\t}\circ\tss\s\t=\id {A_\s}.$
\pf Let $\t\le\u$ and suppose that $(\Tss\r\s\map{A_\s}{A_\r})_{\r\le\s<\t}$ is a
family of contractions such that
\item{(i)} for all $\r\le\s<\t$, $a\in A_\r$, and $b\in A_\s$, $T_{\r,\s}(\tss\r\s(a)b)=aT_{\r,\s}(b)$, and
\item{(ii)} for all $\r\le\s<\t$, $T_{\r,\s}\circ\tss\r\s=\id {A_\r}.$

If $\t=\s+1$ for some $\s<\t$, then $A_\t=E(A_\s,\U_\s)$. By Theorem 1, there exists an
operator $\Tss\s\t\map{A_\t}{A_\s}$ such that $\norm {\Tss\s\t}\le 1$,
$$\Tss\s\t(\tss\s\t(a)b)=a\Tss\s\t(b)\qquad\qquad(a\in A_\s,b\in A_\t),$$
and $T_{\s,\t}\circ\tss\s\t=\id {A_\s}.$ Setting $\Tss\r\t:=\Tss\r\s\circ\Tss\s\t$ for
$\r\le\s$ yields a family of maps $(\Tss\r\s)_{\r\le\s\le\t}$ compliant with
requirements (i) and (ii) above.

Otherwise, $\t$ is a limit ordinal. If $\t=0$ we must set $\Tss00=\id A$. For $\t>0$ we have
$$(A_\t,(\tss\s\t)_{\s<\t})=\dl{}\left((A_\s),(\tss\r\s)\st\r\le\s<\t)\right)\eqno(*)$$
Fix $\r<\t$ and set $\Ss\r\s:=\Tss\r{\s\pri}\circ\tss\s{\s\pri}\map{A_\s}{A_\r}$
where $\max(\r,\s)\le\s\pri<\t$. Thus $\Ss\r\s$ is a well-defined and $A_\r$-linear contraction
since it is a composition of such maps.
Indeed
$$\Ss\r\s=\cases{\Tss\r\s &if $\r\le\s$, and\cr
\tss\s\r &if $\r\ge\s$.\cr}$$
We also have that for all $\s_1\le\s_2<\t$,
$$\eqalign{ \Ss\r{\s_2}\circ\tss{\s_1}{\s_2}
&=\Tss\r{\s_3}\circ\tss{\s_2}{\s_3}\circ\tss{\s_1}{\s_2}
\qquad\hbox{ where }\s_3\ge\max(\r,\s_2)\cr
&=\Tss\r{\s_3}\circ\tss{\s_1}{\s_3}\cr
&=\Ss\r{\s_1}\qquad\hbox{ since }\s_3\ge\s_2\ge\s_1.\cr}$$
Observe that $(A_\t,(\tss\s\t))$ is also the direct limit of ($*$) in the category of normed
spaces and contractions (see Appendix 2). By the universal property
of direct limits, there exists a contraction $\Tss\r\t\map{A_\t}{A_\r}$
such that for all $\s<\t$,
$$\Tss\r\t\circ\tss\s\t=\Ss\r\s=\cases{\Tss\r\s&$(\r\le\s)$\cr
\tss\s\r&($\s\le\r$).\cr}$$
In particular $\Tss\r\t\circ\tss\r\t=\Ss\r\r=\id{A_\r}$. Continuity
shows that $\Tss\r\t$ is $A_\r$-linear as we now check.

Since $D=\bigcup_{\s<\t}
\tss\s\t(A_\s)$ is dense in $A_\t$, it is enough to show that $\Tss\r\t\rest D$
is $A_\r$-linear. Let $b\in D$ and $a\in A_\r$. Since the algebras are
directed, we can assume that $b=\tss\s\t(b_\s)$ for some $\s\ge\r$ and
$b_\s\in A_\s$. Thus
$$\eqalign{\Tss\r\t(\tss\r\t(a)b)&=(\Tss\r\t\circ\tss\s\t)(\tss\r\s(a)b_\s)\cr
&=\Ss\r\s(\tss\r\s(a)b_\s)\cr
&=a\Ss\r\s(b_\s)\qquad\hbox{ since }\Ss\r\s=T_{\r,\s}\hbox{ is $A_\r$-linear}\cr
&=a\Tss\r\t(b).\cr}$$
The result now follows from the transfinite recursion theorem.\eop

\npar We now keep a promise made in Section 1.2.2 and use Proposition 2
to show that every uniform algebra in a system of Cole extensions of
a non-trivial uniform algebra is non-trivial. The following three results are
similar to ones contained in Theorem 3.5 of $\ref{LinIE}$ and Theorem 2.1 of $\ref{Col}$.

Proposition 3 is not essential for our purpose, but we include it
as it is interesting in its own right.

\prop 1.4.2.3. Let $\u$ be an ordinal number and $((A_\t)_{\t\le\u},
(\pss\s\t)_{\s\le\t\le\u})$ be a system of Cole extensions of the
uniform algebra $A=C(X)$ for some compact, Hausdorff space, $X$.
Then for all $\t\le\u$, $A_\t=C(X_\t)$.

\pf Let $\J=\set{\t\le\u\st A_\t=C(X_\t)}$, $\t\le\u$, and suppose
that $[0,\t)\subseteq\J$.

If there exists $\s<\t$ such that $\t=\s+1$ then $A_\t=A_\s^{\U_\s}$
for some set, $\U_\s$, of monic polynomials over $A_\s$. Let $B$ denote a
standard extension of $A_\s$ with respect to $\U_\s$. Since $A_\s=C(X_\s)$
is natural, it follows from Proposition 1.3.1.3 that $A_\t$ is natural, and
so $X_\t=\Om{A_\t}$. By Proposition 1.3.2.1 we have $\overline{B\htt}=A_\t$.
On the other hand, since $A_\s$ is dense in $C(X_\s)$ we have $\overline{B\htt}=C(X_\t)$
by Theorem 3.5 of $\ref{LinIE}$. So $\t\in\J$.

Otherwise $\t$ is a limit ordinal. It is clear that $0\in\J$ so
we may assume that $\t>0$. Then $A_\t=\lim_{\rightharpoondown}(A_\s)_{\s<\t}$.
Explicitly (see Appendix 2), $A_\t$ is the closure in $C(X_\t)$
of the algebra $D=\bigcup_{\s<\t}\pss\s\t(A_\s)$. As noted in $\ref{Col}$,
it follows from the Stone-Weierstrass theorem and inductive hypothesis that
$D$ is dense in $C(X_\t)$. So $A_\t=C(X_\t)$ and $\t\in\J$.

By the transfinite induction theorem, $\J=[0,\u]$ as required.\eop

\noindent Let $((A_\t)_{\t\le\u},
(\pss\s\t)_{\s\le\t\le\u})$ be a system of Cole extensions in which
for all $\t\le\u$, $A_\t$ is a uniform algebra on $X_\t$. Recall that $\pss\s\t$ is the
adjoint $C(X_\s)\to C(X_\t)\;;\;f\mapsto f\circ\Ps\s\t$ where
$\Ps\s\t$ is a continuous surjection $X_\t\to X_\s$. As usual, we shall
not distinguish between $\pss\s\t$ and $\pss\s\t\vert_{A_\s}$.

\lem 1.4.2.4. Let $\u$ be an ordinal number and $\AA=((A_\t)_{\t\le\u},
(\pss\s\t)_{\s\le\t\le\u})$ be a system of Cole extensions of a
uniform algebra $(A, X)$. Then the family
$$\CC=((C(X_\t))_{\t\le\u}, (\pss\s\t)_{\s\le\t\le\u})$$ is a system
of Cole extensions of $C(X)$ with respect to the same sets of

polynomials as are used to generate $\AA$.

\pf
The family $\CC$
clearly forms a direct system.
Let the set of polynomials generating $A_{\s+1}$ be denoted by $\U_\s$
for every $\s<\u$. Thus $\U_\s$ can be regarded as a set of
monic polynomials over $C(X_\s)$ for each $\s<\u$.
We verify the axioms given in Section 1.2.3
for $\CC$ to be a system of Cole extensions of $C(X)$ with respect to $(\U_\s)_{\s<\u}$.

Let $\t<\u$. First suppose that $\t=\s+1$ for some $\s<\t$.
As in the proof of Lemma 3 we have that
$$C(X_\s)^{\U_\s}=C(X_\s^{\U_\s})=C(X_\t).$$

Otherwise $\t$ is a limit ordinal. We can assume $\t>0$. Again, as in the second case
in Proposition 3, $C(X_\t)=\lim_\rightharpoondown(C(X_\s))_{\s<\t}$.

Thus $\CC$ is a system of Cole extensions of $C(X)$ with respect to
$(\U_\s)_{\s<\u}$.\eop

\noindent
Recall that if $E$ is a subset of a normed space, $(F,\norm\cdot)$,
and $f\in F$ then the {\it distance} of $f$ from $E$, is defined to be
$${\rm dist}\,(f, E):=\inf_{e\in E}\norm{e-f}.$$
It is standard (see, for example, p. 69 of $\ref{Sim}$)
that if $E$ is closed and $f\in F-E$ then
$\dist fE>0$.

\thm 1.4.2.5. Let $\u$ be an ordinal number and $\AA=((A_\t)_{\t\le\u},
(\pss\s\t)_{\s\le\t\le\u})$ be a system of Cole extensions of a
uniform algebra $(A, X)$. Then for all $f\in C(X)$,
$$\dist{\pss0\u(f)}{A_\u}=\dist fA.$$
In particular, $A_\u$ is non-trivial if $A$ is.

\pf By Lemma 4, $\CC=((C(X_\t)_{\t\le\u}, (\pss\s\t)_{\s\le\t\le\u})$ is a system
of Cole extensions of $C(X)$ with respect to the same sets of
polynomials as are used to generate $\AA$.
Let $(\Ts\s\t\map{A_\t}{A_\s})_{\s\le\t\le\u}$ and
$(\Tes\s\t\map{C(X_\t)}{C(X_\s}))_{\s\le\t\le\u}$ be the families
of linear contractions associated with $\AA$ and $\CC$ respectively
as in Proposition 2. The bulk of the remainder of this proof
is devoted to checking that $\Tes\s\t\vert_{A_\t}=\Ts\s\t$ for all
$\s\le\t\le\u$.

Let $\J=\set{\t\le\u\st\hbox{ for all }\r<\t,\ \Tes\r\t\vert_{A_\t}=\Ts\r\t}$.
Let $\t\le\u$ and suppose that $[0,\t)
\subseteq\J$.

First suppose that $\t$ is a limit ordinal. It is trivial that $0\in\J$ so
we can assume that $\t>0$. Then $A_\t$ is the subalgebra of
$C(X_\t)$ generated by $D=\bigcup_{\s<\t}\pss\s\t(A_\s)$. Let $\r<\t$.
Since $\Tes\r\t$ is continuous, it is enough to prove that
$\Tes\r\t\vert_D=\Ts\r\t\vert_D$. So let $g\in D$. As the
algebras
$(\pss\s\t(A_\s))_{\s<\t}$ are directed, there exists $\s\ge\r$
and $f\in A_\s$ such that $g=\pss\s\t(f)$. We have
$$\eqalign{\Tes\r\t(\pss\s\t(f))&=\Tes\r\s(f)\qquad(\hbox{by construction of $\Tes\r\t$; see Proposition 2})\cr
&=\Ts\r\s(f)\qquad\hbox{(by inductive hypothesis)}\cr
&=\Ts\r\t(\pss\s\t(f))\qquad\hbox{(by construction of }\Ts\r\t).\cr}$$
Hence $\Tes\r\t\vert_{A_\t}=\Ts\r\t$ and so $\t\in \J$.

Secondly suppose that $\t=\r+1$ for some $\r<\t$. So $A_\t=A_\r^{\U_\r}$
for some set, $\U_\r$, of monic polynomials over $A_\r$. It is enough
to prove that $\Tes\r\t\vert_{A_\t}=\Ts\r\t$ for then if
$\s<\r$ and $f\in A_\t$,
$$\eqalign{\Tes\s\t(f)&= (\Tes\s\r\circ\Tes\r\t)(f)\qquad\hbox{(by construction; see Proposition 2)}\cr
&= \Tes\s\r(\Ts\r\t(f))\cr
&= \Ts\s\r(\Ts\r\t(f))\qquad\hbox{(by hypothesis, since } \Ts\r\t(f)\in A_\r) \cr
&= \Ts\s\t(f)\qquad\hbox{(by construction of }\Ts\s\t),\cr}$$
from which $\Tes\s\t\vert_{A_\t}=\Ts\s\t$.

Consider elements of the form
$$g=\pss\r\t(f)p_{\a_1}^{m_1}\cdots p_{\a_N}^{m_N}$$
where $f\in A_\r$, $\a_1,\ldots,\a_N$ are distinct elements of $\U_\r$, and
$m_j\in\set{0,\ldots,n(\a_j)-1}$ ($j=1,\ldots, N$).
These elements span a dense subalgebra of $A_\t$.
Since $\Tes\s\t$ is linear and continuous, it is sufficient to
prove that $\Tes\s\t(g)=\Ts\s\t(g)$. But, by the construction in
Theorem 1, we have
$$\Tes\s\t(g)={1\over{n(\a_1)\cdots n(\a_N)}}
s_{\a_1}^{m_1}\cdots s_{\a_N}^{m_N}=\Ts\s\t(g).$$
Thus $\t\in\J$.

By the transfinite induction theorem, $\J=\U$. In particular we have
$\Tes0\u(A_\u)\subseteq A$.

Let $f\in C(X)$ and $r=\dist fA$. Since $\Tes0\u$ is a contraction
we have that for all $g\in A_\u$,
$$\eqalign{\norm{\pss0\u(f)-g}&\ge\norm{
(\Tes0\u\circ\pss0\u)(f)-\Tes0\u(g)}\cr
&=\norm{f-\Ts0\u(g)}\qquad \hbox{ (since $\Tes0\u\circ\pss0\u=\id{C(X)}$)}\cr
&\ge r\qquad\hbox{ (since $\Ts0\u(g)\in A$)}.\cr}$$
Hence $s:=\dist{\pss0\u(f)}{A_\u}\ge r$.

Finally, for all $g\in A$ we have, since $\pss0\u$ is isometric, that
$$\eqalign{\norm{g-f}&=\norm{
\pss0\u(g)-\pss0\u(f)}\cr
&\ge s\qquad\hbox{ (since }\pss0\u(A)\subseteq A_\u).\cr}$$
Thus $r\ge s$.\eop

\def\na{{\bf na}}


\def\Go#1{\Gamma_0\bigl(#1\bigr)}
\def\p#1#2{\pi_{#1,#2}}
\def\fibp#1#2#3{\pi_{#1,#2}^{-1}\left({#3}\right)}

\def\ku{\k_\u}

\def\si{{\s_i}}

\def\kl{(\k,\l)}

\def\fib#1{\pi^{-1}\left(#1\right)}
\def\loi{\l_0^{(i)}}
\def\lo1{\l_0^{(1)}}
\def\z{\zeta}
\def\pai{p_{\a_i}}
\def\lai{\l_{0,\a_i}}
\def\tf#1{\widetilde{f_#1}}
\def\T#1#2{T_{#1,#2}}
\def\p#1#2{\ps_{{#1},{#2}}}

\def\ui{{\u_1}}

\vskip 15pt
\noindent {\medtenbf 1.5. The \v Silov and Choquet Boundaries of the Extensions}

\vskip 15pt
\noindent {\bf 1.5.1. Introduction}\vskip 10pt

\noindent  We now use the construction of the previous section to investigate the boundaries
of algebraic extensions of uniform algebras. We shall use the notation $\Nbb X\k$ for the
set of neighbourhoods of a point $\k$ in a topological space, $X$.
Where no confusion is
likely, we shall omit the subscript labelling the space. If $f\map X\C$ is a bounded
function and $E\subseteq X$, we shall use the notation
$$\norm f_E:=\sup_{\k\in E}\ab{f(\k)}.$$

Recall that for a uniform algebra, $(A,X)$, a point $\k\in X$ is a {\it peak point}
for $A$ if there exists $f\in A$ such that
$1=f(\k)$ and for all $\k\pri\in X-\set\k$, $\ab{f(\k\pri)}<1$. The point $\k\in X$ is a
{\it strong boundary point for $A$} if for all neighbourhoods $V\in\Nb\k$,
there exists $f\in A$ such that $1=f(\k)=\norm f$ and $\norm f_{X-V}<1$.
 There are
various striking characterisations of this condition; see, for example, $\ref{Sto}$, pages 47--51.
A {\it boundary} for a set $S\subseteq C(X)$ is any subset of $X$ on which every element of $S$ attains
its maximum modulus.
The {\it Choquet boundary} of a uniform algebra is the set of strong boundary points.
We shall denote the Choquet boundary of $A$ by $\Go A$.

It is standard (see, for example, $\ref{Sto}$, p. 37) that $(A,X)$ has a minimal
closed boundary. It is called the {\it \v Silov boundary} of $A$ and is denoted
by $\S A$. It is often useful to note, as on p. 38 of $\ref{Sto}$, that $\k\in \S A$
if and only if for every $V\in\Nb \k$ there exists $f\in A$ such that
$\norm f_{X-V}<\norm f_V$. Furthermore, by replacing $f$ by $(f/\norm f)^N$ for
sufficiently large $N\in\N$, we can assume that $\norm f_{X-V}<\a<\norm f_V=1$ for
any given $\a\in(0,1)$.

It can be shown (see $\ref{Sto}$, p. 53) that $\Go A$ is dense in $\S A$. Another important
result (see $\ref{Sto}$, p. 54) is that if $X$ is metrizable (that is, the topology
on $X$ comes from a metric in the usual way; see $\ref{Kel}$, p. 118), then $\k\in X$ is a strong
boundary point if and only if it is a peak point for $A$.

We shall investigate the Choquet boundary of Cole extensions. First of all
we use some of
Lindberg's results from $\ref{LinIE}$ to describe the \v Silov boundaries in systems of Cole extensions.

\vskip 15pt
\noindent {\bf 1.5.2. The \v Silov Boundary}\vskip 10pt

\noindent Proposition 1.5.2.2 is surely known to
those who
worked in this area, such as Karahanjan, but we give a detailed statement and proof for the reader's convenience.
We require a lemma, which must have been used by many people.

\lem 1.5.2.1. Let $D$ be a dense subset of a uniform algebra $(A, X)$ and $F$ be a closed
boundary for $D$. Then $\S A\subseteq F$.

\pf Let $s\in\S A$ and $V$ be a neighbourhood of $s$. There exists $f\in A$ such that
$1=\norm f_V>{1\over2}>\norm f_{X-V}$. By the assumption on $D$ there is some $g\in D$
with $\norm{f-g}<{1\over4}$. Hence $\norm g_V>{3\over4}>\norm g_{X-V}$. So $\norm g$ is
attained by $\ab g$ on $V$ and not on $X-V$. Therefore $V\cap F\ne\emptyset$. This shows
that
$s\in\bar F=F$. So $\S A\subseteq F$ as required.\eop

\prop 1.5.2.2. Let $((A_\t,X_\t)_{\t\le\u}, (\pss\s\t)_{\s\le\t\le\u})$
be a system of Cole extensions of $(A,X)$
and $\t\le\u$. Then
$\S{A_\t}=\fibp 0\t{\S A}$.

\pf Let $\J=\set{\t\le\u \st \S{A_\t}=\fibp 0\t{\S A} }$. Let $\t\le\u$ and suppose
that $[0,\t)\subseteq\J$. Clearly $0\in\J$.

Suppose first that $\t=\s+1$ for some $\s<\t$. Then, if $\U_\s$ denotes the
set of monic polynomials over $A_\s$ generating $A_\t$, we have
by Proposition 1.3.2.1 and Theorem 3.5
of $\ref{LinIE}$, that
$$\S{A_\t}=\S{A_\s^{\U_\s}}=\fibp\s\t{\S{A_\s}}.$$
(Theorem 3.5 of $\ref{LinIE}$ asserts, in our notation, that $\S{\overline{(B_{\U_\s})\htt}}
=\pi_{\s,\t}^{-1}(\S {A_\s})$ where $B_{\U_\s}$ is a standard extension of $A_\s$ with
respect to $\U_\s$.) Hence, by hypothesis,
$$\eqalign{\S{A_\t}&=\fibp\s\t{ \fibp0\s{\S A} }\cr
&=(\Ps 0\s\circ\Ps\s\t)^{-1}(\S A)\cr
&=\fibp 0\t{\S A}.\cr} $$
So $\t\in\J$.

Now suppose that $\t>0$ is a limit ordinal,
and so $X_\t=\lim_{\leftharpoondown} X_\s$. (The argument to show that $\t\in\J$ is adapted
from the proof of Theorem 11 of $\ref{Lei}$ on p. 214.) We first show that
$\fibp0\t{\S A}\subseteq\S{A_\t}$. Let $\k\in\fibp0\t{\S A}$.

Let $U$ be a basic, open neighbourhood of  $\k$. Thus there exist $\s_1,\ldots,
\s_n<\t$ and open $U_i\subseteq X_{\si}$ such that
$U=\cap_{i=1}^n \fibp{{\s_i}}\t{U_i}$.

Let $\s=\max(\s_1,\ldots,\s_n)$. Then $\Ps\si\t=\Ps\si\s\circ\Ps\s\t$ and so
$\fibp\si\t{U_i}=\fibp\s\t{W_i}$ where $W_i=\fibp{{\s_i}}\s{U_i}$
$(i=1,\ldots,n)$. Thus
$W_i$ is open in $X_\s$. Moreover $U=\fibp\s\t W$ where
$W=\cap_{i=1}^nW_i$, an open set in $X_\s$.

Since $\k\in U$, $\k_\s=\Ps\s\t(\k)\in W$. We have
$\Ps0\s(\k_\s)=(\Ps0\s\circ\Ps\s\t)(\k)=\k_0:=\Ps0\t(\k)\in\S A$ so
$\k_\s\in\fibp0\s{\S A}=\S{A_\s}$. Since the strong boundary points of $A_\s$ are
dense in $\S{A_\s}$ there exists $s_\s\in W$ and $f\in A_\s$ with
$1=f(s_\s)=\norm f>\norm f_{X_\s-W}$.

Let $s\in\fibp\s\t{s_\s}$ and $g=\pss\s\t(f)$. Then
$\Ps\s\t(s)=s_\s\in W$ so $s\in U$ and $1=g(s)=\norm g$. For
$\o\in X_\t-U$ we have $\Ps\s\t(\o)\not\in W$ so $\ab{g(\o)}<1$.

So $s$ is in $U$ and the Choquet boundary of $A_\t$. Since $U$ is
arbitrary and $\S{A_\t}$ is the
closure of the Choquet boundary, $\k\in\S{A_\t}$.

To obtain the reverse inclusion it is enough to show that
$\fibp0\t{\S A}$ is a (closed) boundary for the dense subset $\bigcup_{\s<\t}
\pss\s\t(A_\s)$ by Lemma 1.

Let $\s<\t$ and $f\in A_\s$. Then $\norm{\pss\s\t(f)}=\norm f
=\ab{f(\k_\s)}$ for some $\k_\s\in\S{A_\s}=\fibp0\s{\S A}$. Let
$\k\in\fibp\s\t{\k_\s}$. We have
$\Ps0\t(\k)=(\Ps0\s\circ\Ps\s\t)(\k)=\Ps0\s(\k_\s)\in\S A$ and
$\norm{\pss\s\t(f)}=\ab{\pss\s\t(f)(\k)}$. Therefore $\S{A_\t}\subseteq
\fibp0\t{\S A}$ and so $\t\in\J$.

By the transfinite induction theorem, $\J=[0,\u]$ as required.\eop

\vskip 15pt
\noindent {\bf 1.5.3. The Choquet Boundary}\vskip 10pt

\noindent
Let $\a(x)$ be a monic polynomial over the uniform algebra $(A,X)$. We shall write $(\aA, X^\a)$,
for the simple Cole extension of $A$ by $\set\a$. We shall show that a similar result
to Proposition 1.5.2.2 holds, namely:
$$\Go\aA=\fib{\Go A}\eqno(1)$$
For extensions by infinitely many
polynomials we only have an inclusion between the Choquet boundary of
the extension and the preimage of the Choquet boundary of the original algebra.
We do not know if
this inclusion can be strict. The kind of relationship in $(1)$ can break down eventually,
as we now explain.

Let $\o_0$ denote the first infinite ordinal.
There is an algebra, $(A_{\o_0},X_{\o_0})$, in a
system of Cole extensions of $(A,X)$
(see $\ref{Col}$) where $X_{\o_0}$ is metrizable, every point of $X_{\o_0}$ is a
peak point for $A_{\o_0}$, but not every point of $X$ is a peak point for $A$.
It follows from the facts mentioned in the introduction to this section that in this
case $\Go{A_{\o_0}}\ne\Ps0{\o_0}^{-1}{\Go A}$.

Some results on the Choquet boundary of algebraic extensions are already known.
The following proposition is obtained from
Lemma 4.2.8(i) of $\ref{FeiThes}$.

\prop 1.5.3.1. Let $E\subseteq X$ be closed and $\U$ be a set of monic polynomials over $A$
of the form $x^2-f$ where each $f(E)\subseteq\set0$. Then
$\Go\UA\cap\fib E=\pi\inv(\Go A\cap E)$.\eop

\npar A general result, which is useful here, is the following.
Recall that a {\it peak set} for a uniform algebra, $(A,X)$, is a subset, $K$, of $X$ for
which there exists some $f\in A$ such that $f(K)=\set 1$ and for all $\k\in{X-K}$,
$\ab{f(\k)}<1$. An intersection of peak sets is called a
{\it peak set in the weak sense}.

\lem 1.5.3.2. Let $A$ and $B$ be uniform algebras on the compact, Hausdorff spaces $X$ and $Y$
respectively. Suppose that $\pi\map Y X$ is a finite-to-one, continuous surjection which
induces an isometric monomorphism $\ps\mapto A Bf{f\circ\pi}$. Then $\fib{\Go A}\subseteq \Go B$.

\pf Let $\k_0\in\Go A$. Then, by Theorem 2.3.4 of $\ref{Bro}$, $\set{\k_0}$ is an intersection of peak sets for $A$. Let $J$ be an index set and $(K_j)_{j\in J}$ be a family of peak
sets for $A$ whose intersection is $\set{\k_0}$.

Now if $K$ is a peak set for $A$ then $\fib K$ is a peak set
for $B$ since for some $f\in A$ we have $f(K)=\set 1$ and for all $\k\in{X-K}$,
$\ab{f(\k)}<1$. Put $g=\ps(f)$.
Then $g\bigl(\fib K\bigr)=\set 1$ and for $y\not\in\fib K$ we have $\ab {g(y)}<1$

Let $\fib{\k_0}=\set{y_1,\ldots,y_u} $. Therefore $\set{y_1,\ldots,y_u} =
\cap_{j\in J}\fib{K_j}$ is a peak set in the weak sense for $B$. It is clear that
$\set{y_1}$ is a peak set (in the weak sense) for $B\vert_{\fib{\k_0}}$. By
Theorem 2.4.4 of $\ref{Bro}$ we have that $\set{y_1}$ is a peak set in the weak sense
for $B$. But this is another characterisation ($\ref{Bro}$, p. 96)
of strong boundary points so $y_1\in\Go B$,
and the result follows.\eop

\noindent The following definition is useful.

\dfn 1.5.3.3. Let $\a(x)$ be a monic
polynomial over a normed algebra $A$.
Let $\o\in\Om A$ and $\e>0$. We shall say that $\e$ {\it separates the roots} of
$\a$ at $\o$ if for each pair of distinct roots $\l_1,\l_2$ of $\o(\a)(x)$ we have
$\ab{\l_1-\l_2}\ge\e$.

\npar We require the following lemma. It is used repeatedly in $\ref{LinFact}$,
and we reproduce Lindberg's justification, with
extra detail, for the convenience of the reader.

\lem 1.5.3.4. Let $\a(x)=a_0+\cdots+a_{n-1}x^{n-1}+x^n$ be a monic polynomial
over the normed algebra $A$, $\k_0\in\O$, $V_0\in\Nb{\k_0}$, and $\r_0>0$.
Let $\pi$ be the usual projection $\Oa\to\O$
where $\Oa$ is
the space of closed, maximal ideals of the Arens-Hoffman extension $\Aa$.
Let $p_\a\mapto\Oa\C\ol\l$.
Let $\pi\inv(\k_0)=\set{y_i\st i=1,\ldots,u}$. Then there exists an open set
$V\in\Nb{\k_0}$ such that $V\subseteq V_0$, $\pi\inv(V)=\udot_{i=1}^uW_i$,
and for each $i=1,\ldots,u$, $W_i$ is an open neighbourhood of $y_i$ such
that $\pi(W_i)=V$ and $p_\a(W_i)\subseteq B_\C(\l_i,\r_0)$.

\pf Let $\l_1,\ldots,\l_u$ be the distinct roots of
$$\k_0(\a)(x)=\k_0(a_0)+\cdots+\k_0(a_{n-1})x^{n-1}+x^n\eqno(*)$$ with
$y_i=(\k_0,\l_i)$ ($i=1,\ldots,u$). Let $\rho\in(0,\r_0)$ be such that
$3\rho$ separates the roots of $\a$ at $\k_0$.
It is well-known from Rouch\' e's theorem (see $\ref{Rao}$, p. 134)
that the roots of ($*$) `depend continuously on the coefficients'. That is,
there exists $\d>0$ such that if
$$\b_j\in B_j:= B_\C(\k_0(a_j),\d)\qquad\qquad(j=0,\ldots,n-1)$$
then for each $k\in\set{1,\ldots,u}$, $B_\C(\l_k,\rho)$
contains $m_k$ roots (counted according to multiplicity)
of $\b(x)=\b_0+\cdots+\b_{n-1}x^{n-1}+x^n$
where $m_k$ is the multiplicity of $\l_k$.
Since $\b(x)$ has at most $n=\sum\nolimits_{k=1}^um_k$ roots, every
root of $\b(x)$ belongs to $B_\C(\l_k,\r)$ for some $k\in\set{1,\ldots,u}$.

Let $V\subseteq V_0\cap\bigcap_{j=0}^{n-1}\widehat{a_j}\inv(B_j)$ be an open
neighbourhood of $\k_0$. Let $W_k=\pi\inv(V)\cap p_\a\inv(B_\C(\l_k,\rho))$.
Then $W_k$ is
an open neighbourhood of $y_k$ ($k=1,\ldots,u$).

Clearly $\udot_{k=1}^u W_k\subseteq\pi\inv(V)$.
We have $\pi\inv(V)\subseteq\cup_{k=1}^uW_k$ for if $\ol\in\Oa$ with
$\o\in V$ then the root, $\l$, of $\o(\a)(x)$, satisfies $\ab{\l-\l_k}<\r$
for some $k\in\set{1,\ldots,u}$. Thus $\ol\in W_k$.

Let $k\in\set{1,\ldots,u}$. Evidently $\pi(W_k)\subseteq V$. We show that $V\subseteq \pi(W_k)$.
If $\o\in V$,
let $\m_1,\ldots,\m_n$ be a list of the roots of $\o(\a)(x)$ repeated
according to multiplicity. By the choice of $V$ and $\d$ we have
$\m_j\in B_\C(\l_k,\rho)$ for at least one $j\in\set{1,\ldots,n}$. Thus $(\o,\m_j)\in W_k$
and so $\pi(W_k)=V$ as required.

Finally it is clear from the construction that $p_\a(W_k)\subseteq B_\C(\l_k,\rho)$.\eop

\lem 1.5.3.5.
Let $\a(x)$ be a monic polynomial over the uniform algebra $(A,X)$.
Let $\k_0\in X$ and the distinct elements of $\fib{\k_0}$ be $y_1,\ldots,y_u$ where
$y_i=(\k_0,\loi)$ and the multiplicity of $\loi$ as a root of
$\eh{\k_0}(\a)(x)$ is $m_i\ (i=1,\ldots,u)$. Let $2\e>0$ separate the roots of $\a$ at
$\k_0$ and $\eta\in(0,1)$. Then there exist $V\in\Nb{\k_0}$ and $g\in\aA$ such that
$$g(y_i)=\cases{1 & $i=1$\cr
	0 & $i>1$,\cr }$$
and, defining $\zeta>0$ by
$${1\over\z}=\sup\left(\ab{g\kl}\st\kl\in\fib V\cap p_\a\inv(B_\C(\lo1,\e)) \right),$$
we have $\eta<\z\le 1$.

\pf
We may assume that $u>1$.
Let $\r_0<\min(\e,\bigl[\bigl({1\over{2\eta}}\bigr)^{1/(u-1)}-1\bigr]\e)$. By
Lemma 4 there exist $V\in\Nb{\k_0}$ and $W_i\in\Nb{y_i}$ such that
$\fib{V}=\udot_{i=1}^uW_i$ and $p_\a(W_i)\subseteq B_\C(\loi,\r_0)$ ($i=1,\ldots,u$).
Note that if $\l\in B(\lo1,\e)$ and $\o(\a)(\l)=0$ for some $\o\in V$
then $\ol\in W_1$ for otherwise $\ol\in W_i$ for some $i>1$ and so
$$2\e\le\ab{\lo1-\loi}\le\ab{\lo1-\l}+\ab{\l-\loi}<\e+\r_0<2\e,$$
a contradiction.

Let
$$g=\prod_{i=2}^u \left( {{p_\a-\loi}\over{\lo1-\loi}} \right).$$
Then $g$ takes the required values on $\pi\inv(\k_0)$.

Since $g(y_1)=1$ we have $\z^{-1}\ge 1$

Suppose $\k\in V,\ \l$ is a root of $\eh\k(\a)(x)$, and $ \ab{\l-\lo1}<\e $. Then
$$\eqalign{ \ab{g\kl}
	&= \prod_{i=2}^u \left\vert { {\l-\lo1+\lo1-\loi} \over{\lo1-\loi} } \right\vert \cr
	&\le \prod_{i=2}^u \left( 1 + { \ab{\l-\lo1} \over {\ab{\lo1-\loi}} } \right) \cr
	&\le \prod_{i=2}^u
\left(1 +  { {\e\biggl[\bigl({1\over{2\eta}}\bigr)^{1/(u-1)}-1\biggr]}\over{\e} }  \right)
\qquad\hbox{(since $\ol\in W_1$)}\cr
	&= {1\over{2\eta}} \cr } $$
Hence $1/\z<1/\eta$ as required.\eop

\npar We can now describe the Choquet boundary of $\aA$ as promised.

\thm 1.5.3.6. Let $\aA$ be the Cole extension of a uniform algebra $(A,X)$
by the monic polynomial $\a(x)$. Then $\Go\aA=\fib{\Go A}$.

\pf By Lemma 2 we have $\fib{\Go A} \subseteq \Go\aA $. To complete the proof we
use Gonchar's `$\a$-$\beta$' condition (see $\ref{Bro}$, p. 99) to show that
$\pi\left( \Go\aA \right)\subseteq \Go A$. Let $(\k_0,\lo1)=y_1\in \Go\aA$ where
the notation for $\fib{\k_0}$ is as in Lemma 4.
To show ${\k_0}$ is a strong boundary point of $A$ it is enough to show that
for all $V_0\in\Nb{\k_0}$ there exists $f\in A$ with $\norm f\le 1, f(\k_0)=3/4$ and
$\norm f_{X-V_0}<1/7$.

Let $2\e$ separate the
roots of $\a(x)$ at $\eh{\k_0}$. By Lemma 5,
there exists $g\in\aA$ and $V\in\Nb{\k_0}$
such that $g(y_k)=1\ (k=1)$, $g(y_k)=0\ (k>1)$,
and $\eta=7/8<\zeta\le 1$ where
$$\zeta\inv=\sup\left(\ab{g\kl}\st\kl\in\fib{V}\hbox{ and }\ab{\l-\l_0^{(1)}}<\e\right).$$
By Lemma 4, we
can further assume
that $V\subseteq V_0$, that $V$ is open, and that $\fib V=\udot_{i=1}^u W_i$
where $W_i$ is an open neighbourhood of $y_i$ such that $\pi(W_i)=V\ (i=1,\ldots,u)$. The
result now follows if we can find $f\in A$ with $\norm f\le 1, f(\k_0)=3/4$ and
$\norm f_{X-V}<1/7$.

By hypothesis there is some $h\in\aA$ with $h(y_1)=1=\norm h$ and $\norm h_{X^\a-W_1}
<1.$ Since $W_1$ is open we can replace $h$ by $h^N$ for some suitably large $N\in\N$
so that $\norm{h^N}_{X^\a-W_1}<1/(7n\norm g)$; we shall assume this has occurred.

Set $F=(3/4)hg$ and let $f=(n/m_1)T(F)$ where $T$ is the operator
$\aA\to A$ as in Section 1.4. Explicitly, for $\k\in X$ we have
$$f(\k)={1\over {m_1}}\sum_{i=1}^n F(\k, \l_i(\k) )$$
where $\l_1(\k),\ldots,\l_n(\k)$ are the roots of $\eh\k(\a)(x)$ repeated
multiplicitywise. Thus
$$f(\k_0)={3\over {4m_1}}\sum_{i=1}^n h (\k_0, \l_i(\k_0) ) g (\k_0, \l_i(\k_0))={3\over4}. $$

Now suppose $\k\in X-V$. Then $\fib\k\cap W_1=\emptyset$ so
$$\eqalign{ \ab{ f(\k) } &\le (3/4)\sum_{i=1}^n \ab{h(\k, \l_i(\k))}\ab{g(\k, \l_i(\k))} \cr
	&< (3/4)\sum_{i=1}^n {1\over {7n}}<1/7 .\cr} $$
So $\norm f_{X-V}<1/7$ (since $V$ is open).

It remains to check that $\norm f_X\le 1$; we need to check that $\norm f_V\le 1$.

Let $\k\in V$ and $\l_1(\k),\ldots,\l_n(\k)$ be relabelled if necessary so that
$\l_i(\k)\in W_1\ (i=1,\ldots,m_1)$ and $\l_i(\k)\not\in W_1$ for $i>m_1$. Thus we have
$$\eqalign{ \ab{ f(\k) } &< (3/4){1\over{m_1}}
\sum_{i=1}^{m_1} \ab{h(\k, \l_i(\k))}\ab{g(\k, \l_i(\k))} +
\left({3\over{4m_1}}\right)\sum_{i=m_1+1}^n{1\over {7n}} \cr
	&< {1\over7}+{1\over{m_1}}\sum_{i=1}^{m_1} (3/4) \ab{g(\k, \l_i(\k))}  .\cr} $$

Now by the choice of $g$ we have $\ab{g(\k, \l_i(\k))} \le 1/\z\ (i=1,\ldots,m_1)$.
But $7/8<\z$ so $1/\z<8/7$ and so
$$ (3/4) \ab{g(\k, \l_i(\k))}  <(3/4)(8/7) =6/7\qquad\qquad(i=1,\ldots,m_1)$$
whence $\norm f_V\le 1/7+6/7=1$.\eop

\npar For multiple extensions we have the following result.

\prop 1.5.3.7. Let $\U$ be a set of monic polynomials over the uniform algebra $A$.
Then we have
$$\fib{\Go A}\subseteq\Go\UA.$$

\pf Let $\k_0\in\Go A$ and $y_0=(\k_0,\l_0)\in\fib{\k_0}$.

Let $W_0\in\Nb{y_0}$. We show there exists $f\in\UA$ with $\norm f\le1=f(y_0)$ and
$\norm f_{X^\U-W_0}<1$.

We may assume $W_0$ is a basic open neighbourhood of the form
$$\fib{V_0}\cap\cap_{i=1}^up_{\a_i}^{-1} \left(B_\C\bigl(\lai,r\bigr)\right) $$
for some $V_0\in\Nb{\k_0}$, distinct $\a_1,\ldots,\a_u\in\U$, $\lai\in\C$ (where
$\eh{\k_0}(\a_i)(\lai)=0\ (i=1,\ldots,u)$), and $r>0$.

From Theorem 6 we have $\Go{A^{\a_i}}=\pi_i^{-1}\left(\Go A\right)\ (1=1,\ldots,u)$ where
$\pi_i$ is the canonical map $X^{\a_i}\to X$. Similarly let
$p_i\mapto {X^{\a_i}}\C\kl\l\ (i=1,\ldots,u)$.

Set $V_i=\pi_i^{-1}(V_0) \cap p_i^{-1}\left(B_\C\bigl(\lai,r\bigr)\right) \
(\in\Nb{\k_0,\lai)}\ (i=1,\ldots,u)$. Then there exists $f_i\in A^{\a_i}$ with
$1=f_i(\k_0,\lai)=\norm{f_i}>\norm{f_i}_{X^{\a_i}-V_i}$.

Now $f_i$ is a uniform limit of polynomials in $\pi^*_i(A)\cup\set{p_i}$ and may
be regarded as an element of $\UA$. More formally, we have for every $\a\in \U$ an
isometric embedding
$$\aA\mapright{\ps_{\a,\U}}(\aA) ^{\U-\set\a} \mapright{\rho^*}\UA$$
in which $\rho^*$ is an isometric isomorphism. It is induced by the homeomorphism
$$\rho\mapto{X^\U}{(X^\a)^{\U-\set\a}}\kl{ \left((\k,\l_\a), (\l_\beta)_{\beta\ne\a}\right)}.$$
Let $\tf i\in\UA$ denote $f_i$ under this identification. Let $f=\tf1\cdots\tf u$. Then
$f\in\UA$ and we have
$$f(y_0)=\prod_{i=1}^u\tf i(\k_0,\l_0)
=\prod_{i=1}^u f_i (\k_0,\lai) =1.$$ Also
$\norm f\le \prod_{i=1}^u\norm{\tf i} =  \prod_{i=1}^u\norm{f_i} =1$.
Finally for $y=\kl\not\in W_0$ we have $\kl\not\in\fib{V_0}$ or
$\kl\not\in \pai^{-1}\left(B_\C\bigl(\lai,r\bigr)\right)$ for some $i\in\set{1,\ldots,u}$.
So there exists $i\in\set{1,\ldots,u}$ such that $ (\k_0,\lai) \not\in V_i$ and
consequently $\ab{ f_i(\k_0,\lai) }<1$. So $\ab{f(y)}<1$.\eop

\def\p{\pi}	

\vfill\eject		


\noindent {\bigtenbf Chapter 2}
\vskip 10pt
\noindent {\bigtenrm Further Properties and Applications of Algebraic Extensions}
\vskip 25pt


\noindent
In this chapter we present some miscellaneous applications of the constructions of
the last chapter.


\def\Se{({\cal S})}

\def\xk{\xi_{\a_k}}

\def\xo{\xi_{\a_0}}
\def\too{t_{\a_0}}
\def\tto#1{\theta_{#1,\oo}}


\vskip 15pt
\noindent {\medtenbf 2.1. Integrally Closed Rings}

\vskip 15pt
\noindent {\bf 2.1.1. Introduction}
\vskip 10pt

\noindent It is  occasionally useful to be aware of the presence of integrally closed,
countable, dense subsets of a normed algebra; see, for example, $\ref{Col}$, Theorem 2.5.
We show that a normed algebra, $A$, with a countable,
dense, subring, $R$, can be extended to a normed algebra, $B$, with an integrally closed,
countable, dense subring. This problem can be solved in the same category, $\na$, $\Ba$,
or $\ua$, as $A$.

\vskip 15pt
\noindent {\bf 2.1.2. Integrally Closed, Countable, Dense Subrings}

\npar The following lemma must be well-known.

\lem 2.1.2.1. Let $R$ be a countable, infinite ring and
$\set{x_1,x_2,\ldots}$ be a sequence of commuting indeterminates. Then
$R[x_1,x_2,\ldots]$ is countable.

\pf Since $R$ is an infinite ring, $\# R[x_1]=\# R=\aleph_0\;(:=\ab\N)$ (see, for
example, $\ref{4DIS}$, Corollary A3.4, or
$\ref{Hal}$ for more general results). By induction we see
that $R[x_1,\ldots,x_n]$ has cardinality $\aleph_0$ for all $n\in\N$. The
result is now obvious since there is a surjection
$\bigcup_{n\in\N}R[x_1,\cdots,x_n]\to R[x_1,x_2,\ldots]$.\eop

\lem 2.1.2.2. Let $A$ be a normed algebra with a
countable, dense subring, $R$. Let $\U$ be a set of monic polynomials over $R$ and $B_\U$ be a
standard normed extension of $A$ with respect to $\U$. Then $R_1$, the
ring of polynomials in the standard roots, $\set{\xi_\a\st\a\in\U}$, over $R$,
is dense in $B_\U$.

\pf (Refer to Section 1.2 for the notation for standard extensions.)
By Lemma 1 we can assume that $\U=\set{\a_0,\a_1,\ldots}$ and take the
well-ordering on $\U$ given by $\a_j\le\a_k$ if and only if $j\le k$ for
all $j,k\in\No$. Set
$$B_{<k}=\cases{\cup_{l<k} B_l &if $k>0$, and\cr
A &if $k=0$.\cr}$$
By definition, $B_\U=\bigcup_{k\in\No} B_k$ where for all
$k\in\No$, there is a canonical, isometric isomorphism
$$\psi_k\map{B_{<k}[x]/(\a_k(x))}{B_k}.$$

Recall that the standard root of $\a_k(x)$ is $\xi_{\a_k}=\psi_k(\xb)$.
If $k=0,\e>0,$ and $b\in B_0$ are given then there exist (unique)
$b_0,\ldots,b_{n-1}\in A$ such that $b=\sum_{j=0}^{n-1}b_j\xo^j$ where
$n=n(\a_0)$ is the degree of $\a_0(x)$. By hypothesis there exist
$c_0,\ldots,c_{n-1}\in R$ such that $\norm{b-\sum_{j=0}^{n-1}c_j\xo^j}
=\sum_{j=0}^{n-1}\norm{b_j-c_j}\too^j<\e$. Since $\sum_{j=0}^{n-1}c_j\xo^j\in R_1\cap B_0$
we have that $R_1\cap B_0$ is dense in $B_0$.

Now let $k\in\N$ and suppose that $R_1\cap B_l$ is dense in $B_l$
for all $l<k$. Let $b\in B_k$ and $\e>0$. Let
$n=n(\a_k)$. There exist $b_0,\ldots,b_{n-1}\in B_{<k}$ such that
$b=\sum_{j=0}^{n-1}b_j\xk^j$. By the comments above there is some $l<k$
such that $b_0,\ldots,b_{n-1}\in B_l$. A similar argument to the case when
$k=0$ now shows that $R_1\cap B_k$ is dense in $B_k$.

By the principle of mathematical induction, $R_1\cap B_k$ is dense in $B_k$
for all $k\in\No$ and therefore $R_1$ is dense in $B_\U$ as required.\eop

\npar We can now establish the result mentioned in the introduction.

\prop 2.1.2.3. Let $A$ be a normed algebra with a countable,
dense, subring, $R$. Then there exists a normed algebra, $B$, in the same category, $\na$, $\Ba$,
or $\ua$, as $A$, such that $B$ has an integrally closed,
countable, dense subring.

\pf Let $R_1$ be as in Lemma 2.
Note that there exists a surjection $R[x_1,x_2,\ldots]\to R_1$ so
$R_1$ is countable. It follows from Proposition 1.3.2.1
that $R_1$ and $\widehat{R_1}$ are dense in the Narmania and Cole extensions
(where applicable)
of $A$ by $\U$, $A_\U$ and $A^\U$, respectively. Let $A_1$ and $\ths1$ denote
the extension and embedding in the appropriate categories. (See Section 1.2.3.)

By induction we can construct a direct sequence
$$A=A_0\ \mapright{\ths1}\ A_1\ \mapright{\ths2}\ A_2\ \mapright{\ths3}\ \cdots\eqno\Se$$
of normed algebras and isometric embeddings in each of the categories where
for all $k\in\No$, $A_k$ has a countable, dense subring, $R_k$ and every monic
polynomial over $R_k$ has a solution in $R_{k+1}$.

It only remains to take the direct limit, $A_\oo$, of $\Se$ in the appropriate
category. (See Appendix 2.) Thus $A_\oo$ has a dense subalgebra $D=
\bigcup_{k\in\No}\tto k(A_k)$ where the $\tto k$ are the canonical maps
$A_k\to A_\oo$.

Let $R_\oo$ be the countable, dense set $\bigcup_{k\in\No}\tto k(R_k)\subseteq D$.
It is easy to check that $R_\oo$ is an integrally closed ring.\eop


\def\BA{ {B_{<\a}} }
\def\fib#1{\pi\inv\left(#1\right)}
\def\fibp#1#2#3{\p{#1}{#2}^{-1}\left({#3}\right)}
\def\D#1{\Delta_{#1}}

\def\gt{\tilde g}

\def\OBb{\Om{B_\b}}

\def\tBb{\widetilde{B_\b}}
\def\p#1#2{\pi_{#1,#2}}

\def\ku{\k_\u}


\vskip 15pt
\noindent {\medtenbf 2.2. Regularity and Localness}

\vskip 15pt
\noindent {\bf 2.2.1. Regularity}\vskip 10pt

\noindent
The notion of regularity applies to any algebra; see, for example, $\ref{Pal}$, Section 7.2.
A unital normed algebra, $A$ is {\it regular} if for each closed set $E\subseteq\O$ and
$\o\in\O-E$ there exists $a\in A$ such that $\ha(E)\subseteq\set0$ and $\ha(\o)=1$,
where $\O$ is the space of closed, maximal ideals of $A$.

Algebraic extensions in the category $\Ba$ seem to
behave well with respect to regularity; by one of Lindberg's theorems in $\ref{LinAE}$ any
Arens-Hoffman
extension of a regular Banach algebra is regular. Results in
$\ref{Kar}$ and $\ref{FeiNTSR}$ show
that Cole's universal root extensions for adjoining square roots are regular if the
base algebras are. In contrast,
in the category of normed algebras, regularity is not always preserved by
Arens-Hoffman extensions. We shall give an example to demonstrate this
in Section 2.2.2.

First of all we generalise the result of Lindberg mentioned above to
standard extensions. Please refer to Section 1.2.2 for the notation
for standard extensions.

\prop 2.2.1.1. Let $A$ be a regular, commutative, unital, normed algebra
and $B_\U$ a standard normed extension with respect to a set of monic
polynomials $\U\subseteq A[x]$. Then the completion of $B_\U$
is regular.

\pf Without loss of generality we may assume that
$\U$ has a maximum element, $\a_1$. Let
$${\cal J}=\set{\a\le\a_1\st \widetilde{B_\a}\hbox{ is regular }}.$$
Let the least element of $\U$ be $\ao$. Let
$\a\in\U$ and suppose that $[\ao,\a)\subseteq{\cal J}$. If $\a=\ao$
then
$B_\ao$ is isometrically isomorphic to $A[x]\bigm/(\ao(x))$ and so its
completion can be identified with $\tilde A[x]\bigm/(\ao(x))$. This
follows easily from the universal property of Arens-Hoffman extensions and the
uniqueness of completions. (The fact was stated for restricted sets of
polynomials in $\ref{Heu}$; for a proof of the general case see
Theorem 3.13 of $\ref{4DIS}$.) Now $\tilde A$ is regular since $A$ is regular
and dense in $\tilde A$ (so they have the same continuous-character space). Hence
by Lindberg's theorem in $\ref{LinAE}$, $\widetilde{B_\ao}$ is regular.

Now suppose $\a>\ao$. Then (up to isometric isomorphism, as usual)
$B_\a=\BA[x]\bigm/(\a(x))$ where $\BA=\bigcup_{\b<\a} B_\b$. By
Lindberg's result, it is enough to prove that $\widetilde{\BA}$
is regular.

Let
$$\eqalign{y,y\pri\in\Om{\widetilde\BA}&=\Om\BA\cr
&=\set{ (\o,\l)\in\OA\times \C^{[\ao,\a)}\st\hbox{ for all }\b<\a\quad
\o(\b)(\l_\b)=0}.\cr}$$
See Section 1.3 for a description of the space of closed, maximal ideals
of a standard extension. Note that for any $\b<\a$, $B_{<\a}$ and $B_\b$
are standard extensions of $A$. There is therefore a continuous
surjection
$$\p\a\b\mapto{\Om\BA}{\OBb}{(\o,\l)} { (\o,(\l_\g)_{\g\le\b}) }.$$
The effect of $\Ps\a\b$ is simply restriction, under the identifications made in
Section 1.3.

Let $y=(\o,\l)\ne(\o\pri,\l\pri)=y\pri$. It is enough to show that
there exists $a\in\widetilde{B_{<\a}}$ such that $\ha(y\pri)=1$ and $\ha\inv(0)$ is
a neighbourhood of $y$. If $\o\ne\o\pri$ then, since $A$ is regular
there is an $a\in A$ with $\ha(\o\pri)=1$ and $\ha\inv(0)\in
\Nbb {\OA} \o$. Now, if $a$ is considered as an element $a_{\BA}$ of $\BA$ we have
$\widehat{a_{\BA}}=\ha\circ\pi$  where $\pi$ is the usual projection
$\Om\BA\to\O\;;\;(\ph,\m)\mapsto\ph$.
So $\widehat{a_{\BA}}(y\pri)=1$ and $\widehat{a_{\BA}}\inv(0)
=\fib{ \ha\inv(0) }\in\Nbb{\O(B_{<\a})} y$.

Now suppose that $\o=\o\pri$ so that $\l_\b\ne\l\pri_\b$ for some $\b<\a$.
Let $z=\Ps\b\a(y)$ and $z\pri=\Ps\b\a(y\pri)$. Then $z\ne z\pri$ so by
hypothesis there exists $b\in\tBb$ such that $\hat b(z\pri)=1$ and
$\hat b\inv(0)\in\Nbb\OBb z$.

We may realise $\tBb$ as the closure in $\widetilde\BA$ of $B_\b$ and
therefore regard $b$ as an element $a:=b_{\widetilde\BA}\in\widetilde\BA$
whose Gelfand transform
is a map $\ha=\hat b\circ \Ps\b\a\map{\Om{\widetilde\BA}=\Om\BA}\C$.
So again $\ha(y\pri)=1$ and $\ha\inv(0)\in\Nbb{\O(\widetilde\BA)}y$.

Therefore the completion of $\BA$ is regular and $\a\in{\cal J}$. By the transfinite
induction theorem, ${\cal J}=\U$. This completes the proof because
$B_\U=B_{\a_1}$.\eop

\cor 2.2.1.2. Let $((A_\s),(\tss\s\t)_{\s\le\t}\st\s,\t\le\u)$
be a system of Narmania or Cole extensions
of the regular Banach algebra $A$. Then $A_\tau$ is regular for all $\tau\le\u$.

\pf  The proof is routine and has been given in special cases in $\ref{FeiNTSR}$
and $\ref{Kar}$. We supply the details for the benefit of the reader.

Let $\J=\set{\t\le\u\st A_\t\hbox{ is regular}}$, $\t\le\u$, and suppose that
$[0,\t)\subseteq\J$. By assumption $0\in\J$. If $\t=\s+1$ for some $\s<\t$ then
by  Proposition 1.3.2.1, the Narmania extension $A_{\s, \U_\s}$ equals $\widetilde{B_{\s,\U_\s}}$
and
is regular by the proposition above. If $A_\s$ is a uniform algebra then the Cole extension
$A_\s^{\U_\s}$ is regular since, by Proposition 1.3.2.1 again, $A_\s^{\U_\s}$
can be identified with the sup-norm completion of the regular algebra $A_{\s, \U_\s}$.

It remains to check the case where $\t$ is a non-zero limit-ordinal. To do this
we show that the dense subalgebra $D=\cup_{\s<\t}\tss\s\t(A_\s)\subseteq A_\t$ contains enough
functions to satisfy the regularity condition.

Now, as in Proposition 1.3.2.2,
$$(\O_\t,(\Ps\s\t)_{\s<\t})=\il{} \left((\O_\s),(\Ps\r\s)\st\r\le\s<\t\right)$$
where $\Ps\r\s=\th^*_{\r,\s}\colon\o\mapsto\o\circ\tss\r\s$ for all $\r\le\s\le\u$.

Again it is convenient to work with the concrete space
$$\O_\t=\set{\k\in\times_{\s<\t}\O_\s\st\hbox{ for all }\r\le \s<\t\quad \k_\r=\Ps \r\s(\k_\s)}$$
so that $\Ps\s\t$ is the restriction to $\O_\t$ of the coordinate map for all $\s<\t$.

Let $y$ and $y\pri$ be distinct points of $\O_\t$. Recall that $E\int$
denotes the interior of a subset $E$ of a topological space. We show that there exists $b\in D$ with
$$y\in\hat b\inv(0)\int\qquad\hbox{ and }\qquad\hat b(y\pri)=1.$$
Since $y\ne y\pri$ there exists $\s<\t$ such that $y_\s\ne y\pri_\s$. By hypothesis
there exists $b_\s\in A_\s$ with $y_\s\in\widehat{b_\s}\inv(0)\int$ and $\widehat{b_\s}(y\pri_\s)=1$.

Let $b=\tss\s\t(b_\s)$. The following calculation shows that $\hat b=\widehat{b_\s}\circ\Ps\s\t$:
for all $\k\in\O_\t$,
$$\eqalign{\hat b(\k)&=\tss\s\t(b_\s)\htt(\k)\cr
&=\k(\tss\s\t(b_\s))=(\k\circ\tss\s\t)(b_\s)\cr
&=\Ps\s\t(\k)(b_\s)\cr
&=\widehat{b_\s}(\Ps\s\t(\k))\cr
&=(\widehat{b_\s}\circ\Ps\s\t)(\k).\cr}$$
Therefore $\hat b\inv(0)=\Ps\s\t\inv\left(\widehat{b_\s}\inv(0)\right)\supseteq U$
where $U:= \Ps\s\t\inv\left(\widehat{b_\s}\inv(0)\int\right)$ is open.
We have $y\in U$ since
$\Ps\s\t(y)=y_\s\in\widehat{b_\s}\inv(0)\int$. Finally,
$\hat b(y\pri)=\left(\widehat{b_\s}\circ\Ps\s\t\right)(y\pri)=\widehat{b_\s}(y\pri_\s)=1$.
\eop

\vskip 15pt
\noindent {\bf 2.2.2. Localness}\vskip 10pt

\noindent We remind the reader about a notion which is related to regularity. It is
of interest in the theory of uniform algebras.

\dfn 2.2.2.1. Let $A$ be a unital normed algebra.
The function $f\map\O\C$ is said to {\it belong locally to $A$ at $\o\in\O$}
if there exists $a\in A$ and a neighbourhood, $V$, of $\o$, such that
$f\rest V=\ha\rest V$. We say that $f$ is {\it locally in $A$} if
$f$ belongs locally to $A$ at every point of $\O$.
The normed algebra $A$ is said to be {\it local} if $\hA$ contains every function on $\O$
which belongs locally to $A$.

\npar The first example of a non-local
uniform algebra
was given in $\ref{Kal}$. It is standard that if $A$ is a regular normed algebra
then $A$ is local (for a proof of this,
see $\ref{Pal}$, Lemma
7.2.8). We give an example of a regular normed algebra whose Arens-Hoffman extension, $\Aa$,
is not local.

\ex 2.2.2.2. Let $A$ be the algebra of continuous functions on $I=[0,1]$ which are
piecewise polynomial. To be explicit, for each $p\in A$ there exists a partition
$(0=t_0<t_1< \ldots< t_n=1)$ of $I$ and polynomials $(p_1,\ldots,p_n)$ such that
$$p(t)=p_i(t)\qquad\qquad t\in [t_{i-1}, t_i],\qquad i=1,\ldots,n.$$

The algebra $A$ is dense in $C(I)$ in the sup-norm by the Stone-Weierstrass theorem.
Therefore its space of continuous characters, $\O$, can be identified with the
character space of $C(I)$, which is known to be $I$ (by
the example after Definition 1.1.3.1). The algebra is clearly
regular (on $I$) since it contains all piecewise-linear, continuous
maps.

Let $a_0(t)=t\ (t\in I)$ so that $\a(x)=x^2-a_0\in A[x]$. Let $\Aa$ denote the Arens-Hoffman
extension of $A$ with respect to $\a(x)$ for some choice of norm parameter (for
example, 1). The continuous-character space of $\Aa$ is
$\Oa=\set{(t,\l)\in I\times\C\st\l^2=t}$
which is homeomorphic to $J:=[-1,1]$ under the map
$$\ph\mapto J\Oa\l{(\l^2,\l)}.$$
Let $E=\ph([0,1])$ and suppose $b\in\Aa$ satisfies $\hat b(E)=\set0$.
There exist $b_0,b_1\in A$ such that $b=b_0+b_1\xb$. By repeating some of the polynomials if
necessary we can assume that $b_0$ and $b_1$ are defined by the polynomials
$(p_1,\ldots,p_n)$ and $(q_1,\ldots,q_n)$ respectively and the common partition
$(t_0, \ldots, t_n)$ of $I$.

By assumption we have for all $\l\in [0,1]$ that
$$(\hat b\circ\ph)(\l)=\widehat{b_0}(\l^2) +\l \widehat{b_1}(\l^2)=0.$$
Consider a subinterval $(t_{i-1} , t_i )$ for $i\in\set{1,\ldots,n}$. Then for all
$\l\in (\sqrt{t_{i-1}} , \sqrt{t_{i}} )$,
$$p_i(\l^2)+\l q_i(\l^2)=\hat b(\ph(\l))=0.$$
But the functions $\l\mapsto1,\l,\l^2,\ldots$ are linearly independent on
$[\sqrt{t_{i-1}} , \sqrt{t_{i}} ]$ so $p_i=q_i=0$. Since $i$ was arbitrary,
$b_0=b_1=0$. Therefore $b=0$. So $\hat b(\ph([0,1]))=\set0$ implies
that $b=0$.

We now show that $\Aa$ is not local. Consider
$$g\mapto J\C\l{\cases{({4\over3} )(\l^2 - {1\over4}) & $\l\in [-1,-{1\over2}]$\cr
	0  & $\l\in [-{1\over2}, 1]$\cr} }$$
Then $\gt:=g\circ\ph\inv\in C(\Oa)$. We show that $\gt$ belongs locally to $\Aa$. We
know that $\gt$ does not belong to $\Aa$ because $\gt(
\ph([0,1]))=\set0$ and $\gt\ne0$.

Let $U_1=\ph([-1,0))$ and $U_2=\ph((-{1\over2}, 1])$. Then $\set{U_1,U_2}$ is
an open cover of $\Oa$. Also $\gt\vert_{U_2}$ is constant and so
belongs to $\widehat\Aa\vert_{U_2}$.

Let $$a\mapto I\C t{\cases{0 & $t\in [0,{1\over4}]$\cr
	({4\over3} )(t - {1\over4})  & $t\in [{1\over4}, 1]$\cr} }$$
Then $a\in A$ and for all $\l\in[-1,0)$,
$$\eqalign{(a+0\xb)\htt(\ph(\l))=a(\l^2)&=
\cases{({4\over3} )(\l^2 - {1\over4}) & $\l\in [-1,-{1\over2}]$\cr
	0  & $\l\in [-{1\over2}, 0)$\cr}\cr
	&= \gt(\ph(\l)).\cr}$$
So $\gt\vert_{U_1}=\ha\vert_{U_1}$.\eop

\vskip 10pt
\noindent We now combine methods of Cole and Sidney (see $\ref{Col}$
and $\ref{Sid}$ respectively) to produce an extension of a given
uniform algebra, $(A, X)$, which is local and integrally closed. Of course,
it is easy to produce a trivial solution to this: we can take a Cole
extension of $C(X)$ which is an integrally closed algebra, $C(X)_\oi$, and,
by Proposition 1.4.2.3, $C(X)_\oi=C(X_\oi)$, which is local. However in the
construction below, the \v Silov boundary of the extension will be proper if
the same is true of $A$.

\thm 2.2.2.3. Let $(A, X)$ be a uniform algebra. Then there exists an
integrally closed, local, uniform algebra which is an isometric extension
of $A$ and which has a proper \v Silov boundary if the \v Silov boundary
of $A$ is proper.

\pf First let us fix some notation. For a uniform algebra, $B$, let
$L(B)\subseteq C(\Om B)$ denote the closure of the algebra of the functions locally
in $B$. Let $\U(B)$ denote the set of monic polynomials over $B$.

By first passing to $\hA$ (this preserves the \v Silov boundary of
course; see $\ref{Lei}$, p. 52) if necessary, we may assume that $A$ is natural.

Let $\ui>0$ be an ordinal, $\u\le\ui$, and suppose we have chosen natural uniform algebras,
$(A_\t, X_\t)_{\t<\u}$, and continuous, open surjections $(\Ps\s\t\map{X_\t}{X_\s})_{\s\le\t<\u}$ such
that
\item{(i)} $(A_0,X_0)=(A,X)$;
\item{(ii)} if $\t=\s+1<\u$ for some $\s<\t$ then $A_\t=L(A_\s)^{\U(A_\s)}$,
the Cole extension of $L(A_\s)$ generated by $\U(A_\s)$, the set of all monic polynomials over $A_\s
\subseteq L(A_\s)$;
\item{(iii)} if $0<\t<\u$ and $\t$ is a limit ordinal then $((X_\s)_{\s<\t}, (\Ps\r\s)_{\r\le\s<\t})$
is an inverse system with inverse limit $(X_\t,(\Ps\s\t)_{\s<\t})$ and $(A_\t,X_\t)$ is the
direct limit (see Appendix 2) of the direct system
$\bigl( (A_\s)_{\s<\t}, (\pss\r\s\colon f\mapsto f\circ\Ps\r\s)_{\r\le\s<\t} \bigr)$;
\item{(iv)} for all $\t<\u$, $\S{A_\t}=\fibp0\t{\S A}$.

By the transfinite recursion theorem we can obtain a
system satisfying the above conditions with $\u=\ui$ if we can replace
`$<\u$' everywhere by `$\le\u$'. (This is obvious if $\u=0$.)

First suppose that $\u=\t+1$ for some $\t<\u$. Then set $A_\u=L(A_\t)^{\U(A_\t)}$.
By a result of Rickart ($\ref{RicMIS}$) we have $\Om{L(A_\t)}=\Om{A_\t}=X_\t$.
So, since $L(A_\t)$ is natural, the Cole extension is too by Proposition 1.3.1.3.
Let $\p\t\u$ be the usual
projection $X_\u\to X_\t$ where $X_\u$ is the space, $X_\t^{\U(A_\t)}$,
underlying $A_\u$, and define $\p\s\u=\p\s\t\circ\p\t\u$ for
$\s\le\t$. All these maps are continuous, open surjections by Lemma 1.3.2.2 and
the inductive hypothesis.

We have $\S{A_\u}=\fibp\t\u{\S{L(A_\t)}}$ by Proposition 1.5.2.2. It follows from a remark
in $\ref{StoMIS}$, $\S{L(A_\t)}=\S{A_\t}$ and so by hypothesis
$$\S{A_\u}=\fibp\t\u{ \fibp0\t{ \S A} }
=(\p0\t\circ\p\t\u)^{-1}(\S A)=\fibp0\u{\S A}.$$

Now suppose that $\u$ is a non-zero limit ordinal. Let
$$\eqalign{(X_\u, (\p\t\u)_{\t\le\u})&=
\il{}\bigl( (X_\t)_{\t<\u}, (\p\s\t)_{\s\le\t<\u} \bigr),\qquad\hbox{ and}\cr
(A_\u, (\pss\t\u)_{\t\le\u})&=
\dl{}\bigl( (A_\t)_{\t<\u}, (\pss\s\t)_{\s\le\t<\u} \bigr)\cr}$$
(see Appendix 2). By Proposition 1.3.1.2, $A_\u$ is again natural.
The proof of Lemma 1.3.2.2 shows that the maps $(\Ps\t\u)_{\t\le\u}$
are open surjections. The same arguments as
in Proposition 1.5.2.2 show that $\S{A_\u}=\fibp0\u{\S A}$.

Therefore we can obtain a family $(A_\u, X_\u)_{\u\le\ui}$ satisfying
conditions (i)-(iv) above where $\ui=\o_1$, the first uncountable
ordinal.

A similar argument to the one used in Section 1.2.4
shows that $$A_{\o_1}=
\bigcup_{\u<\o_1}\pss\u{\o_1}(A_\u)$$ and that $A_\oi$ is integrally closed.

To prove that $A_{\o_1}$ is local let $f\in C(X_{\o_1})$ belong locally to
$A_{\o_1}$. It has been noted in Proposition 1.4.2.3 that
$$C(X_\oi)=\dl{}\bigl( (C(X_\u))_{\u<\oi}, (\pss\t\u)_{\t\le\u<\oi} \bigr)
=\bigcup_{\u<\oi}\pss\u\oi(C(X_\u)).$$ Hence
$f\in\pss{{\u_0}}\oi(C(X_{\u_0}))$
for some $\u_0<\oi$.

Let $(V_1,\ldots,V_n)$ be an open cover of $X_{\o_1}$ and
$f_1,\ldots,f_n\in A_\oi=\bigcup_{\u<\o_1}\pss\u{\o_1}(A_\u)$ with
$$f\vert_{V_i}=f_i\vert_{V_i}\qquad(i=1,\ldots,n).$$

Since $\bigl( \pss\u{\o_1}(A_\u) \bigr)_{\u<\oi}$ is an
increasing family of algebras we can take $\u<\oi$ with
$f_1,\ldots,f_n\in\pss\u\oi(A_\u)$ and $f\in\pss\u\oi(C(X_\u))$. Let $g_1,\ldots,g_n
\in A_\u$ and $g\in C(X_\u)$ be such that $f=\pss\u\oi(g)$
and $f_i=\pss\u\oi(g_i)\ (i=1,\ldots,n)$.

Since $\Ps\u\oi$ is an open map,
$(U_i)_{i=1}^n$ is an open cover of $X_\u$ where $U_i=\p\u\oi(V_i)\
(i=1,\ldots,n)$. Fix $i\in\set{1,\ldots,n}$ and let $\ku\in U_i=\p\u\oi(V_i)$ and
$\k\in V_i\cap\fibp\u\oi\ku$. Then
$$\eqalign{ g(\ku)&=(g\circ\p\u\oi)(\k)=f(\k)=f_i(\k)\cr
&=(g_i\circ\p\u\oi)(\k)=g_i(\ku).\cr}$$
So $g\vert_{U_i}=g_i\vert_{U_i}$ and $g$ is locally in $A_\u$.

Therefore $g\in L(A_\u)$ and $\pss\u{\u+1}(g)\in A_{\u+1}$. Let $h=
\pss\u\oi(g)\in A_\oi$. For $\k\in X_\oi,\ \k\in V_i\ (i\in\set
{1,\ldots,n})$ we have $f(\k)=g(\ku)=h(\k)$ so $f\in A_\oi$.\eop

\def\D{{\bf D}}
\def\p{\pi}

\vfill\eject


\noindent {\bigtenbf Chapter 3}
\vskip 10pt
\noindent {\bigtenrm Algebraic Extensions and Topology}
\vskip 25pt


\noindent
The well-known interaction between topology and Banach algebras (see Appendix 1)
has borne many fruit in
the literature on algebraic extensions.
See $\ref{Gor}$, $\ref{Cou}$, $\ref{Zam}$, and $\ref{Hos}$ for further details
and surveys of this. The last two references use algebraic extensions to classify
certain classes of covering spaces of compact, Hausdorff spaces, but we shall not pursue this subject.
However, the first two papers are in the same in the same vein as the results we shall present in
this chapter.

We investigate which properties of $\O$ are shared by $\Om B$, where $B$ is an algebraic
extension of $A$.

In the latter half of the chapter
we focus on \v Cech cohomology. This has been used in $\ref{Hat}$ (among
other papers) to characterise conditions on the compact, Hausdorff space $X$ which ensure that
certain algebraic equations over $C(X)$ have solutions in that Banach algebra.

Another purpose of this chapter is to prepare the ground for Chapters 4 and 5.


\def\tf{\tilde f}

\def\tg{\tilde g}

\vskip 15pt
\noindent {\medtenbf 3.1. Dimension}

\vskip 15pt
\noindent {\bf 3.1.1. Introduction}\vskip 10pt

\noindent The notion of the dimension of a topological space generalises that of the dimension
of Euclidean space. We shall only be concerned with compact, Hausdorff spaces in which case dimension
can be characterised by the condition given in Definition 3.1.1.1. We set
$$S^n=\set{x\in\R^{n+1}\st\sum_{j=1}^{n+1}x_j^2=1}\qquad(n\in\No).$$

\dfn 3.1.1.1 (see $\ref{Pea}$, Theorem 3.2.2). Let $n\in\No$ and $X$ be a compact, Hausdorff
space. We say
that $\dim X\le n$ if and only if
for each closed set $E\subseteq X$ and continuous map $f_0\map E{S^n}$, there exists
a continuous map $f\map X{S^n}$ such that $f\rest E=f_0$. If the statement is false for
all $n\in\No$ then we write $\dim X=\infty$. For $n\in\No$, we say that $\dim X=n$ if
$\dim X\le n$ but it is false that $\dim X\le n-1$.

\npar Dimension is interesting for us because, by Proposition 3.3.2 of $\ref{Pea}$,
the condition $\dim X\le 1$ is equivalent to $C(X)$ having a dense invertible group, the
subject of the next chapter.

The dimension of $X$ is also related to algebraic equations: it has been characterised (see
$\ref{Pea}$, Chapter 10) entirely in
terms of the solvability of algebraic equations over subalgebras of the real Banach algebra $C(X,\R)$, the
continuous, real-valued functions on $X$.

Furthermore, dimension, together with
\v Cech cohomology (see the next section in this chapter), was used in $\ref{Hat}$
to characterise all locally-connected, compact, Hausdorff spaces $X$ for which the following holds:
$$\hbox{ for all }f\in C(X)\ \hbox{there exists a root in }C(X)\hbox{ of }x^2-f.$$

\vskip 15pt
\noindent {\bf 3.1.2. Dimension and Algebraic Extensions}\vskip 10pt

\noindent  We establish the simple result, some cases of which are
proved in $\ref{Kar}$, that forming
algebraic extensions of a normed algebra $A$ cannot lead to a normed algebra whose
continuous-character space has dimension greater that $\dim\O$. Note that when
$\t$ is a non-zero limit ordinal and $\O_\t$ is the maximal ideal space
of the $\t$th algebra in a system,
$((A_\s)_{\s\le\u}, (\tss\r\s)_{\r\le\s\le\u})$, of algebraic extensions
of $A$ then $\O_\t=\lim_{\leftharpoondown}(\O_\s)_{\s<\t}$.
To see this,
the same argument as is used in the proof of Proposition 1.3.2.2 shows that
$$\cases{ \O_\t\to\left\{ \k\in\times_{\s<\t}\O_\s\st\hbox{ for all }\r\le\s<\t,
\ \k_\r=\k_\s\circ\tss\r\s\right\}&\cr
\o\mapsto(\o\circ\tss\s\t)_{\s<\t}&\cr}$$
is a homeomorphism.

\prop 3.1.2.1 (from $\ref{Kar}$, Theorem 4).
Let $\O$ be the continuous-character space of $A$ and $\u$ be an ordinal number. Let
$\O_\t$ be the continuous-character space of the $\t$th algebra in a system,
$((A_\s)_{\s\le\u}, (\tss\r\s)_{\r\le\s\le\u})$, of algebraic extensions
generated by $A$. Then $\dim \O_\t\le\dim \O$.

\pf (This is a rewording, in greater generality, of Karahanjan's proof.)
The result is obvious if $\dim \O=\infty$ so we may assume that
$\dim \O=n\in\No$.
Let $\J=\set{\t\le\u\st \dim \O_\t\le n}$ and
suppose that $[0,\t)\subseteq\J$. Clearly $0\in\J$.

Suppose that $\t$ is a non-zero limit-ordinal. Then, by Corollary $8.1.7$ in $\ref{Pea}$, since
$\O_\t$ is the inverse limit of $(\O_\s)_{\s<\t}$ it has dimension not more
that $\sup_{\s<\t}\dim \O_\s\le n$.

Now suppose that $\t=\s+1$ for some $\s<\t$. Then $A_\t=E(A_\s,\U)$,
the algebraic extension of $A_\s$ obtained by adjoining roots of each
element of $\U$, a set of
monic polynomials over $A_\s$. By Proposition 1.3.2.1, $A_\t$ can be realised as a
direct limit of normed algebras $E(A_\s,U)$ obtained by
applying the appropriate construction to $A_\s$ for finite
subsets, $U$, of $\U$. By Lemma 9.2.16 of $\ref{Pea}$, we have that if $X$ and $Y$ are
compact, Hausdorff spaces and $f\map XY$ is a
continuous, open surjection with finite fibres then $\dim X=\dim Y$.
By Proposition 1.3.2.2 we therefore have that for all $U\in\U\fs$
$$ \dim\Om{E(A_\s,U)}=\dim\O_\s\le n.$$
Then by
the result used in the previous paragraph we again have
that $$\dim \O_\t=\dim{\lim_{\leftharpoondown
}} \bigl(\Om{E(A_\s,U)}\st{U\in\U\fs}\bigr)\le n.$$
Therefore $\t\in\J$ in all cases and by the transfinite induction theorem
$\J=[0,\u]$.\eop

\npar It would be interesting to know if the change from
$\O$ to $\Ou$ did not lower dimension. For then we
would be able to construct topological spaces of arbitrarily high dimension with the
properties (such as `hereditary unicoherence'; see $\ref{Cou}$) known to be
enjoyed by spaces underlying certain algebraic extensions.


\def\tg{{\tilde g}}
\def\uis{{\U\in S\fs}}
\def\pp#1#2{\pi_{#1,#2}}
\def\tpp#1#2{\tilde\pi_{#1,#2}}
\def\pps#1#2{\pp#1#2^*}
\def\aiu{{\a\in U}}

\def\tr#1{\theta_{#1}}
\def\d{\delta}

\def\dl{{\lim_\rightharpoondown}\mathstrut}	
\def\il{{\lim_\leftharpoondown}\mathstrut}		

\vskip 15pt
\noindent {\medtenbf 3.2. \v Cech Cohomology}

\vskip 15pt
\noindent
Every compact, Hausdorff space, $X$, gives rise to a sequence, $(\Hn n X)_{n\in\No}$, of abelian groups
called
the {\it \v Cech cohomology groups} of $X$.
An excellent introduction to \v Cech cohomology is $\ref{Tay}$.

Let $G(A)$ denote the group of invertible elements of our generic algebra, $A$, and
$e^A$ denote the subgroup $\set{b\in A\st\hbox{there exists }
a\in A\hbox{ such that }b=e^a}$.

It is well-known (see, for example, $\ref{Tay}$, p. 137) that
$\Ho X\cong G(C(X))\bigm/ e^{C(X)}$.
If $\O$ is the maximal ideal space of a Banach algebra, $A$, then
the celebrated Arens-Royden theorem asserts that
$$G(A)/e^A\to G(C(\O))/e^{C(\O)}\;;\;e^Aa\mapsto e^{C(\O)}\ha$$
is an isomorphism (see $\ref{Pal}$, Theorem 3.5.19).


\vskip 15pt
\noindent {\bf 3.2.1. Introduction}\vskip 10pt

\noindent The following example shows that the
the condition $\Hn1 X=0$ need not be preserved when $X$ changes from $\O$ to $\Oa$,
where $\Oa$ denotes the space of closed, maximal ideals of an Arens-Hoffman extension,
$\Aa$, of $A$. Let
$A=C(I)$ and $\a(x)=(x-e^{\p \ri t})(x-e^{-\p \ri t})$. Then
$$I_\a=\set{ (t,\l)\in I\times\C\st \l=e^{\pm\p \ri t}\,},$$
and the map $\f\mapto{I_\a}{S^1}{(t,\l)}\l$ is a homeomorphism.
It is well-known that every continuous function
on $I$ with no zeros has a continuous logarithm ($\ref{Lei}$, p. 198) while
$\id{S^1}\in G(C(S^1))$ has not (see, for example, $\ref{Rao}$, p. 28).

However, we shall prove, in Section 3.2.3, that the converse is true.
That is, if $\O_\u$ is the closed, maximal ideal space
of the $\u$th algebra in a system of algebraic extensions then
$$\Hn 1\Ou=0\quad\hbox{implies}\quad\Hn 1\O=0.\eqno(1)$$
The result has
interesting consequences for Chapter 5.

It is not easy to work with the
purely topological definition of \v Cech cohomology groups and our methods make strong use of the
Arens-Royden theorem. Along the way we obtain an interesting result, Lemma 3.2.3.2, about
the \v Cech cohomology of general direct limits of Banach algebras.

\vskip 15pt
\noindent {\bf 3.2.2. Logarithms in Algebraic Extensions}\vskip 10pt

\noindent We now give the detail behind $(1)$. The key to doing
this is to establish the result for simple Arens-Hoffman extensions, in the next lemma,
and then to extend the result to direct limits (Lemma 3.2.3.2).

\lem 3.2.2.1. Let $A$ be a commutative, unital Banach algebra and $\a(x)$
be a monic polynomial over $A$. Let $\Oa$ denote the maximal
ideal space of $\Aa$, the Arens-Hoffman extension of $A$ by $\a(x)$. Then
\item{(i)} if $a\in A$ and $a\in e^{\Aa}$ then $a\in e^A$, and
\item{(ii)} $H^1(\Oa,\Z)=0$ implies that $H^1(\O,\Z)=0$.

\pf We prove (i); the second statement quickly follows from this.
Note that since $\Oa$ is also the maximal ideal space of the Arens-Hoffman
extension of $C(\O)$ by $\hat\a(x)=\widehat{a_0}+\cdots+\widehat{a_{n-1}}x^{n-1}+x^n$
we might as well assume that $A=C(\O)$. (This avoids using the \v Silov idempotent
theorem in the following argument.)

Suppose that $f\in A=C(\O)$ is invertible and that there exists $h\in C(\O)_\a$ such that
$f=e^h$.

Let $g=T(h)$ where $T$ is the operator $\Aa\to A$ of Section 1.4.2. Then $g\in A$ and
$$g(\k)={1\over n}\sum_{k=1}^n h(\k, \l_k(\k) )\qquad\qquad(\k\in\O)$$
where $\l_1(\k),\ldots,\l_n(\k)$ are the roots of $\eh\k(\a)(x)$ repeated
by multiplicity. So for each $\k\in\O$,
$$\eqalign{\left( e^{g(\k)}\right)^n&=\exp\left(\sum_{k=1}^n h(\k, \l_k(\k) \right)\cr
&=\prod_{k=1}^n e^{h(\k, \l_k(\k) }=f(\k)^n.\cr}$$
So for each $\k\in\O$ there exists a unique $m(\k)\in\set{0,\ldots,n-1}$ such that
$f(\k)=\zeta^{m(\k)}e^{g(\k)}$ where $\zeta=e^{2\pi\ri/n}$.

Now $\eta=e^{-g}f\in A$ so the sets $\O_k=\eta^{-1}(\zeta^k)$ are open, disjoint in, and
cover $\O$. Hence there exists $v\in A$ such that $v\vert_{\O_k}=2\pi\ri k/n\ (k=0,\ldots,
n-1)$. Put $\tg=g+v$. Then for all $\k\in\O$, if $\k\in\O_k$,
$$e^{\tg(\k)}=e^{g(\k)}\zeta^k=e^{g(\k)}\eta(\k)=f(\k).$$
So $f=e^\tg$ as required.\eop

\vskip 15pt
\noindent {\bf 3.2.3. Logarithms in Direct Limits of Normed Algebras}\vskip 10pt

\noindent In this section, $A$ need not be commutative.
Recall that if
the unital Banach algebra, $A$, is not
commutative then $G_1(A)$, the component of $G(A)$ containing the identity,
is the normal subgroup of
all finite products of elements of $e^A$. For other facts about
$G_1(A)$ see Chapter 10 of $\ref{RudFA}$. In particular we shall need the following
elementary observation.

\lem 3.2.3.1. Let $A$ be a unital Banach algebra and $x\in G(A)-G_1(A)$. Then
for all $y\in G_1(A)$ we have $\norm{x-y}\ge\norm{x^{-1}}^{-1}$.

\pf It is enough to show that $B_A(1,1)\subseteq G_1(A)$ for if this holds
and $\norm{x-y}<\norm{x^{-1}}^{-1}$ for $x\in G(A),y\in G_1(A)$ then we have
$$\norm{1-yx^{-1}}=\norm{(x-y)x^{-1}}<1.$$
From this we have $yx^{-1}\in G_1(A)$ so that $y\in G_1(A)x$ and $x\in G_1(A)$,
a contradiction.

We have the elementary
fact that $B_A(1,1)$ is contained in the path component of the identity, $1$, of $A$
for, if $x\in B_A(1,1)$ then
$$F\mapto IAt{(1-t)x+t1}$$
is a path in $A$ connecting $x$ and $1$. Moreover, $\norm{1-F(t)}<1$ for all
$t\in I$ so $\im F\subseteq G(A)$.\eop

\npar For the standard facts and constructions of direct limits we refer the reader
to Appendix 2.

\lem 3.2.3.2. Let $(A, (\tr i)_{i\in J})$ be the direct limit of a system
of unital Banach algebras and unital homomorphisms of norm 1,
$$(( A_i)_{i\in J},(\tss ij)_{i\le j}).$$
Fix $i\in J$ and suppose that for all $j\ge i$,
$$\tss ij(G(A_i)-G_1(A_i))\subseteq G(A_j)-G_1(A_j)\eqno(*)$$
Then $$\tr i(G(A_i)-G_1(A_i))\subseteq G(A)-G_1(A).$$

\pf  Suppose $x_i\in G(A_i)-G_1(A_i)$.
By ($*$) we have that for all $j\ge i$, $\tss ij(x_i)\not\in G_1(A_j)$. Thus for
all $j\ge i$ and $y_j\in G_1(A_j)$ we have $\norm{y_j-\tss ij(x_i)}\ge
\norm{\tss ij(x_i)^{-1}}^{-1}.$
Also, $\norm{\tss ij(x_i)^{-1}}=\norm{\tss ij(x_i^{-1})}$ and
$\norm{\tss ij(x_i^{-1})}\le \norm{x_i^{-1}}$ so we have that
$$\hbox{for all }j\ge i,\ y_j\in G_1(A_j)\qquad\norm{y_j-\tss ij(x_i)}\ge\norm{x_i\inv}^{-1}
\eqno(\dag)$$
Suppose $\tr i(x_i)\in G_1(A)$. There exist $c_1,\ldots, c_n\in A$ such that $\tr i(x_i)
=e^{c_1}\cdots e^{c_n}$. The map $A^n\to A\;;\;(b_1,\ldots,b_n)\mapsto
e^{b_1}\cdots e^{b_n}$ is continuous and so there exists $\d>0$ such that
for all $d_1,\ldots, d_n\in A$, $\norm{ e^{d_1}\cdots e^{d_n}-
\tr i(x_i) }<\e:= \norm{x_i\inv}^{-1}$ if $\sum_{k=1}^n\norm{d_k-c_k}<\d$.
Now $\cup_{j\in J}\tr j(A_j)$ is dense in $A$ so there exists $j\ge i$ and
$b_1,\ldots, b_n\in A_j$ such that $\sum_{k=1}^n\norm{\tr j(b_k)-c_k}<\d$.
Note that, since $\tr j$ is continuous, we have $y_j:=\tr j(w_j)
=e^{\tr j(b_1)}\cdots e^{\tr j(b_n)}\in G_1(A)$ where $w_j=e^{b_1}\cdots e^{b_n}\in
G_1(A_j)$.

By the construction of the normed direct limit,
$$\eqalign{\norm{y_j-\tr i(x_i)}&=\norm{\tr j(w_j-\tss ij(x_i))}\cr
&=\limsup_{k\ge j}\norm{ \tss jk(w_j-\tss ij(x_i))}<\e.\cr}$$
Thus there exists $k\ge j$ such that
$\sup_{j\pri\ge k}\norm{ \tss j{j\pri}(w_j-\tss ij(x_i))}<\e.$ Therefore
$$\eqalign{ \norm{ \tss jk(w_j-\tss ij(x_i))}&=\norm{\tss jk(w_j)-\tss ik(x_i)}\cr
&=\norm{ e^{\tss jk(b_1)}\cdots e^{\tss jk(b_n)}-\tss ik(x_i)}<\norm{x_i\inv}^{-1},\cr}$$
and this contradicts (\dag).\eop

\npar We remark that, in the category of commutative, unital Banach algebras and isometric
isomorphisms, the lemma above follows from the Arens-Royden theorem and the
`continuity of \v Cech cohomology'.

The following special case
of Lemma 2 helps in Corollary 5.

\cor 3.2.3.3. Suppose the commutative, unital Banach algebra $A$ has a dense subalgebra
$\cup_{i\in J} A_i$ where $J$ is a well-ordered set and for each $i\in J$,
$A_i$ is a unital, Banach subalgebra of $A$. Let $0$ be the minimum of $J$.
Then, provided that
$$G(A_0)-e^{A_0}\subseteq G(A_0)-\bigcup_{j\ge0} e^{A_j},$$
we have
$$G(A_0)-e^{A_0}\subseteq G(A)-e^A,$$
and so if $\Ho{\Om A}=0$ then $\Ho{\Om{A_0}}=0$.\eop

\noindent The next two results could be easily be extended to Narmania extensions, using the Arens-Royden
theorem, but in this thesis our only interest is in their consequences for uniform algebras.

\cor 3.2.3.4. Let $A$ be a uniform algebra with maximal ideal space $\O$. Let
$\UA$ denote the Cole extension generated by a
set of monic polynomials, $\U$, over $A$. Then $\Ho\OU=0$ implies that
$\Ho\O=0$ where $\OU$ is the maximal ideal space of $\UA$. More generally,
$$\ps\left(G(A)-e^A\right)\subseteq G(\UA)-e^{\UA}.$$

\pf We express $A^\U$ as a direct limit of Cole extensions by finite subsets
of $\U$ to obtain the result. In order to do this we need to make some
preliminary observations. For all $S\subseteq\U$, let $(A^S,\O_S)$ denote
the Cole extension of $A$ with respect to $S$ and let the canonical maps
be denoted as follows:
$$\eqalign{&p_{\a,S}\mapto{\O_S}\O\ol{\l_\a}\qquad\qquad(\a\in S);\cr
&\pi_S\mapto{\O_S}\O{\ol}\o;\cr
&\ps_S\mapto{C(\O)}{C(\O_S)}f{f\circ\pi_S}.\cr}$$
Notice that for all $U\subseteq V\subseteq \U$ there is a continuous surjection
$$\Ps UV\mapto{\O_V}{\O_U}\ol{(\o,(\l_\a)_{\a\in U}).$$
It induces an isometric monomorphism
$$\pss UV\mapto{C(\O_U)}{C(\O_V)}f{f\circ\Ps UV}.}$$
It is easy to check that $\pss UV(p_{\a,U})=p_{\a,V}$ and that $\pss UV\circ\ps_U=\ps_V$.
Since $A^U$ is generated by $\ps_U(A)\cup\set{p_{\a,U}\st \a\in U}$ we have
$\pss UV(A^U)\subseteq A^V$.

It is clear that $D:=\bigcup_{U\in\U\fs}\pss U\U(A^U)$ is dense
in $A^\U$. By Corollary 3 it is enough to show that for all $a\in G(A)-e^A$ and
$U\in\U\fs$ we have $\ps(a)=\pss\emptyset\U(a)\not\in e^{\pss U\U(A^U)}$.
Suppose that $b\in A^U$ is such that $\pss\emptyset\U(a)=e^{\pss U\U(b)}$. Since
$\pss U\U$ is continuous, $e^{\pss U\U(b)}={\pss U\U(e^b)}$.
We also have $\pss\emptyset\U(a)=\pss U\U(\pss \emptyset U(a))$ and since $\pss U\U$
is injective this gives $\pss\emptyset U(a)=e^b\in A^U$.

An argument similar to the one at the end of the proof of Proposition 1.5.3.7
shows that $A^U$ is isometrically isomorphic to $A$ extended finitely many
times by the Cole construction. It now follows inductively from Lemma 3.2.2.1
that $a\in e^A$, a contradiction. Since $U$ was arbitary, this
completes the proof.\eop

\cor 3.2.3.5. Let $((A_\t,\O_\t),(\pss\s\t)_{\s\le\t\le\u})$ be a system of natural Cole
extensions of $(A,\O)$ in the usual notation. Then for all $\t\in[0,\u]$,
$$\pps0\t(G(A)-e^A)\subseteq G(A_\t)-e^{A_\t}.$$
In particular, $\Ho{\O_\t}=0$ implies that $\Ho\O=0$.

\pf Let $\J=\set{\t\le\u\st \pps0\t(G(A)-e^A)\subseteq G(A_\t)-e^{A_\t}}$, and
suppose that $[0,\t)\subseteq \J$.

First suppose that $\t$ is a limit ordinal. Clearly $0\in\J$ so we can assume
$\t>0$. Then $A_\t=\dl(A_\s\st \s<\t)$ and $\bigcup_{\s<\t}\pss\s\t(A_\s)$
is dense in $A_\t$. We have
$$\eqalign{G(\pps0\t(A))-e^{\pps0\t(A)} &= \pps0\t(G(A)-e^A)\qquad\hbox{(since $\pss0\t$
is injective)} \cr
&\subseteq\bigcap_{\s<\t}\left[ \pps0\t(G(A))-\pps\s\t(e^{A_\s}) \right]\qquad
\hbox{(by inductive hypothesis)}\cr
&=\pps0\t( G(A))-\bigcup_{\s<\t}\pps\s\t(e^{A_\s})
=G(\pps0\t(A))-\bigcup_{\s<\t}e^{\pps\s\t(A_\s)}.\cr}$$
So $\t\in \J$ by Corollary 3.

Secondly suppose that $\t=\s+1$ for some $\s<\t$. By Corollary 4 we have
$$\pps\s\t(G(A_\s)-e^{A_\s})\subseteq G(A_\t)-e^{A_\t}.$$
By hypothesis $\pps0\s(G(A)-e^A)\subseteq G(A_\s)-e^{A_\s}$. Thus
$$\pps0\t(G(A)-e^A)=\pps\s\t\left( \pps0\s(G(A)-e^A) \right) \subseteq G(A_\t)-e^{A_\t},$$
so $\t\in \J$.

By the transfinite induction theorem, $\J=[0,\u]$ as required.\eop

\npar The results in this section suggest that if $B$ is an integral extension of $A$
then $\Ho\O$ vanishes whenever $\Ho{\Om B}$ does, but we do not know if this is true.

\def\na{{\bf na}}	
\def\dl#1{{\lim_\rightharpoondown}\mathstrut_{#1}}	
\def\il#1{{\lim_\leftharpoondown}\mathstrut_{#1}}		

\vfill\eject


\noindent {\bigtenbf Chapter 4}
\vskip 10pt
\noindent {\bigtenrm  The Invertible Group of an Extension of a Normed Algebra}
\vskip 25pt


\noindent We study the group of invertible elements of $A$ in this chapter. Notation
for this group was introduced in Section 3.2.

The reason for investigating $G(A)$ was that it yielded an insight into a famous open
problem of Gelfand. He asked (see $\ref{GelSCF}$) if every natural uniform algebra
on the unit interval is trivial. For a brief history of this problem see
$\ref{Rud}$, Chapter 22. We describe how our work applies to the problem
in Section 3. This leads to a natural question, which is taken
up in Chapter 5.

In the first section we collect some elementary results and examples about $G(A)$.

In the following section we pursue our programme of investigating which properties are common to
$A$ and its algebraic extensions. We show that, in $\Ba$, if $A$ has dense invertible group then so does
every algebraic extension. This allows us to simplify some arguments of Karahanjan
in the examples in $\ref{Kar}$.
We also note, in this section, two other types
of extension of $A$ to which our methods apply.


\def\tA{\widetilde{A}}
\def\sp#1{\sigma\left(#1\right)}

\def\T#1{\Theta(#1)}

\vskip 15pt
\noindent {\medtenbf 4.1. Introductory Results and Examples}

\vskip 15pt
\noindent {\bf 4.1.1. The Group of Quasi-Invertible Elements}\vskip 10pt

\noindent Many attributes of the spectra of elements of a (normed) algebra, $A$, are more naturally
expressed by means of the quasi-invertible group, $qG(A)$, of $A$, which we proceed to define.
For $a,b\in A$ let $a\circ b:=a+b-ab$. This is a binary operation on $A$ with identity
element $0$. If $a\circ b=0$ we say $a$ is a {\it left quasi-inverse} for $b$ and
$b$ is a {\it right quasi-inverse} for $a$. We set $qG(A)=\set{a\in A\st a \hbox{
has (right) quasi-inverse}}$.
For the definition
and further standard properties of $qG(A)$ we refer the reader to $\ref{Pal}$, p. 192. In particular
$(qG(A),\circ)$ is a group and
there is a homeomorphism
$$\theta\,\colon\, (\C-\set0)\times qG(A) \longrightarrow G(A_1)  \;;\;(\m,b)\longmapsto (-\m b,\m),$$
where $A_1$ denotes the standard unitisation (see Appendix 3) of $A$. We also have that
$qG(A)$ is open if $A$ is complete.

\prop 4.1.1.1. The group $G(A_1)$ is dense in $A_1$ if and only if $qG(A)$ is dense in $A$.

\pf Suppose that $qG(A)$ is dense in $A$. Let $(a,\l)\in A_1$ and let $\e>0$. Choose $r>0$
so that $r[1+(r+\ab\l)]<\e$. Take $\m\in B_\C(\l,r)-\set0$.
By assumption here exists $b\in qG(A)$ such that $\norm{b+\m^{-1}a}<r$.
Then $ (-\m b,\m) \in G(A_1) $ and
$$\norm{ (-\m b,\m) -(a,\l) }=\norm{\m b+a}+\ab{\m-\l}<\ab\m r+r\le
(r+\ab\l)r+r<\e.$$

Conversely, suppose that $\overline{G(A_1) }=A_1$. Let $b\in A$ and let $\e>0$.
If $a=0$ then $a\in qG(A)$ so
we may assume that $a\ne0$. Choose $r\in(0,1)$ such that $r+\left({r\over{1-r}}\right)\norm a
<\e$.

Now take $(a,\l)\in G(A_1)$ with $\norm{ (a,\l)-(-b,1)}=\norm{a+b}+\ab{1-\l}<r$. Then
$(\l,-\l^{-1}a)\in (\C-\set0)\times qG(A)$. We have
$$\eqalign{\norm{-\l^{-1}a-b}&=\norm{\l^{-1}a+b}=\norm{a+b+(\l^{-1}-1)a}\cr
	&\le r+\ab{{1-\l}\over\l}\norm a\cr
	&\le r+(r/\ab\l)\norm a\cr
	&< r+\left( {r\over{1-r}}\right)\norm a<\e.\cr}$$
Hence $\overline{qG(A)}=A$ as required.\eop

\npar Now suppose that $A$ is a unital normed algebra.
By Proposition A.1.6, for every closed, maximal
ideal, $M$, of $A$, the map
$$M_1\to A\;;\;(a,\l)\mapsto a+\l 1_A$$
is a topological isomorphism. Also recall that, since $A$ is unital, it has at least one
maximal ideal ($\ref{Dal}$, p. 40). The following result quickly follows from
these observations.

\cor 4.1.1.2. Let $A$ be a unital normed algebra. Then the following are
equivalent
\item{(i)} $G(A)$ is dense in $A$;
\item{(ii)} $qG(M)$ is dense in $M$ for some closed, maximal ideal, $M$, of $A$;
\item{(iii)} $qG(M)$ is dense in $M$ for every closed, maximal ideal,
$M$, of $A$.

%
%

\npar Let $A$ be a Banach algebra. If $A$ is unital then it is shown in
$\ref{Cor}$
and $\ref{Fal}$ that $\overline{G(A)}=A$ if and only if the set of
elements of $A$ whose spectra have empty interiors is dense in $A$.
We now describe an extension of this spectral characterisation
to non-unital Banach algebras. Let
$$\qs A a:=\set{\l\in\C-\{0\}\st \l^{-1}a\not\in qG(A)}\qquad(a\in A).$$
Then (see $\ref{Pal}$, p. 196) we have that, if $A$ is unital, for all $a\in A$,
$$\set0\cup \qs A a=\set0\cup\s_A(a),$$
where $\s_A(a)$ denotes the spectrum of $a$ in $A$. In particular, $\s_A(a)$
has empty interior if and only if $\qs A a$ has. Define for any Banach algebra, $A$,
$$\T A:=\set{a\in A\st \qs Aa\int=\emptyset}.$$
Thus if $A$ is unital then $\T A= \set{a\in A\st \s_A(a)\int=\emptyset}$.
We can now give the following result. The proof is a modification of the
proof in $\ref{Fal}$ for the special case of unital Banach algebras.

\prop 4.1.1.3. Let $A$ be a Banach algebra which is not necessarily unital. Then
$$\ov{qG(A)}=A\quad\hbox{if and only if}\quad\ov{\T A}=A.$$

\pf First suppose that $qG(A)$ is dense in $A$. Let $\set{\l_n\st n\in\N}$ be a
countable dense subset of $\C-\set0$. Since $qG(A)$ is open in $A$, we have
from Baire's category theorem that $T=\bigcap_{n\in\N}qG(A)\l_n$ is dense in $A$. Now
$T\subseteq\T A$ for if $b\in T$ and $\qs A b\int\ne\emptyset$ then for some
$n\in\N$ we must have $\l_n\in\qs A b$. But then $\l_n^{-1}b\not\in qG(A)$,
a contradiction.

Conversely suppose that $\T A$ is dense in $A$. Let $a\in A$ and $\e>0$. Take
$b\in\T A$ with $\norm{b-a}<\e/2$. If $b=0$ then $b\in qG(A)$ and we are done.
Otherwise let $\d>0$ be such that $1-\d>0$ and for all $\l\in B_\C(1,\d)$ we have
$\ab{1-1/\l}<{\e\over{2\norm b}}$. By assumption there exists some
$\l\in B_\C(1,\d)-\qs A b$. Since $\l\ne0$ we must have $\l^{-1}b\in qG(A)$.
Now $\norm{a-\l^{-1}b}\le\norm{a-b}+\norm{b-\l^{-1}b}<\e$.\eop

\vskip 15pt
\noindent {\bf 4.1.2. Simple Examples}\vskip 10pt

\noindent Much of the next two sections appears in $\ref{DF}$. The following
definition, which seems quite standard in the literature, is useful.

\dfn 4.1.2.1. A subalgebra, $B$, of an algebra,
$A$, is called {\it full} if we have $G(B)=B\cap G(A)$.

\npar For example, $R_0(X)$, the algebra of rational functions on a compact plane set, $X$, with
poles off $X$ and $C^n(I)$, the algebra of functions $I\to\C$ with $n$ continuous derivatives,
are both full subalgebras of their completions
with respect to the uniform norm.

\ex 4.1.2.2. Let $B$ be the restriction to $S^1$ of the
algebra of complex, rational functions with poles off $S^1\cup\set2$.
Then $B$ is easily seen to be
uniformly dense in $A=C(S^1)$ (by the
Stone-Weierstrass theorem). However, $B$ is not a full subalgebra of $A$ because, for
example, the function given by $z-2\in B$ is invertible in $A$. \eop

\npar For future reference we note the following elementary result.

\prop 4.1.2.3. Let $A$ be a Banach algebra with maximal ideal space $\O$
such that $\hA$
is dense in $C(\O)$ (for example, if $A$ is symmetric). Give $\hA$ the uniform
norm. Then $\hA$ has dense invertible
group if and only if $C(\O)$ has.

\pf It is elementary that $\hA$ is a full subalgebra of $C(\O)$. Moreover, if
$\hA$ is dense in $C(\O)$ then, since $G(C(\O))$ is open in $C(\O)$, we have from
elementary topology that
$$\overline{G(\hA)}=\overline{\hA\cap G(C(\O))}=\overline{G(C(\O))}$$
where the closures are taken in $C(\O)$. This shows that if $G(C(\O))$ is dense in
$C(\O)$ then the invertible elements of $\hA$ are dense in $\hA$. The converse
is trivial.
\eop

\npar It follows from facts mentioned in Section 1.2 of $\ref{Cor}$ that if $A$
is a commutative Banach algebra with maximal
ideal space $\O$ and $A$ is regular (see Section 2.2.1)
then the denseness of $G(A)$ in $A$ implies that $C(\O)$ also
has dense invertible group. We do not know if the converse is true for
regular Banach algebras. Neither do we know any examples of non-regular Banach
algebras for which $G(A)$ is dense in $A$ but $C(\O)$ does not have dense invertible
group.

\ex 4.1.2.4. In $\ref{Rie}$ it is shown that if $A$ is a unital $C^*$-algebra
(see $\ref{Mur}$ for this class of Banach algebras) then $A$ has
dense invertible group if and only if $M_n(A)$ has for every $n\in\N$. We remark that
the proof given in $\ref{Rie}$ applies to all Banach algebras.

\prop 4.1.2.5. Suppose that $A$ is a commutative, unital Banach algebra which is
rationally generated by a single element, $a$. (See $\ref{Dal}$, Definition 2.2.7:
we mean that the algebra
of elements obtained
by applying the rational functions with poles off $\sp a$ to $a$ is dense in $A$.) If
$\sp a$ has empty interior then $A$ has dense invertible group.

\pf This follows from the spectral mapping theorem (see $\ref{RudFA}$, Theorem 10.28)
and the fact that if $K$ is a subset of $\C$ with
empty interior and $\phi$ is a function holomorphic on a neighbourhood
of $K$ then $\phi(K)$ also has empty interior (see $\ref{Fal}$, Lemma 2.2).
Thus $A$ has a dense set
of elements whose spectra have no interior, which, as is noted in
$\ref{Fal}$, is equivalent to the condition $\bGA=A$.\eop

\npar Recall that a {\it weight sequence} on $\Z$ is a function
$\o\map\Z{(0,+\infty)}$ such that for all $m,n\in\Z$,
$\o(m+n)\le\o(m)\o(n)$. A weight function, $\o$, induces the
{\it Beurling algebra},
$$\ell^1(\Z,\o):=\set{a\in\C^\Z\st\norm a:=\sum_{n\in\Z}\ab{a_n}\o_n<+\infty},$$ under
pointwise addition, multiplication given by
$$ab:=\left(\sum_{k\in\Z}a_kb_{n-k}\right)_{n\in\Z}\qquad\qquad(a,b\in \ell^1(\Z,\o)),$$
and the norm indicated. We refer the reader to $\ref{Dal}$, p. 501, for
further properties of the Banach algebra $\ell^1(\Z,\o)$.

\cor 4.1.2.6. Let $\o$ be a weight sequence on $\Z$ and let
$A$ be the associated Beurling algebra, $\ell^1(\Z,\o)$. Then $A$ has dense invertible
group if and only if $\rho_+:=\lim_{n\to+\infty}\o_n^{1/n}=\lim_{n\to+\infty}\o_{-n}^{1/{-n}}
=:\rho_-$.
\pf It is standard (see for example $\ref{Dal}$, Theorem 4.6.7) that $A$ is
rationally generated by ${\d_1}$
where $\d_1$ is the characteristic function of $\set1$. It is also standard that
the maximal ideal space of $A$ can be identified with the annulus
$\set {w\in\C\st\rho_-\le\ab w\le \rho_+}$. By elementary results in dimension
theory, $\O(A)$ has dimension not more than 1 if $\rho_-=\rho_+$ and dimension 2 otherwise.
If $A$ has dense invertible group then so has $\hA$ and we must have
$\rho_-=\rho_+$ by Proposition 3, as $A$ is plainly symmetric.

Conversely suppose that $\rho_-=\rho_+$. Then the spectrum of $\d_1$ has
empty interior and Proposition 5 applies.\eop


\def\res{{\rm res}\,}
\def\bp#1{b^\prime_{#1}}

\def\sp#1{\sigma\left(#1\right)}

\def\tc#1{\tilde c_{n-#1}}
\def\tcb#1{\tilde c_{n-(#1)}}
\def\bpr{b^\prime}
\def\tA{\widetilde{A}}
\def\bppr{b^{\prime\prime}}

\vskip 15pt
\noindent {\medtenbf 4.2. Algebraic Extensions}

\vskip 15pt
\noindent {\bf 4.2.1. Introduction}\vskip 10pt

\noindent
As hinted at in Section 3.1.1, denseness of invertible elements in
uniform algebras seems to be connected with the existence of roots of
monic polynomials over the algebra. Accordingly we mention the following result
stated by Grigoryan in $\ref{Gri}$: a uniform algebra $A$ on a compact, Hausdorff space,
$X$, is trivial
if $C(X)$ is
finitely generated and is an integral extension of $A$

Our main result in this section is that an Arens-Hoffman extension of a Banach algebra $A$
has a dense invertible group if $A$ has.
The proof of our main result, Theorem 4.2.2.1, relies on resultants; please refer to Appendix 3 for
a discussion of these.
To set notation here, let $\a(x)$ be a monic polynomial over $A$ and $\b(x)\in A[x]$.
We shall usually fix $\a(x)$ and allow $\b(x)$ to vary. We then
denote their resultant $\res(\a(x),\b(x))$ by $R_\a(\b(x))$.

Sometimes in the literature (for example, $\ref{Col}$ and $\ref{FeiNTSR}$) the algebraic
extensions considered are obtained by adjoining square roots. In this
case we have $\a(x)=x^2-a_0$ and then $R_\a(b_0+b_1x)=b_0^2-a_0b_1^2$.

Before we embark on proving Theorem 4.2.2.1, let us check there is
not an obvious reason why it might be true. We now explain this.

Recall that the Arens-Hoffman extension $\Aa$ is the quotient of $A[x]$ by
the principal ideal $\prin$.
Arens and Hoffman prove in $\ref{Are}$ that
when the polynomials over $A$ are given the norm
$$\norm{\sum_{k=0}^n b_kx^k}=\sum_{k=0}^n\norm{b_k}t^k\qquad(b_0,\ldots,b_{n-1}\in A)$$
for sufficiently large, fixed $t>0$ then $\prin$
is closed in $A[x]$. The Arens-Hoffman norm is defined in $\ref{Are}$
to be the quotient norm on $\Aa$. Thus if $A[x]$ retained the property of having
a dense group of invertible elements then so would Arens-Hoffman extensions. We
give an example to show that this is not the case.

Let $t>0$ be fixed. Consider the following Banach algebra (a special case
of Example 2.1.18(v) in $\ref{Dal}$)
$$B:=\set{\b=\sum_{k=0}^{+\infty}b_kx^k\in A[[x]]\st
\norm\b:=\sum_{k=0}^{+\infty}\norm{b_k}t^k<+\infty},$$
where $A[[x]]$ is the algebra of formal power series and the norm is as indicated.
We can clearly regard $A[x]$ as a dense subalgebra of $B$. Thus
we can take $B$ to be the completion of $A[x]$ (see Appendix 1).
It follows from this and Lemma A.1.4
that there is a homeomorphism between the character space of $B$
and the space of continuous characters of $A[x]$.

Recall (see, for example, $\ref{Are}$) that the space of continuous characters
of $A[x]$ can be identified with $\O\times\D$ where $\D$ is the closed
unit disc centre $0$ and radius $t$. Thus the characters on $B$ are given by
$$(\o,\l)\left(\sum_{k=0}^{+\infty}b_kx^k\right)=\sum_{k=0}^{+\infty}\o(b_k)\l^k
\qquad\left((\o,\l)\in\O\times\D,\sum_{k=0}^{+\infty}b_kx^k\in B\right).$$

Now let $A=C(I)$ and $t$ be a suitable parameter. Then the character space of $B$ is
now identifiable with $I\times\D$ and therefore of dimension $3$. The algebra $B$
has some `analytic properties' as we now explain. Let $q=\sum_{k=0}^{+\infty}b_kx^k\in B$
and fix $t_0\in I$. This pair defines a continuous function on $\D$ given by
$$g(\l):=\hat q(t_0,\l)=\sum_{k=0}^{+\infty}b_k(t_0)\l^k\qquad\qquad(\l\in\D).$$
Moreover we see that $g$ is analytic on $\D\int$.

This connection with analytic function has two consequences. Firstly, note that
$B$ is not symmetric because $\bar p\not\in\hat B$ where $p=\hat x$. To see this,
suppose $q\in B$ satisfies $\hat q=\bar p$. Let $g$ be associated with $q$ as
above (with $t_0=0$ for definiteness). Then for all $\l\in\D$
we would have $g(\l)=\sum_{k=0}^{+\infty}b_k(0)\l^k=\bar\l$. This
contradicts the fact that $z\mapsto\bar z$ is
not analytic anywhere. Therefore we cannot apply Proposition
4.1.2.3 to $B$.

Secondly it is now easy to show that $x\not\in\overline{G(B)}$. If we had $x\in\overline{G(B)}$
then for every
$\e>0$ there would exist $q=\sum_{k=0}^{+\infty}b_kx^k\in G(B)$ with $\norm{q-x}<\e$. Let
$f=\id\D$ and $g(\l)=\hat q(0,\l)$. Then $g\in G(C(\D))$ and $\norm{g-f}<\e$.
Again, $g$ is analytic on $\D\int$. Hence there exists a
sequence, $(g_n)$, of invertible, analytic functions on $\D\int$ with $g_n\to g_0:=\id{\D\int}$
uniformly on compact subsets of $\D\int$. Since $g(0)=0$, this contradicts
Hurwitz' Theorem ($\ref{Dal}$, p. 797).

\vskip 15pt
\noindent {\bf 4.2.2. The Main Results}

\thm 4.2.2.1. Let $A$ be a commutative, unital Banach algebra with dense invertible
group and let $\a(x)$ be a monic polynomial over $A$. Then $G(\Aa)$ is dense in $\Aa$.
\pf
We may assume that $n$, the degree of $\a(x)$ is at least $2$,
or else $\Aa$ is isometrically isomorphic to $A$.
Let $\b(x)=b_0+b_1x+
\cdots+b_{n-1}x^{n-1}\in A[x]$ and $\e>0$. We have to show that there exists
$\tilde\b(x)\in A[x]$ with $\norm{\tilde\b(\xb)-\b(\xb)}<\e$ and $\tilde
\b(\xb)\in G(\Aa)$. In fact we shall show that by slightly perturbing
$b_0$ only, we can obtain a polynomial $\tilde\b(x)$ with $R_\a(\tilde\b(x))
\in G(A)$.

By the comments about resultants in Appendix 3,
$P(c)=R_\a(c+b_1x+\cdots+b_{n-1}x^{n-1})$ is a polynomial
$p_0+p_1c+\cdots+p_{n-1}c^{n-1}+c^n$ where $p_0,\ldots,p_{n-1}\in A$ are independent
of $c$. Consider the $n$ formal derivatives of $P$ as maps $A\to A$:
$$\eqalign{&P^{(0)}(c) = P(c);\cr
&P^\prime(c) = p_1+2p_2c+\cdots+nc^{n-1};\cr
&\qquad\vdots\cr
&P^{(n-1)}(c)=n!c.\cr}$$
Set $\tc 1=b_0$. Note that, trivially, $P^{(n-1)}$ is a local homeomorphism at
$\tc 1$.
Now let $1\le k<n$ and suppose that $\tc1,\ldots,\tc k\in A$ have been chosen so that
\item{(i)} $P^{(n-j)}$ is a local homeomorphism at $\tc j\ (1\le j\le k)$, and
\item{(ii)} $\norm{\tc j-\tcb {j-1}}<\e/n$ for $1<j\le k$.

We shall now show how to choose $\tcb{k+1}$ so that (i) and (ii) become
true with `$k$' replaced by `$k+1$'.
It is easy to see from the inverse function theorem for Banach spaces
(see, for example, Theorem 8 of Chapter 7 in $\ref{Bol}$) that $P^{(n-(k+1))}$ is a local
homeomorphism at $a\in A$ if $P^{(n-k)}(a)$ is invertible. (This fact is also stated
in $\ref{Fal}$.)

By hypothesis, $P^{(n-k)}$ is a local homeomorphism at $\tc k$ so some open
neighbourhood of $\tc k$ is mapped onto an open set in $A$. Since $G(A)$ is
dense in $A$ there is some $\tcb{k+1}\in B_A(\tc k,\e/n)$ with
$P^{(n-k)}(\tcb{k+1})\in G(A)$.

Thus $\tcb{k+1}$ has the required
properties and we can choose
$\tilde c_{n-1},\ldots,\tilde c_0\in A$
with $\norm{\tilde c_k-\tilde c_{k-1}}<\e/n\ (k=1,\ldots,n-1)$, $\tc 1=b_0$,
and $P$ a local homeomorphism at $\tilde c_0$.

Again, since $P$ is a local homeomorphism at $\tilde c_0$, we can find $\tilde b_0
\in B_A(\tilde c_0,\e/n)$ with $P(\tilde b_0)\in G(A)$. Since
$$\norm{\tilde b_0-b_0}\le\norm{\tilde b_0-\tilde c_0}+
\norm{\tilde c_0-\tilde c_1}+\cdots+\norm{\tilde c_{n-2}-\tilde c_{n-1}}<\e,$$
the result is proved.\eop

\npar
We also remark that the method of proving Theorem 1 gives another way to see
that if $A$ is a commutative, unital Banach algebra with $\bGA=A$ then for
every $n\in\N$, $M_n(A)$ also has dense invertible group. The role
of the resultant is taken by the determinant map which links the invertible
groups of $A$ and $M_n(A)$. We shall show that we can approximate
any $B=[b_{ij}]\in M_n(A)$ by an invertible matrix by perturbing only $n$ entries
$b_{i_1,j_1},\ldots, b_{i_n,j_n}$ provided that $n\mapsto i_n,j_n$ are both
permutations.

To see this, suppose that we have chosen the $n$ entries to perturb; let their indices be
$$(1, v(1)),\ldots,(n, v(n))$$ for some permutation $v\in S_n$, the
group of permutations on $\set{1,\ldots,n}$. Consider the map
$$P\mapto A A c{{\rm det}\,[b_{ij}+c\d_{v(i),j}]}$$
(where $\d_{ij}$ is the Kronecker delta). We see from the usual expansion that
$$P(c)=\sum_{\p\in S_n}\pm(b_{1,\p(1)}+c\d_{v(1),\p(1)})\cdots
(b_{n,\p(n)}+c\d_{v(n),\p(n)}),$$
a polynomial of degree $n$, $p_0+p_1c+\cdots+p_{n-1}c^{n-1}+c^n$ over $A$ where
the $p_j$ do not depend on $c$. By Proposition 1.3.9 on p. 132 of $\ref{Dal}$,
$G(M_n(A))={\rm det}\inv(G(A))$, and we may now proceed as above.

\npar Theorem 1 can easily be applied to integral extensions (of which,
as noted in $\ref{LinIE}$, Arens-Hoffman
extensions are an example).

\cor 4.2.2.2. Let $A$ and $B$ be commutative, unital
normed algebras and suppose that $B$ is an integral
extension of $A$. Suppose that $A$ is a full subalgebra of its completion, $\tA$,
and that $A$ has dense invertible group. Then $B$ has dense invertible group.

\pf It is sufficient to prove the case when $B$ is an Arens-Hoffman extension
of $A$ for if $C$ is a normed, integral extension of $A$ and $c\in C$ is a root
of the monic polynomial $\a(x)\in A[x]$ then by Proposition 1.2.1.2 there is a continuous, unital
homomorphism $\theta\;\colon\,\Aa\to C$ with $\theta(\xb)=c$; the result quickly
follows from this.

Let $b\in \Aa$ and let $\e>0$.
It is not hard to show that the universal property of Proposition 1.2.1.2 allows us to identify
$(\tilde A)_\a$ with $\widetilde{\Aa}$. By Theorem 1
there exists $b\pri\in G((\tA)_\a)$ with $\norm{\bpr-b}<\e/2$.
We can write $\bpr=\sum_{k=0}^{n-1}\bpr_k\xb^k$ where $n$ is the degree of $\a(x)$ and
$\bpr_0,\ldots,\bpr_{n-1}\in \tA$. Since the Arens-Hoffman
extension of a Banach algebra is complete (see $\ref{Are}$),
$G((\tA)_\a)$ is open and we can perturb $\bpr_0,\ldots,\bpr_{n-1}$
to obtain $\bppr_0,\ldots,\bppr_{n-1}\in A$ so that $\bppr=\sum_{k=0}^{n-1}b^{\prime\prime}_k\xb^k$
is invertible in $(\tA)_\a$ and $\norm{\bppr-\bpr}<\e/2$. But now
$R_\a(\bppr)\in G(\tA)\cap A= G(A)$ so $\bppr$ is invertible in $\Aa$
and $\norm{b^{\prime\prime}-b}<\e$. \eop

\npar The following theorem follows directly from this corollary.

\thm 4.2.2.3. Let $B$ be a commutative Banach algebra which is an
integral extension of the commutative Banach algebra $A$. If
$A$ has dense invertible group then so has $B$.

\npar The converse of Theorem 1 seems harder to investigate. However the method of
the proof gives at least partial information in this direction.

\prop 4.2.2.4. Suppose that $A$ is a commutative, unital Banach algebra and $\a(x)$
is a monic polynomial of degree $n$ over $A$. If $\Aa$ has dense invertible group
then $\set{b\in A\st\hbox{there exists } a\in A\hbox{ such that } b=a^n}$, the set
of $n$th powers, is contained in the closure of $G(A)$.

\pf Fix $a\in A$ and let $\e>0$. Let the norm parameter for the Arens-Hoffman
extension be $t>0$. Then there exists $\b_\e(x)
=b_{\e,0}+\cdots+b_{\e,n-1}x^{n-1}\in A[x]$ such that $R_\a(\b_\e(x))\in G(A)$
and
$$\norm{\b_\e(x)-a}=\norm{b_{\e,0}-a}+\sum_{j=1}^{n-1}\norm{b_{\e,j}}t^j<\e.$$
Now as mentioned in Appendix 3, $R_\a(b_0+\cdots+
b_{n-1}x^{n-1})$ is homogeneous of degree $n$ in $b_0,
\ldots,b_{n-1}$. Therefore, writing $R_\a(\b(x))=P(b_0)=\sum_{j=0}^{n-1}p_jb_0^j
+b_0^n$ as in Theorem 1, we have that each coefficient $p_j\ (j=0,\ldots,
n-1)$ is a sum of elements of $A$ each of which has $b_k$ as a factor
for some $k\in\set{1,\ldots,n-1}$.

Thus, letting $\e\to0$, we obtain invertible elements $R_\a(\b_\e(x))$ in $A$ which
tend to $a^n$.\eop

\npar It is clear that the property of having a dense invertible group passes to (non-zero)
quotients and completions of normed algebras. The following is also clear.

\prop 4.2.2.5. Let $A$ be a direct limit of normed algebras with dense invertible groups
and in which the connecting homomorphisms are unital, isometric monomorphisms. Then
$A$ has dense invertible group.

\pf (Refer to Appendix 2 for notation.) Let
$$(A,(\ths i)_{i\in I})=\dl{}((A_i)_{i\in I},(\tss ij)_{i\le j}).$$
As we see from Appendix 2, the subset $D=\cup_{i\in I}\ths i(A_i)$ is dense
in $A$. It is therefore enough to show that for each $i\in I$, $\e>0$, and
$a\in A_i$, there exists an invertible element $b\in D$ such that $\norm{\ths i(a)-b}<\e$.

Since $G(A_i)$ is dense in $A_i$, there exists $c\in G(A_i)$ with
$\norm{a-c}<\e$. Since $\ths i$ is unital, $b:=\ths i(c)\in G(A)$
and since $\ths i$ is a contraction, $\norm{b-\ths i (a)}<\e.$\eop

\cor 4.2.2.6. Let $\u$ be an ordinal number and $(A_\s,(\tss\s\t)\st\s\le\t\le\u)$
be a system of algebraic extensions in $\Ba$ or $\ua$, in the notation of
Section 1.2.3. If $A_0$ has dense invertible group, then so has $A_\s$ for all
$\s\le\u$.

\pf This is a straightforward application of the transfinite induction theorem.

By assumption, $0\in\J:=\set{\t\le\u\st A_\t\hbox{ has dense invertible group}}$.
Let $0<\t\le\u$ and suppose that $[0,\t)\subseteq\J$.

First suppose there exists $\s<\t$ such that $\t=\s+1$. Then $A_\t=E(A_\s,\U_\s)$
for some set of monic polynomials, $\U_\s$, over $A_\s$. For each finite subset,
$S\in\U_\s\fs$, we have $G(A_S)$ dense in $A_S$ by Theorem 1 and induction.
(The extension $A_S$ is, as noted in $\ref{Nar}$,
easily checked to be isomorphic an algebra obtained by extending $A$ finitely-many
times by the Arens-Hoffman construction.) By Proposition 5, the Narmania extension,
$A_{\s,\U_\s}=\lim_{\rightharpoondown}(A_S\st S\in\U_\s\fs)$ has dense invertible
group.

If we are working in the category of uniform algebras, then the Cole extension
is $A_\t=A_\s^{\U_\s}=\overline{(A_{\s,\U_\s})\htt}$ by Proposition 1.3.2.1, and
has dense invertible group.

If $\t$ is a limit ordinal then $A_\t=\lim_{\rightharpoondown}(A_\s\st \s<\t)$
and has dense invertible group by hypothesis and Proposition 5.
By the transfinite induction theorem, $\J=\U$.\eop

\npar We mention the consequences of these results
for some of the examples of uniform algebras constructed by Lindberg, Karahanjan, Cole, and Feinstein
($\ref{LinIE}$, $\ref{Kar}$, $\ref{Col}$, $\ref{FeiNTSR}$).
Specifically we are referring to those obtained
by taking Cole extensions. In many of these the initial algebra has dense invertible group
and so the final algebra also has this property.

To illustrate this, consider Theorem 4 of $\ref{Kar}$ which we quote below.

\thm (Karahanjan). There is a non-trivial uniform algebra,
$A$, such that (1) $A$ is integrally closed, (2) $A$ is regular, (3) $\O$
is hereditarily unicoherent, (4) $G(A)$ is dense in $A$, (5) $A$ is
antisymmetric, and
(6) the Choquet boundary of $A$ is equal to $\O$.

\npar We refer the reader to $\ref{Kar}$ and the
literature on uniform algebras for the definitions of terms in the statement of this theorem
which we have
not defined here.

The proof of Karahanjan's theorem (p. 693-694 of $\ref{Kar}$) consists of checking that
the final algebra, $(A_\oi,\O_\oi)$, in a system of Cole extensions has the required properties
when $A_0$ is a suitable example, like McKissick's example ($\ref{McK}$) of a non-trivial, regular, uniform algebra
and the set of all monic polynomials  is used to pass from $A_\s$ to $A_{\s+1}$ for each
$\s<\oi$ (see Section 1.2.3).

Each of these specific choices of $A_0$ has a dense invertible group.
It follows immediately from Corollary 6
that $A_\oi$ also has dense invertible group.
It is therefore unnecessary to use Theorem 3 of $\ref{Kar}$ (on the complicated, auxiliary notion
of a `dense thin system') to deduce this in the proof of Theorem 4 of $\ref{Kar}$.


\def\Co{{\cal C}}
\def\tK{\tilde K}
\def\sf{SF(K)}
\def\SF#1{SF_{#1}(K)}
\def\ff{f_1,\ldots,f_n}
\def\gs{{\g^*}}

\vskip 15pt
\noindent {\medtenbf 4.3. Context: Uniform Algebras}

\vskip 15pt
\noindent {\bf 4.3.1. Introduction}\vskip 10pt

\noindent  The theorem below, which
first appeared in $\ref{DF}$,
prompted our investigation of the invertible group. In this section we extend the result.
Further context for this work will
be given in Section 4.3.3.

\thm 4.3.1.1. Let $A$ be a natural uniform algebra on $I$ with dense invertible group.
Then $A=C(I)$.
\pf By results mentioned in Section 3.2.1, every invertible element of $A$ has a continuous
logarithm in $A$ and therefore a square root in $A$. Therefore the set of elements of $A$
which have square roots in $A$ is dense in $A$. Since $I$ is locally connected, it
follows from \v Cirka's theorem (see
$\ref{Sto}$, p. 131-134) that $A=C(I)$.\eop

\vskip 15pt
\noindent {\bf 4.3.2. Generalisations of Theorem 4.3.1.1}\vskip 10pt

\noindent  In Theorem 1.7 of $\ref{DF}$, the hypotheses of Theorem 4.3.1.1 are relaxed so that
the norm on $A$ need not be the uniform norm and the maximal ideal space of $A$
can be replaced by any locally-connected space, $X$, such that $\Ho X=0$.
The conclusion is then that $\hA$ is dense in $C(\O)$.

We further show here that under certain conditions is suffices to restrict attention
to the components of the space. We require two elementary lemmas which are surely
well-known. Various parts of them are stated in in $\ref{Lei}$, but not collected
together as below.

Recall that if $(A,X)$ is a uniform algebra and $E\subseteq X$, then the
{\it $A$-convex hull} of $E$, denoted $\tilde E$, is the maximal ideal space of, $A_E$, the
uniform algebra on $E$ generated by $A\rest E$. We write $\chi_E$ for
the characteristic function of $E$. We shall call the set $E$ {\it clopen} if $E$
is both closed and open.

\lem 4.3.2.1. Let $(A,X)$ be a natural uniform algebra on a compact,
Hausdorff space $X$. Let the set of components of $X$ be denoted by $\Co$. Then
\item{(i)} for all $K\in\Co$, $K$ equals its $A$-convex hull, and
\item{(ii)} if $\Ho X=0$ then for all $K\in\Co$, $\Ho K=0$.

\pf (i) It is standard that, $\tK$, the $A$-convex hull of $K\in\Co$ is a subset
of $\O$, the maximal ideal space of $A$. Since
$K$ is connected, $A_K$
has no non-trivial idempotents. Therefore $\tK$ is connected by \v Silov's idempotent
theorem (see $\ref{Dal}$, p. 224). Now $K\subseteq\tK$ and $K$ is a component so it follows that $\tK=K$.

(ii) Let $\Ho X=0$, $K\in\Co$, and $f_0\in G(C(K))$. We show that $f_0$ has a continuous
logarithm. By Tietze's extension theorem (see,
for example, Proposition 1.5.8 in $\ref{Ped}$), there exists $f_1\in C(X)$ such
that $f_1\rest K=f_0$. By the compactness
of $K$ there exists an open set, $U$, containing $K$ and on which $f_1$ does
not vanish. Now since $X$ is compact and Hausdorff, $K$ is the intersection of
all those clopen subsets of $X$ which contain $K$. (See $\ref{Lei}$, p. 182, for
a proof of this fact.) A standard compactness argument
yields a clopen set, $E$, such that $K\subseteq E\subseteq U$. Since $E$ is clopen
we have $f_2=\chi_{X-E}+f_1\chi_E\in C(X)$. Now $f_2$ is clearly invertible so there
exists $g\in C(X)$ such that $f_2=e^g$. The restricted function $g\vert_K$ is
a logarithm of $f_0$.\eop

\lem 4.3.2.2. Let $A, X$, and $\Co$ be as in the lemma above. Then
\item{(i)} for all $K\in\Co$, $A_K=A\vert_K$, and
\item{(ii)} for all $f\in C(X)$, if for each $K\in\Co$, $f\vert_K\in A_K$ then $f\in A$.

\pf (i) As in the last lemma, each $K\in\Co$ is an intersection of
clopen sets. By the \v Silov idempotent theorem,
each clopen set is a peak set for $A$
and so (for example, by Theorem 3 on p. 163 of $\ref{Lei}$) we have $A_K=A\vert_K$.

(ii) Let $f$ satisfy the conditions stated above. By Bishop's theorem ($\ref{Lei}$, p. 39) it is enough
to show that $f\vert_E\in A_E=A\vert_E$ for every `maximal antisymmetric set', $E$, of
$A$. Let $E$ be a maximal antisymmetric set of $A$.
It is not necessary for us to define such sets; we merely need to
note from the remark on p. 170 of $\ref{Lei}$ that if
$E$ is a maximal antisymmetric set of $A$ then $\tilde E$ must be connected.
Therefore there exists $K\in\Co$ such that $E\subseteq\tilde E\subseteq K$.
Let $g\in A$ with $g\vert_K=f\vert_K$.
Then $g\vert_E=f\vert_E\in A_E$. The result follows.\eop

\noindent
The following result is the promised generalisation of Theorem 4.3.1.1.

\prop 4.3.2.3. Let $A$ be a commutative, unital Banach algebra with maximal
ideal space $X$ such that $\Ho X=0$. If $A$ has dense invertible group and every component
of $X$ is locally connected then $\hA$ is dense in $C(X)$.

\pf We assume that $(A,X)$ is a natural uniform algebra and show that $A=C(X)$. The
more general result then follows from applying this to the natural uniform algebra
$\overline{A\htt}$ (which has dense invertible group if $A$ has).

Let $f\in C(X)$. By Lemma 2, it is enough to show that $f\vert_E\in A_E$ for every
connected component, $E$, of $X$.

By Lemma 1, we have $\Ho E=0$ and the $A$-convex hull of $E$ is $E$ so $A_E=A\vert_E$
is a natural uniform algebra on $E$. By hypothesis $E$ is locally connected and $A_E$
has a dense exponential group. The
same argument as in the proof of Theorem 4.3.1.1 shows that $A_E=C(E)$.
The result now follows.
\eop

\npar There are apparently some classes of spaces in the literature (see $\ref{Cou}$)
for which components are automatically locally-connected.

We end this section by
using Proposition 3 to obtain a result about natural uniform algebras on the unit circle.
It is useful to recall some standard terminology and results beforehand. These
concepts will also be referred to in Chapter 5.

Let $A$ be a unital, normed algebra with space
of continuous characters $\O$. For $E\subseteq \O$ we define the {\it kernel} of $E$
to be $ k(E):=\set{a\in A\st\hbox{ for all }\o\in E\quad \o(a)=0}$. The {\it hull} of
$S\subseteq A$ is $h(S):=\set{\o\in \O\st\hbox{ for all }a\in S\quad \o(a)=0}$.
It is a standard fact that all such kernels are closed ideals of $A$.
Since $k(E)$ is a closed ideal of $A$ we can form the
normed algebra $C=A/k(E)$ with the quotient norm. Let the natural map be $\natural
\map AC$. As mentioned in Proposition 4.1.11 of $\ref{Dal}$, the following
map is a homeomorphism
$$\natural^*\mapto{\O(C)}{hk(E)}\o{\o\circ\natural}.$$
(This can be proved for normed algebras in the same way as for
Banach algebras.) Thus if $a\in A$
then under this identification, for all $\o\in hk(E)$, we have
$$(k(E)+a)\htt(\o)=\hat a(\o).$$

In order not to interrupt the flow of
the following proof, we describe a special case of this situation now.
Let $(B,X)$ be a natural uniform algebra and $E\subseteq X$. Then
the ideal $k(E)$ is the kernel of the restriction homomorphism $B\to B\vert_E$.
By the first isomorphism theorem we may therefore regard the
algebra of functions $B\vert_E$ as a Banach algebra under the quotient
norm on $B/k(E)$. Under this identification, the maximal ideal space of $B\vert_E$
is $hk(E)$, and for all $f\in B$ and $\o\in hk(E)$,
$$\widehat{f\vert_E}(\o)=(k(E)+f)\htt(\o)=\hat f(\o).$$
Finally, recall that a regular uniform algebra is
natural (see $\ref{Sto}$, Theorem 27.2 and Theorem 27.3).

\prop 4.3.2.4. Let $(A,S^1)$ be a regular uniform algebra. If $A$
has dense invertible group then $A=C(S^1)$.

\pf Suppose for a contradiction that $A$ is non-trivial and has dense invertible group.
We show that for each proper, closed set $E\subset S^1$ we have $A_E=C(E)$.
Uniform algebras with this property are called
`pervasive'. It follows from Theorems 20 and 21 on p.\ 173 of $\ref{Lei}$
that a non-trivial, pervasive uniform algebra is an integral domain. This will
contradict the regularity of $A$.

Let $E$ be a proper, closed subset of $S^1$. Then $E$ is contained in some
proper, closed interval, $I\subset S^1$. It is enough to show that $A_I=C(I)$.
To see this, let $f\in C(E)$ and $\e>0$.  Let $\tilde f$ be a Tietze extension
of $f$ to $I$. Assuming $A_I=C(I)$ then there exists $g\in A$
such that $\norm{g\rest I-\tilde f}_I<\e$. Then $\norm{g\rest E-f}<\e$.

Let $I$ be a proper, closed interval of $S^1$. As discussed above, we can regard
$A\rest I$ as a Banach algebra isomorphic to $A/k(I)$. Since $A$ has dense invertible
group then so has $A\rest I$ under this quotient norm. We also have
$\O(A\rest I)\approx hk(I)$. It follows from Proposition 4.1.20 of $\ref{Dal}$ and the regularity
of $A$ that $hk(I)=I$. As $\Ho I=0$ we have that $A\rest I$ is uniformly dense in
$C(I)$ by Proposition 3. This completes the proof.\eop


\vskip 15pt
\noindent {\bf 4.3.3. Uniform Approximation of Continuous Functions by Holomorphic Ones}
\vskip 10pt

\noindent Recall that a complex-valued function is called {\it space-filling}
if its image has non-empty interior. A consequence of Theorem 4.3.1.1
and Proposition 4.1.1.3 is that there
would be
an abundance of space-filling functions in a non-trivial, natural uniform algebra
on $I$ in the sense that there would be a non-empty, open subset of the algebra
consisting of such curves. For convenience we shall denote the set of
space-filling curves on a compact Hausdorff space, $K$, by $\sf$. If $A$
is a subset of $C(K)$ then we let $\SF A$ denote $\sf\cap A$.

It will also be useful to introduce the following notation for some
standard examples of uniform algebras. For a compact
subset $K\subset\C^n$ let $H_0(K)$ denote the algebra of functions holomorphic
on a neighbourhood of $K$, $R_0(K)$ the algebra of rational functions
which have poles off $K$, and $P_0(K)$ the algebra of polynomials
on $K$. The closures of $H_0(K)$, $R_0(K)$, and $P_0(K)$ in $C(K)$ are
denoted $H(K)$, $R(K)$, and $P(K)$ respectively. Also recall that a {\it simple
Jordan arc} in $\C^n$ is a continuous injection $\g\map I\C^n$.

The following result appears in $\ref{Cir}$ and should be compared with Proposition 4.3.2.3.

\thm 4.3.3.1 (\v Cirka). Let $\g$ be a simple Jordan arc in $\C^2$. Let $K=\im\g$.
If $H(K)\ne C(K)$ then for all $f\in H_0(K)-\C.1$, $f\in\sf$.\eop

\noindent It is easy to see that (unless $K$ is a singleton)
$H_0(K)-\C.1$ is dense in $H(K)$, so the theorem above asserts
the existence of a dense set of space-filling curves in $H(K)$, provided
$H(K)$ is non-trivial. At this point it is convenient to remind the reader about
the standard connections between polynomially-convex
arcs in $\C^n$ and finitely-generated subalgebras of $C(I)$. Recall that a set $K\subseteq\C^n$
is called {\it polynomially or rationally convex} provided
that $P(K)$ or $R(K)$ respectively is natural.

Suppose that $f_1,\ldots,f_n\in C(I)$ separate the points of $I$. Then
$\g\mapto{I}{\C^n}t{(f_1(t),\ldots,f_n(t))}$ gives a homeomorphism onto a compact
subset, $K$, of $\C^n$. As $K$ is the image of a simple Jordan arc and $\g$ induces an isometric isomorphism
$\gs\mapto{P(K)}Af{f\circ\g}$ where $A$ is the uniform subalgebra of $C(I)$ generated by
$\ff$. All these facts are discussed on p. 200 of $\ref{Lei}$.
Since $\ff$ generate $A$, the maximal ideal space of $A$ is naturally identified
(see, for example, $\ref{Sto}$, p. 25)
with $K$ by the homeomorphism $\g$ if $A$ is natural. In particular $K$ is polynomially convex.
It is clear that $P(K)=C(K)$ if and only if $\gs(A)=C(K)$ or, equivalently, $A={\gs}^{-1}(C(K))=C(I)$.
If we  take $A$ to be a natural, finitely-generated uniform algebra on $I$, then
Theorem 4.3.1.1 yields the following corollary.

\cor 4.3.3.2. Let $\g$ be a simple Jordan arc in $\C^n$ such that
$K=\im\g$ is polynomially convex. If $H(K)\ne C(K)$ then there
is a non-empty open set $U\subseteq H(K)$ such that $U\subseteq \sf$.\eop

\noindent
This makes the differences clearer. In Theorem 1, the arc need not be
polynomially convex, but is restricted to $\C^2$, and $H(K)\ne C(K)$ implies
$\SF{H(K)}$ is dense. However
in Corollary 2 the conclusion is that $\SF{H(K)}$ has non-empty interior.

\npar We end this section by noting that it is possible to improve Corollary
2 by using a form of the Oka-Weil approximation theorem which is more
general than the version usually quoted. The following version is taken from
the appendix in $\ref{StoCon}$.

\thm 4.3.3.3 (Oka-Weil Approximation Theorem). Let $K$ be a compact subset of $\C^n$.
Then $H(K)=P(K)$ if $K$ is polynomially convex and $H(K)=R(K)$ if $K$ is
rationally convex.\eop

\npar There is also (see $\ref{Sto}$, p. 372) a notion of the `holomorphic
convex hull' of a compact subset $K\subseteq\C^n$: it is equal by definition to the
maximal ideal space of $H(K)$. We call $K$ {\it holomorphically convex} if
$H(K)$ is natural. It is remarked in $\ref{Sto}$ that this hull
need not be a subset of $\C^n$ if the original space is.

\cor 4.3.3.4. Let $\g$ be a simple Jordan arc in $\C^n$ such that $K=\im\g$ is
holomorphically convex. If $H(K)\ne C(K)$ then
$SF_{H(K)}(K)$ has non-empty interior (in $H(K)$).

\pf The homeomorphism $\g\map IK$ induces a commutative diagram in
which the horizontal arrows are isometries:
$$\matrix{ H(K)&\mapright\gs & A:=\gs(H(K)) \cr
	\mapdown\subseteq & &\mapdown\subseteq \cr
	C(K) &\mapright\gs & C(I) \cr}$$
The assumption on $K$ implies that $H(K)$ is natural and so $A$ is natural.
If not every continuous function of $K$ can be approximated by holomorphic
functions on a neighbourhood of $K$ then $A$ is non-trivial and so $U:=A-\bGA$ is a non-empty,
open subset of $SF_A(I)$. Hence $\gs^{-1}(U)$ is a non-empty, open subset
of $\SF{H(K)}$.\eop

\npar Note that, by Theorem 3, polynomially convex sets and rationally
convex sets are holomorphically convex so the
conclusion of Corollary 4 also holds if $K$ is polynomially
or rationally convex.

\vfill\eject


\noindent {\bigtenbf Chapter 5}
\vskip 10pt
\noindent {\bigtenrm Uniform Algebras with Connected, Dense Invertible Group}
\vskip 25pt

\noindent
The main result, Proposition 4.3.2.3, of Section 4.3 makes strong use of \v Cirka's theorem
(see $\ref{Sto}$, p. 131-134).

\thm (\v Cirka). Let $A$ be a uniform algebra on locally-connected space, $X$, such that
the set,
$S$, of squares of elements of $A$ is dense in $A$.
Then $A=C(X)$.
\npar
However in the application in Chapter 4, a condition stronger than $\overline{S}=A$
held, namely that $e^A$ was dense in $A$.
The only other condition in the hypothesis of \v Cirka's theorem is that $X$ be locally connected.
Could it be that this is redundant if we have $\overline{e^A}=A$? That is, does the condition
$\overline{e^A}=A$ imply that $A=C(X)$ for any uniform algebra on a compact, Hausdorff space, $X$?

The statement is false if the conditions on $A$
are weakened so that its norm is not necessarily the uniform norm. For example,
$A=C^1(I)$ (see Section 4.1.2) under the norm $\norm f:=\norm f_I+\norm{f\pri}_I$ is not
equal to $C(I)$. To see that $\overline{e^A}=A$ in this case, we have from Theorem 4.4.1
of $\ref{Dal}$ that the polynomials in $\id I$ are dense in $A$ and so by Proposition
4.2.1.5, $A$ has dense invertible group. We also have from Theorem 4.4.1 that the
maximal ideal space of $A$ is homeomorphic to $I$ so by the Arens-Royden theorem, $G(A)=e^A$.

Note that for
a general commutative Banach algebra $A$, $\overline{e^A}=A$ is equivalent to $\overline{G(A)}=A$ and $\Ho\O=0$. It is clear that $G(A)$ is dense in $A$ if $\overline{e^A}=A$. It follows from
Lemma 3.2.3.1 that every invertible element of $A$ belongs to $e^A$ in this case. Therefore
we also have $\Ho\O=0$. Conversely, if $\Ho\O=0$ then $G(A)=e^A$ so if $\overline{G(A)}=A$
then $\overline{e^A}=A$. Compare this pair of conditions
with those in Lavrentieff's theorem (see, for example, $\ref{Ale}$, p. 12).
This theorem states that if $X$ is a compact, plane-set
and $A=P(X)$ then $A=C(X)$ if and only if $X\int=\emptyset$ and $\Ho X=0$.

In this chapter we search for a counterexample to the following conjecture.

\npar CONJECTURE. Let $(A,X)$ be a uniform algebra such that  $\overline{e^A}=A$. Then $A=C(X)$.

\npar We have not been able to find such an example. Various classes of uniform algebras are
ruled out in Sections 1 and 2. The bulk of this chapter is devoted to modifying the methods of
Chapter 1 to construct `logarithmic extensions'. We show how this
might lead to an example of a non-trivial uniform algebra with dense
exponential group. Various simple properties of
these extensions will be established. For example, in certain situations they are non-trivial if the
original algebra is.


\def\inte{^\circ}
\def\D{\Delta}
\def\Dp{\Delta\pri}

\vskip 15pt
\noindent {\medtenbf 5.1. Uniform Algebras on Compact Plane-Sets}

\vskip 15pt
\noindent {\bf 5.1.1. Introduction}\vskip 10pt

\noindent In this section $X$ denotes a compact subset of $\C$. We show that the conjecture
holds when $A=P(X)$ or $R(X)$. The latter algebra has been an important source of counterexamples
to conjectures about uniform algebras in the past. (See $\ref{McK}$ and $\ref{Col}$.)

Recall that a multiplicative group, $G$, is said to be {\it divisible by $n\in\N$} provided
that for all $g\in G$ there exists $h\in G$ such that $g=h^n$.

\vskip 15pt
\noindent {\bf 5.1.2. Generalisations of Lavrentieff's Theorem}\vskip 10pt

\noindent We begin with a lemma which is valid in more generality than is necessary
for the main result in this section.

\lem 5.1.2.1. Let $F$ be an entire function
with no zero of order one.
Let $A=P(X)$ or $R(X)$ and suppose that $F(A)$ is dense in $A$. Then $X$
has empty interior.

\pf Suppose for a contradiction that $\l\in X\inte$, the interior of $X$. Then
$f=z-\l\in A$ where $z:=\id X$. Let $r>0$ be such that $\D=B[\l,r]\subseteq X\inte$. Then
$\norm{f}_\D=r>0$. By assumption there exists $g_2\in A$ such that
$\norm{F\circ g_2-f}_X<({1\over3})\norm{f}_\D$.

An exercise in uniform continuity shows that $F_*\mapto{C(\D)}
{C(\D)}g{F\circ g}$ is continuous so there exists a rational function
$g_1\in A$ such that $\norm{F\circ g_1-f}_\D<({2\over3})\norm{f}_\D$.

By Rouch\' e's theorem, $F\circ g_1$ has exactly one zero in $\D\inte$,
and its multiplicity is one.

Let $\m\in\C$ and $s>0$ be such that $\Dp= B[\m,s]\subseteq\D\inte$ and $(F\circ g_1)(\m)=0$.
By Taylor's theorem, since $F\circ g_1$ is analytic on an open disc, $U$,
centre $\l$ and containing $\D$, the power-series expansion
$$F\circ g_1=a_1(z-\m)+a_2(z-\m)^2+\cdots\qquad\qquad(a_1\ne0)\eqno(*)$$
is valid on an open disc $V\supseteq\Dp$, centre $\m$. Since $F$ is entire
with no zero of order 1, there is an entire function $h$ such that
$F=(z-g_1(\m))^2h$. But now we have, for some function, $k$, which is
analytic at $\m$, that, for all $w$ near $\m$,
$$\eqalign{F(g_1(w))&=(g_1(w)-g_1(\m))^2h(g_1(w))\cr
&=(w-\m)^2k(w),\cr}$$
contradicting ($*$).\eop

\noindent We can now give the main result of this section.

\cor 5.1.2.2.
Let $A=P(X)$ or $R(X)$. Suppose that
there exists $n\in\N$ ($n>1$)
such that the set of $n$th powers of elements of $A$ is dense in $A$.
Then $A=C(X)$.

\pf We show that $P(X)$ is equal to $C(X)$. By Lavrentieff's theorem
this is equivalent to showing that
$X$ has no interior and $\C-X$ is connected. The first condition
follows from Lemma 1 above. We prove that $\C-X$ has only one component.
It is elementary that $\C-X$ has countably-many components, of which
exactly one is unbounded.

Let
$M$ be a transversal for the bounded components of $\C-X$.
Recall (or see $\ref{BroTop}$) that $\Ho X\cong\oplus_{m\in M}\Z$ (the
group, under pointwise addition, of elements
$h\in \Z^M$ such that $h$ has a finite number of non-zero entries).
We must
show that $\Ho X=0$. Since $G(A)$ is open, it follows from the denseness
of the $n$th powers in $A$ that $G(A)$ is a subset of the $n$th powers.
For let $g\in G(A)$. By Theorem 10.43 of $\ref{RudFA}$,
the coset $e^Ag$ is open so there exists $h\in A$ such that $h^n\in e^Ag$.
Thus $g\in e^Ah^n$ and so $g$ has an $n$th root in $A$. Therefore
$\Ho X$ is divisible by $n$ and must be trivial.\eop

\noindent The next result is a corollary of the last proof. The situation
described contrasts strongly with that in Section 5.2.

\prop 5.1.2.3. Let $A$ be a uniform algebra whose maximal ideal
space, $\O$, is homoeomorphic to a compact subset of $\C$ (for example,
if $A$ has a single generator). Then
$\Ho\O=0$ if and only if for every $a\in G(A)$ there exists
$b\in A$ such that $b^2=a$.
\pf Suppose that every invertible element of $A$ is square of another.
It follows as in the proof of Corollary 2 that $\Ho\O$ is divisible
by $2$ and therefore that $\Ho\O=0$. The converse is trivial.\eop

 \def\D{{\bf D}}       

\def\inte{^\circ}
\def\D{\Delta}
\def\Dp{\Delta\pri}

\vskip 15pt
\noindent {\medtenbf 5.2. Systems of Cole Extensions}
\vskip 15pt

\noindent A necessary condition for $\overline{e^A}=A$ is that $\Ho\O=0$, and in the previous section
we used (in a special case) the divisibility of the group $\Ho\O$ to show that it was the trivial
group. Even if a non-trivial uniform algebra, $(A,\O)$, with dense invertible group,
does not satisfy $\overline{e^A}=A$
one might imagine that we could apply the Cole construction to it to obtain a such an algebra.
The reasoning is that an integrally closed extension,
$A_\oi$, in a system of Cole extensions shares the properties of being non-trivial and having dense
invertible group. (This follows from Theorem 1.4.2.5 and Corollary 4.2.2.6.) Additionally, $\Ho{\Om{A_\oi}}$
is now unrestrictedly divisible, and an easy way for this to occur is if it is trivial.

Further grounds for suspicion are that for every $n\in\N$ and $a\in G(A_\oi)$,
the polynomial $P_n(x):=(1-a)+\sum_{k=1}^n
x^k/k!$ has a solution in $A_\oi$. (We have $P_n(x)\to e^x-a$ in $A_\oi[[x]]$ under the
norms discussed in Section 4.2.)

However, if $\Ho{\Om{A_\oi}}=0$ then Corollary 3.2.3.5 shows that $\Ho\O=0$, and
so the original algebra would have to have dense exponential group, which we assumed was untrue.
Therefore Cole extensions can not generate a counterexample out of an algebra with
dense invertible group which is
not already a counterexample.

In the next section, we vary the Cole construction by adjoining logarithms of invertible
elements rather than
roots of monic polynomials. This eventually produces extensions in
which every invertible element has a logarithm.

 \def\D{{\bf D}}       


\vskip 15pt
\noindent {\medtenbf 5.3. Logarithmic Extensions}

\vskip 15pt
\noindent {\bf 5.3.1. Introduction}\vskip 10pt

\def\na{{n(\a)}}
\def\ff{f_1,\ldots,f_n}
\def\gs{{\g^*}}
\def\Ala{A_{\log a}}
\def\laA{A^{\log a}}
\def\Ola{\O_{\log a}}
\def\bp{\bar p}
\def\hl#1{\widehat{l_{#1}}}
\def\ci#1{\overline{G(#1)}}
\def\tb{\tilde b}

\def\hT{\hat T}
\def\whu#1{\widehat{u_{#1}}}
\def\r{\rho}
\def\oo{{\o_0}}
\def\fs{^{<\oo}}
\def\yb{\bar y}
\def\V{{\cal V}}
\def\AlU{A_{\log \U}}
\def\lUA{A^{\log \U}}
\def\OlU{\O_{\log \U}}
\def\pp#1#2{\pi_{#1,#2}}
\def\tpp#1#2{\tilde\pi_{#1,#2}}
\def\pps#1#2{\pp#1#2^*}
\def\lSA{A^{\log S}}
\def\OlS{\O_{\log S}}
\def\dl{{\lim_\rightharpoondown}\mathstrut}	
\def\il{{\lim_\leftharpoondown}\mathstrut}		
\def\mapne#1{\nearrow${$\scriptstyle#1$}$}
\def\mapsw#1{ \swarrow${$\scriptstyle{#1}$}$  }
\def\uis{{\U\in S\fs}}

\def\het{\hat \theta}


\noindent Consider the following problem. Given a natural uniform
algebra $(A,\O)$ and $a\in G(A)$ find an extension
$(\laA,\Ola)$ such that:
\item{(i)} $\laA\ne C(\Ola)$ if $A\ne C(\O)$,
\item{(ii)} $\laA$ has dense invertible group if $A$ has,
\item{(iii)} there exists $b\in \laA$ such that $e^b=a$, and
\item{(iv)} the process can be repeated (similarly to in Section 1.2.3)
so that (i) and (ii) are preserved by the final algebra and
the final uniform algebra has dense exponential group.

If we can solve this, then we can apply the construction to a non-trivial uniform algebra
with dense invertible group (such algebras exist: McKissick's $R(X)$ in $\ref{McK}$ is one)
and obtain a non-trivial uniform algebra with dense exponential group. In this chapter
we give a very promising construction. It is known to satisfy all the above requirements
except that we do not yet know if the final algebra must be non-trivial, even for specific
examples.


\vskip 15pt
\noindent {\bf 5.3.2. Adjoining a Logarithm to a Banach Algebra}\vskip 10pt


\noindent The problem of adjoining logarithms to semisimple, commutative Banach algebras has
been studied by Lindberg. We summarise some of this work briefly.

The terms `local' and `semisimple' are defined in Section 2.2 and Appendix 1 respectively. We
further require the following standard terminology. Let $A$ be a unital, normed algebra with space
of continuous characters, $\O$.
Let $a\in G(A)$ and $\o\in\O$. We say that {\it $e^x=a$ has local solutions at $\o$}
if there exists $V\in\Nb\o$ and $b\in A$ such
that $\exp\left(\hat b \right)\rest V=\ha\rest V$. We shall also make use of the
notions of the hull and kernel; these were defined in Section 4.3.2.
Finally, for use later in the chapter, we
recall that if $(V_k)_{k=1}^n$ is an open cover
of a topological space, $X$, then a {\it partition of unity} subordinate to this cover
is a sequence of continuous functions, $(f_k\map X\C)_{k=1}^n$, such that
$$\sum_{k=1}^nf_k=1\qquad\hbox{and}\qquad{\rm supp}\,(f_k):=\overline{\set{\k\in X\st
f_k(\k)\ne0}}
\subseteq V_k\quad(k=1,\ldots,n).$$

\thm 5.3.2.1 ($\ref{LinEANA}$). Let $A$ be a local, semisimple, commutative, unital Banach algebra
and $a\in G(A)$. Suppose that $e^x=a$ has local solutions at each point of $\O$.
Then
$e^x=a$ has a solution in a complete, normed, semisimple extension, $\Ala$, of $A$. Furthermore,
\item{(i)} the extension is $\Ala=A(x;t)\bigm/ kh((e^x-a))$ where $A(x;t)$
is the completion of the algebra of polynomials over $A$ in the norm
$$\norm{\sum_{k=0}^n a_kx^k}=\sum_{k=0}^n \norm{a_k}t^k,$$
for a certain fixed value of $t>0$ (see Section 4.2.1);
\item{(ii)} the norm on $\Ala$ is equivalent to and dominates the quotient norm on the factor algebra;
\item{(iii)} the isometric embedding is given by $\th:A\to\Ala\;;\;b\mapsto J+b$
where $J=kh((e^x-a))$;
\item{(iv)} the maximal ideal space of $\Ala$ is $$\Ola:=\Om{\Ala}=\set{
(\k,\l)\in\O\times B_\C[0,t]\st e^\l=\ha(\k)},$$ where for $\ol\in\Ola$,
$$\ol\left(J+\sum_{k=0}^{+\infty}b_kx^k\right)= \sum_{k=0}^{+\infty}\o(b_k)\l^k
\qquad\qquad\left(\sum_{k=0}^{+\infty}b_kx^k\in A(x;t)\right);$$
\item{(v)} the map $\pi\mapto{\Ola}\O\ol\o$ is a continuous surjection.

\npar The value of $t$ defining the extension is chosen to be at least the
maximum of the norms of some choice of $b_1,\ldots,b_m\in A$ such that
$\widehat{b_k}$ is a local solution
of $e^x=a$ valid on an open set $V_k\subseteq\O$ and $\cup_{k=1}^m V_k=\O$. In
future, if we speak of `a norm parameter $t$' for $\Ala$ then it will be
understood that $t$ is a bound for some set of local solutions to $e^x=a$
which are valid on open sets covering $\O$.
Clearly $J+x$ is a solution to $e^x=a$ in $\Ala$.

We can use the description of $\Ola$ to extend Lindberg's
result to uniform algebras in the following minor way.

\prop 5.3.2.2. Let $A$ be a local uniform algebra and let $a\in A$ be such that $e^x=a$ has
local solutions in $A$. Set $\laA=\overline{\Ala\htt}$. Then $\laA$ is a uniform-algebra
extension of $A$ in which $e^x=a$ has a solution, namely $p:=(J+x)\htt$.

The isometric embedding $\het\map A{\laA}$ is given by the restriction to $\laA$
of the adjoint map $\ps\mapto{C(\O)}{C(\Ola)}f{f\circ\pi}$.

\pf Of course $(\laA,\Ola)$ is a natural uniform algebra by Lemma A.1.4.
Define $\het\mapto A\laA b{(J+b)\htt}$.
Thus $\het=\ps$, if we identify $\hA$ with $A$. We show that this map is isometric.
For $b\in A$
we have
$$\eqalign{\norm{\het(b)}_{\laA}&=\sup_{(\o,\l)\in\Ola}\ab{(J+b)\htt(\o,\l)}
=\sup_{(\o,\l)\in\Ola}\ab{\hat b(\o)}\cr
&=\sup_{\o\in\O}\ab{\hat b(\o)}=\norm b.\cr}$$
The last line is true because the canonical map $\pi\mapto\Ola\O{(\o,\l)}\o$ is
surjective.\eop

\npar From now on we shall regard $\pi$, $\het$, and $p$ as
objects associated with the extension $\laA$.
In contrast to the Arens-Hoffman construction, the map $\pi$ need not be open
as the following example shows.

\ex 5.3.2.3. Let $A=C(S^1)$. (The unit circle, $S^1$, is one of the simplest spaces
on which not every invertible,
continuous function has a continuous logarithm; see, for
example, p. 28 of $\ref{Rao}$.) Let $a=\id{S^1}$. Then $A$
satisfies the conditions for the construction. We can clearly take
$$\eqalign{
V_1&=\set{e^{\ri\th}\st -3\pi/4<\th<3\pi/4},\cr
V_2&=\set{e^{\ri\th}\st \th\in(-\pi,-\pi/4)\cup(\pi/4,\pi]},\cr}$$
and continuous functions $f_1,f_2\in C(S^1)$ such that
$$\eqalign{ &f_1(w)=\arg w\qquad(w\in V_1),\cr
&f_2(w)=\pi-\arg(\ri w)\qquad(w\in V_2),\cr}$$
and $\norm{f_1},\norm{f_2}\le 2\pi$. We may take $t=2\pi$.

If $\pi$ were an open map then $\pi(W)$ would be open where
$$W=\pi^{-1}(\set{e^{\ri\th}\st -\pi/2<\th<\pi/2})\cap p^{-1}(B_\C(2\pi\ri,\pi/4)).$$
Note that $(1,2\pi\ri)\in\Ola$. However simple computations show that
$$\pi(W)=\set{e^{\ri\th}\st -\pi/4<\th\le0}.$$ So $\pi$ is not an open map.\eop


\vskip 15pt
\noindent {\bf 5.3.3. Properties of Simple Logarithmic Extensions}\vskip 10pt


\noindent We refer to algebras constructed by adjoining a logarithm of a
single invertible element as {\it simple}.
In this section we show that the logarithmic extension $\laA$
must share the properties of non-triviality (under certain conditions),
denseness of the invertible group, and
regularity with $A$. We focus on regularity because Lindberg's construction requires
a local base-algebra, localness is implied by regularity, and, as we shall see,
regularity is shared by $\Ala$. Furthermore, as remarked in
$\ref{LinEANA}$, for every invertible element, $a$, of a regular Banach algebra $A$
there are local solutions to $e^x=a$ in $A$. Thus, in this case,
there is no need to restrict the elements of $G(A)$ we
attempt to provide logarithms for.

\npar
First we give a criterion for the construction to satisfy requirement (i)
of Section 5.3.1. Let $(A,\O)$ be a natural uniform algebra, $f\in C(\O)$,
and $\o\in\O$. Recall that $f$ is {\it locally approximable by $A$ at $\o$}
if there exists an open neighbourhood, $V$,
of $\o$ such that $f\vert_{\bar V}\in A_{\bar V}$ (see Section 4.3 for
this notation). The algebra $A$
is called {\it holomorphically closed} if $f\in A$
whenever $f$ is locally approximable by $A$ at every point of $\O$. For example, if $K$ is a
compact plane-set then $R(K)$ is holomorphically closed.

To see this standard result, let $f\in C(K)$ and $V_1,\ldots,V_n$ be an open cover of $K$ such
that $f\vert_{\bar V_k}\in A_{\bar V_k}\ (k=1,\ldots,n)$. Thus for $k\in
\set{1,\ldots,n}$ and $\e>0$ there exists a rational function, $q$, with poles
off $K$ such that $\norm{f-q}_{\bar V_k}<\e$. Hence $f\vert_{\bar V_k}\in
R({\bar V_k})$. By Theorem 3.2.13 on p. 170 of $\ref{Bro}$ this implies that $f\in
R(K)$.

Clearly localness is implied by being holomorphically closed. For if $A$
is holomorphically closed and $f\in C(\O)$ belongs locally to $A$, let
$(V_k)_{k=1}^n$ be an open cover of $\O$ and $f_1,\ldots,f_n\in\hA$ be such
that $f_k\rest{V_k}=f\rest{V_k}$ for $k=1,\ldots,n$. Then $f\rest{\overline{V_k}}
\in A\vert_{\overline{V_k}}$  ($k=1,\ldots,n$) so $f\in\hA$.

It is also standard that a regular uniform algebra is holomorphically closed; this
follows from the existence (see $\ref{Dal}$, p. 413) of partitions of unity
in the algebra subordinate to any open
cover of the maximal ideal space.

\prop 5.3.3.1. Let $(A,\O)$ be a holomorphically closed, natural uniform algebra. Let
$a\in G(A)$ be such that $e^x=a$ has local solutions in $A$ at each point of $\O$.
Then $\laA\ne C(\Ola)$ if $A\ne C(\O)$.

\pf Suppose $A\ne C(\O)$. There must exist
$f\in C(\O)-A$. We show that $\ps(f)$ does
not belong to $\laA$. Suppose for a contradiction that $\ps(f)\in\laA$.

By the choice of norm parameter there exist $g_1,\ldots,g_n\in B_A[0,t]$ and open sets
$V_1,\ldots, V_n$ in $\O$ which cover $\O$ such that
$$e^{g_k(\o)}=\ha(\o)\qquad(\o\in\bar V_k,\quad k=1,\ldots,n).$$
Since $\Ala\htt$ is uniformly dense in $\laA$, for every $\e>0$ there exist
$f_1,\ldots,f_m\in A$ such that
$$\norm{ \sum_{j=1}^m\ps(f_j)p^j-\ps(f)}<\e$$
where $p$ is as above. Hence for
$k\in\set{1,\ldots,n}$ and all $\o\in \bar V_k$ we have,
setting $F_\e=\sum_{j=1}^m f_jg_k^j\in A$, that
$$\ab{F_\e(\o)-f(\o)}= \ab{ \sum_{j=1}^m\left(\ps(f_j)p^j\right)(\o,g_k(\o))-\ps(f)(\o,g_k(\o))}
<\e.$$
Since $\e$ is arbitrary it follows that $f\vert_{\bar V_k}\in A_{\bar V_k}$. As $A$ was assumed to be
holomorphically closed this implies that $f\in A$, a contradiction.\eop

\npar The converse result can be proved using the Arens-Calder\' on theorem.

\prop 5.3.3.2. Suppose that $A=C(\O)$ and $a\in G(A)$. Then $\laA=C(\Ola)$.

\pf We show that $\Ala\htt$ is self-adjoint; the result then follows from
the Stone-Weierstrass theorem since this subalgebra clearly separates the
points of $\Ola$ and contains the constant functions. It is clearly enough to show that
$\bar p$, the pointwise complex-conjugate of $p$ belongs to $\left(\Ala\right)\htt$.

We have that for all $\ol\in\Ola$,
$$e^{\bp\ol}=(e^{p\ol})^-=\ha(\o)^-=(\bar a\htt)(\o)=(\bar a\htt\circ\pi)\ol=\het(\bar a)\ol$$
and $\th(\bar a)\in \Ala$ so by the Arens-Calder\' on theorem (see, for example,
$\ref{Dal}$, p. 223, Corollary 2.4.31(i)) there exists (a unique) $b\in\Ala$
such that $e^b=\th(\bar a)$ and $\hat b=\bp$.\eop

\def\hB{\hat B}
\def\hb{\hat b}


\npar We now show that regularity is preserved by simple logarithmic extensions
of uniform algebras. In fact we shall prove a more general result. The
proof of the next proposition is a modification of an argument of Feinstein from
$\ref{FeiNTSR}$. Note that the map $\pi$ is not required to be open or surjective.

\prop 5.3.3.3. Let $\pi\map YX$ be a continuous map where
$X$ and $Y$ are the maximal ideal spaces of commutative, unital Banach algebras
$A$ and $B$ respectively. Assume that
\item{(a)} $\pi$ induces a map $\ps\mapto{\hA}\hB\ha{\ha\circ\pi}$,
\item{(b)} for each $x\in X$, $\pi\inv(x)$ is totally disconnected,
\item{(c)} $A$ is regular.

Then $\overline{B\htt}$ is a regular uniform algebra.

\pf Let $y_1$ and $y_2$ be distinct points of $Y$. We show that there exists
$h$ in the closure of $\hB$ such that $h(y_1)=1$ and $y_2$ belongs to the
interior of $h\inv(0)$. Set $x_j=\pi(y_j)\ (j=1,2)$.

First suppose that $x_1\ne x_2$. Then by the regularity of $A$ there exists $a\in A$
with $\ha(x_1)=1$ and $x_2\in\ha\inv(0)\int$. Take $b\in B$ with $\hat b=
\ha\circ\pi$. Then $\hb(y_1)=\ha(x_1)=1$ and $y_2\in\pi\inv(\ha\inv(0)\int)
\subseteq \pi\inv(\ha\inv(0))=\hb\inv(0)$.

Now suppose that $x_1=x_2=:x$. Consider the commutative Banach algebra
$C=B\bigm/k(\pi\inv(x))$. By results mentioned in Section 4.3.2,
the maximal ideal space of
$C$ is homeomorphic to $hk(\pi\inv(x))$. It is easy to see that $hk(\pi\inv(x))
=\pi\inv(x)$: we have $\pi\inv(x)\subseteq hk(\pi\inv(x))$ ($\ref{Dal}$, p. 411, again).
Now let $y\in hk(\pi\inv(x))$ and suppose that $\pi(y)=:w\ne x$. By the
case proved in the last paragraph there exists $b\in B$ such that
$\hat b(\pi\inv(x))=\set0$ and $\hat b(y)=1$. But then
$b\in k(\pi\inv(x))$ and this contradicts the
fact that $y\in hk(\pi\inv(x))$. Hence $hk(\pi\inv(x))= \pi\inv(x)$
is totally disconnected by
assumption (b). It follows that we can find a clopen set $G\subseteq\Om C$
with $y_1\in G$ and $y_2\not\in G$. By the \v Silov idempotent theorem
($\ref{Dal}$, p. 224) there exists some $c\in C$ such that $\hat c=\chi_G$.
Hence, under the identification of $\O(C)$ with $hk(\pi\inv(x))$, there exists
$b\in B$ with $\hb(G)=\set1$ and $\hb(\Om C-G)=\set0$.

Set $B_1=B_\C(1,1/3)$, $B_2=B_\C(0,1/3)$, and
$$\eqalign{&U = \hat b\inv(B_1\cup B_2);\cr
&K = Y-U;\cr
&C_j = \hat b\inv(\overline{B_j})\qquad (j=1,2).\cr}$$
Then $U$ is an open neighbourhood of $\pi\inv(x_1)$ and $C_1$ and $C_2$ are compact.
We also have that $K\subseteq Y-\pi\inv(x_1)$ and so by what was
proved in the second paragraph and the compactness
of $K$ there exists $g\in \hB$ and an open superset, $V$, of $K$ such that
$g(V)\subseteq\set0$ and $g(y_1)=1$.

By Runge's theorem ($\ref{RudRCA}$, p. 272) there exists a sequence of complex polynomials, $(q_k)$, such
that
$$q_k\to\cases{0 & on $\overline{B_1}$\cr
1 & on $\overline{B_2}$\cr}\qquad\hbox{ uniformly as }k\to+\infty.$$
Set $f_k=(q_k\circ \hat b)g\ (k\in\N)$. Then $(f_k)$ is a Cauchy sequence in $\hB$
and has some limit $f\in\overline{B\htt}$. We see that
$$f(y)=\cases{0 &if $y\in Y-C_1$, \quad and\cr
g(y) &if $y\in C_1$.\cr}$$
Now $y_2\in C_2\subseteq Y-C_1\subseteq f\inv(0)$ and $f(y_1)=g(y_1)=1$
so the result is proved.\eop

\npar Now we show that the denseness of the invertible group is preserved by the construction.
We retain the
notation used after Theorem 5.3.2.1.

\lem 5.3.3.4. Let $A$ be a local, semisimple, commutative, unital Banach algebra.
Let $a\in G(A)$ have local solutions at each point of $\O$. Then there exist
finitely many $l_1,\ldots,l_m\in A$ such that the maximal ideal space of the
Lindberg extension, $\Ala$, is a subset of the union of the graphs of
$\hl1,\ldots,\hl m.$ That is, for each $\ol\in\Ola$ there exists $j\in\set{1,
\ldots,m}$ such that $\l=\hl j(\o)$.

\pf By the assumption on $A$, $a$, and the norm-parameter, $t$, there exist
$l_1,\ldots,l_m\in B_A[0,t]$
and open sets $V_1,\ldots,V_m\subseteq \O$ which cover $\O$ such that
$\hl j$ is a local solution of $e^x=\ha$ on $V_j\ (j=1,\ldots,m)$. Let
the sequence $(l_1,\ldots,l_m)$ be extended by $l_{m+1},\ldots,l_M\in A$
where the additional items are obtained by adding $2\pi\ri k$ to all
the $l_1,\ldots,l_m$ for all $k\in\Z$ such that $\ab{k}\le t/\pi$.

For $\ol\in\Ola$ we have $e^\l=\ha(\o)$ and $\ab\l\le t$. Let $j\in\set
{1,\ldots,m}$ with $e^{\hl j(\o)}=\ha(\o)$. Then there exists $k\in\Z$ with
$\hl j(\o)=\l-2\pi\ri k$. Since $t$ was chosen to satisfy
$t\ge\norm{l_1},\ldots,\norm{l_m}$,
$$\ab k ={ \ab{\hl j(\o)-\l}\over{2\pi} }\le {t\over\pi}$$
and so there is some $j\pri\in\set{1,\ldots, M}$ with
$\l=\hl{j\pri}(\o)$.\eop

\prop 5.3.3.5. Let $A$ and $a$ be as above. Suppose that $\ci A=A$. Then
$\Ala$ has dense invertible group.

\pf Suppose $A$ is as in the statement. Since the
polynomials over $A$ are dense in $A(x;t)$ it is sufficient to find,
for each $b_0,\ldots,b_n\in A$ and $\e>0$, an invertible element $g\in\Ala$
with $\norm{g-b(\xb)}<\e$ where $b(x)=\sum_{k=0}^nb_kx^k$ and
$\xb:=J+x$. In fact we shall show that the element
$g$ can be obtained by perturbing $b_0$ only.

By Lemma 4 there exist $l_1,\ldots,l_m\in A$ such that
$$\Ola\subseteq\cup_{j=1}^m{\rm graph}\,(\hl j).$$ Consider the map
$$P\mapto AAc{\prod_{j=1}^m(c+b_1l_j+\cdots+b_nl_j^n)}.$$
Then there are constants $C_0,\ldots,C_{m-1}\in A$ such that
$P(c)=c^m+C_{m-1}c^{m-1}+\cdots+C_0$. The argument
(without modification) used
in the proof of Theorem 4.2.2.1 shows that
there exists $\widetilde{b_0}\in A$ such that $P(\widetilde{b_0})\in G(A)$
and $\norm{\widetilde{b_0}-b_0}<\e$. Let $\tb(x)=\widetilde{b_0}+b_1x+
\cdots+b_nx^n$. Then
$$\norm{\tb(\xb)-b(\xb)}\le\norm{\tb(x)-b(x)}=\norm{\tilde{b_0}-b_0}<\e.$$
Moreover for each $\ol\in\Ola$ there exists $j\in\set{1,\ldots,m}$ such that
$\l=\hl j(\o)$. Since $P(\widetilde{b_0})$ is invertible,
$\widetilde{b_0}+b_1l_j+\cdots+b_nl_j^n$ is invertible and therefore
$$\ol(\tb(\xb))=\o(\tilde{b_0})+\cdots+\o(b_n)\l^n=\o(
\tilde{b_0}+\cdots+b_nl_j^n)\ne0.$$
So $\tb(\xb)\htt$ does not vanish on $\Ola$; that is $\tb(\xb)$ is
invertible in $\Ala$.\eop

\npar We end this section by strengthening Proposition 1. It was hoped that the
following auxiliary
construction could be iterated along with extensions, as in Section 1.4.
However there is an obstacle as we shall see.

\prop 5.3.3.6. Let $A$ be a regular, semisimple, commutative, unital Banach algebra and
$a\in G(A)$. There exists a norm parameter for $\Ala$ such that there is a unital, linear contraction
$T\map\Ala A$ with $T\circ\th=\id A$ where $\th$ is the canonical
embedding $A\to\Ala$. Moreover, $T$ is continuous with respect to the spectral
norms, and there exists a continuous extension, $\hT$, such that the following
diagram, in which the vertical arrows are the Gelfand transforms
and the horizontal arrows are surjections, is commutative:
$$\matrix{ \Ala &\mapright{T} & A\cr
\mapdown{} & &\mapdown{}\cr
\laA &\mapright{\hT} &\overline{A\htt}\cr}$$

\pf As remarked above, $e^x=a$ has local solutions in $A$ in view of the regularity of $A$,
so there exist pairs of open subsets of $\O$ and elements of $A$, $(V_1,l_1),\ldots,
(V_m,l_m)$, such that
$$e^{\hl j\vert_{V_j}}=\ha\vert_{V_j}\quad(j=1,\ldots,m)\qquad\hbox{ and }
\qquad \O=\bigcup_{j=1}^m V_j.$$

Again since $A$ is regular, there is a partition of unity $(u_j)_{j=1}^m\subseteq A$
subordinate to $(V_j)_{j=1}^m$.
In particular, since $A$ is semisimple, $\sum_{j=1}^m u_j=1$.

Choose a parameter, $t>0$, for the Lindberg extension such that $t\ge\sum_{j=1}^m
\norm{u_j}\norm{l_j}$ and $$t\ge\max_{j=1,\ldots,m}\norm{l_j}.$$ Set
$s_k:=\sum_{j=1}^mu_jl_j^k$ $(k=0,1,2,\ldots)$.

Let $T_0$ be defined from $A[x]$ to $A$ by
$$a(x)=\sum_{k=0}^n a_kx^k\mapsto \sum_{k=0}^n a_ks_k.$$
Then $T_0$ is a linear map and
$$\eqalign{ \norm{T_0(a(x))} & \le \sum_{k=0}^n \norm{a_k}\norm{s_k}\cr
&\le \sum_{k=0}^n \norm{a_k}\sum_{j=1}^m\norm{u_j}\norm{l_j}^k\cr
&\le \sum_{k=0}^n \norm{a_k}t.t^{k-1}=\norm{a(x)}.\cr}$$
Hence there is a continuous, norm-preserving, linear extension $T_1\map{A(x;t)}A$.
We see that it is given by $T_1(q(x))=\sum_{k=0}^{+\infty}q_ks_k$
where $q(x)=\sum_{k=0}^{+\infty}q_kx^k\in A(x;t)$.
Clearly
$T_1$ is surjective since for all $a_0\in A$ we have
$T_1(a_0)=(\sum_{j=1}^m u_j)a_0=a_0$. Since $T_0$ is a contraction,
so is $T_1$.

Now we check that $T_1(J)=\set0$. This is almost obvious for if
$q(x)=\sum_{k=0}^{+\infty}q_kx^k\in J=kh((e^x-a))$ then for all
$\ol\in h((e^x-a))$ we have $\ol(q(x))=0$.
(See Section 4.2.1 for the character space of $A(x;t)$.)
Since $\ol\in h((e^x-a))$ we have $e^\l=\ha(\o)$. Thus for all $\o\in\O$,
$$\eqalign{\o(T_1(q(x)) )&= \o\left( \sum_{k=0}^{+\infty}q_ks_k\right)\cr
&=\sum_{k=0}^{+\infty}\o(q_k)\o(s_k)\cr
&=\sum_{k=0}^{+\infty}\o(q_k)\sum_{j=1}^m\whu j(\o)\hl j(\o)^k\cr
&=\sum_{j=1}^m\whu j(\o)\sum_{k=0}^{+\infty}\o(q_k)\hl j(\o)^k.\cr}$$
Now if $\o\in V_j$ then $e^{\hl j(\o)}=\ha(\o)$ and
$\ab{\hl j(\o)}\le\norm{l_j}\le t$ so $(\o,\hl j(\o))\in h((e^x-a))$. If
$\o\not\in V_j$ then $\whu j(\o)=0$. Thus
$$\o(T_1(q(x)) )=\sum_{\scriptstyle j=1\atop\scriptstyle\whu j(\o)\ne0}^m\whu j(\o)
(\o,\hl j(\o))(q(x))=0.$$
Since $A$ is semisimple, $T_1(q(x))=0$ as required.

Thus $T_1$ induces a continuous, linear map $T\map{\Ala=A(x;t)/J}A$ such that
$T\circ\th=\id A$. It is a contraction with respect to the quotient norm on
$\Ala$ and therefore with respect to the final norm on $\Ala$, by Theorem 5.3.2.1(ii).

Now $\Ala$ is semisimple by Theorem 5.3.2.1 so we can also consider the
spectral radius norms (which we
shall denote by $\ab\cdot$ as usual) on $\Ala\htt$ and $\overline{A\htt}$ and
the induced map $S\map{(\Ala)\htt}{\overline{A\htt}}$. Let
$q(x)=\sum_{k=0}^n q_kx^k\in A[x]\ (q_0,\ldots,q_n\in A)$; elements
of the form $q(\xb)\htt$ are dense in $\Ala\htt$ and so in its
uniform completion, $\laA$. Then $S(q(\xb)\htt)=T(q(\xb))\htt=
\left(\sum_{k=0}^n q_ks_k\right)\htt$. Let $\o\in\O$. Then
$$\eqalign{\ab{ S( q(\xb)\htt )(\o)} &=
\ab{ \o\left(\sum_{k=0}^n q_ks_k\right)}\qquad\hbox{(by definition of $S$)}\cr
&= \ab{ \sum_{k=0}^n \o(q_k) \sum_{j=1}^m\whu j(\o)\hl j(\o)^k}
\qquad\hbox{(by definition of $s_k$)}\cr
&= \ab{ \sum_{\scriptstyle j=1\atop\scriptstyle\whu j(\o)\ne0}^m
\sum_{k=0}^n \o(q_k)\whu j(\o)\hl j(\o)^k} \cr
&= \ab{ \sum_{\scriptstyle j=1\atop\scriptstyle\whu j(\o)\ne0}^m
\whu j(\o)(\o,\hl j(\o))\Biggl( \sum_{k=0}^n q_k\xb^k\Biggr)}
\qquad\hbox{(since $(\o, \hl j(\o))\in\Ola$)} \cr
&= \ab{\sum_{\scriptstyle j=1\atop\scriptstyle\whu j(\o)\ne0}^m
\whu j(\o)q(\xb)\htt(\o,\hl j(\o))} \cr
&\le c\ab{q(\xb)}\cr}$$
where $c=\sum_{j=1}^m\ab{u_j}$. Hence $S$ is continuous and $\norm{S}\le c$.
(We also have $c\ge\ab{\sum_{j=1}^m u_j}=1$.) Therefore $S$ extends to a continuous map,
$\hT\map\laA{ \overline{A\htt}}$, with the required
properties.\eop

\npar The obstruction mentioned above is that no bound can be assumed on the elements of the partitions of
unity unless $A$ is uniformly dense in $C(\O)$.
This follows from a theorem of Bade and Curtis; see $\ref{Dal}$, p. 414.

\cor 5.3.3.7. Let $A$ be a regular, semisimple, commutative, unital Banach algebra.
Let $a\in G(A)$. If the norm-parameter for $\laA$
is sufficiently large then $\laA=C(\Ola)$ implies that $\overline{A\htt}=C(\O)$.

\pf Let $(V_j)_{j=1}^m$ be an open cover of $\O$ and $l_j\in A$
be a local solution of $e^x=a$ valid on $V_j$ ($j=1,\ldots,m$). Let $(u_j)_{j=1}^m
\subseteq A$ be a partition of unity subordinate to $(V_j)_{j=1}^m$, and choose
the norm-parameter, $t$ as in Proposition 6. By Proposition 2, $C(\Ola)=C(\O)^{\log a}$
and so we can form the continuous operator
$\hT\map{C(\Ola)}{C(\O)}$ as in Proposition 6. Writing $s_k=\sum_{j=1}^m u_jl_j^k$
($k\in\No$) we have for all $q_0,\ldots,q_n\in A$,
$$\hT\left( \left(\sum_{k=0}^n q_k\xb^k\right)\htt\right)=\left(\sum_{k=0}^n q_ks_k\right)\htt.$$
Since elements of the form $\left(\sum_{k=0}^n q_k\xb^k\right)\htt$ are dense in $\laA$
and $\hT$ is continuous this shows that $\hT(\laA)\subseteq \overline{A\htt}$.

We also have $\hT\circ\het=\id{C(\O)}$ and so for all $g\in\laA$ and $f\in C(\O)$,
$$\norm{f-\hT(g)}\le\norm\hT\norm{\het(f)-g}.$$
Hence if $\laA=C(\Ola)$ then $\hA$ is dense in $C(\O)$.\eop


\vskip 15pt
\noindent {\bf 5.3.4. Infinite Logarithmic Extensions and Systems of Extensions}\vskip 10pt


\noindent
We now show how to iterate logarithmic extensions. We shall assume from
now on that $A$ is regular and semisimple as well as commutative. We have shown
that these properties are preserved by simple extensions. It is also clear that if
$a_1,a_2\in G(A)$ then $\th(a_2)$ is invertible in $A_{\log a_1}$ where $\th$
is the embedding $A\to A_{\log a_1}$. Moreover if
$l$ is a local solution to $e^x=a_2$ then $\th(l)$ is a local solution
to $e^x=\th(a_2)$ in $A_{\log a_1}$ and we may use the
same norm parameter, $t_{a_2}$ for $A_{\log a_2}$ as in $(A_{\log a_1})_{\log a_2}$.
We shall from now on assume that for each $a\in G(A)$ a specific choice of
norm parameter $t_a>0$ has been made.

A further benefit of assuming regularity is that, as noted in Section 4.3.2, regular uniform algebras
are natural.

Recall that given a set, $S$, $S\fs$ denotes the set of all finite subsets
of $S$.

Let $\U=\set{a_1,\ldots,a_n}\in G(A)\fs$. In the next result we identify
$$A_{\log a_1,\ldots,\log a_n}:=(\cdots((A_{\log a_1})_{\log a_2})\cdots)_{\log a_n}.$$

We write $x^\nu$ as an abbreviation for $x_1^{\nu_1}\cdots x_n^{\nu_n}$ in
the following. Let $B(\U)=A(x_1,\ldots,x_n;t_{a_1},\ldots,t_{a_n})$ denote the following Banach algebra,
which appears in Example 2.1.18(v) of $\ref{Dal}$,
$$\set{\sum_{\nu\in\No^n}a_\nu x^{\nu}\in A[[x_1,\ldots,x_n]]\st
\norm{\sum_{\nu\in\No^n} a_\nu x^{\nu}}=\sum_{\nu\in\No^n}
\norm{a_\nu}t_{a_1}^{\nu_1}\cdots t_{a_n}^{\nu_n}<+\infty  \quad(a_\nu\in A)}$$
We may clearly view $B(\U)$ as the completion of the
algebra of polynomials $A[x_1,\ldots,x_n]$.
Let $J(\U)$ be the ideal $kh((e^{x_1}-a_1,\ldots,e^{x_n}-a_n))$ of
$B(\U)$ and let $\AlU=B(\U)/J(\U)$. It is easy to see, using induction,
and writing $\approx$ for a homeomorphism, that
$$\Om{B(\U)}\approx\O\times{\scriptstyle n\atop{\displaystyle\times\atop{j=1}} }B_\C[0,t_{a_j}]$$
where for $\ol\in\O\times{\scriptstyle n\atop{\displaystyle\times\atop{j=1}} }B_\C[0,t_{a_j}]$,
$$\ol\left( \sum_{\nu\in\No^n} a_\nu x_1^{\nu_1}\cdots x_n^{\nu_n}\right)=
\sum_{\nu\in\No^n} \o(a_\nu) \l_1^{\nu_1}\cdots\l_n^{\nu_n}
\qquad\left(\sum_{\nu\in\No} a_\nu x_1^{\nu_1}\cdots x_n^{\nu_n}\in B(\U)\right).$$
Hence
$$\Om{\AlU}\approx h(J(\U))\approx\set
{\ol\in\O\times\C^n\st \ab{\l_j}\le t_{a_j},\ e^{\l_j}=\ha_j(\o)\ (j=1,\ldots,n)}.$$

\lem 5.3.4.1. Let $A$ and $\U\in G(A)\fs$ be as above. There exists a (unique)
topological isomorphism $\phi\map {\AlU}{A_{\log a_1,\ldots,\log a_n}}$ such
that $\phi(\xb_j)=\yb_j\ (j=1,\ldots,n)$ and $\phi(J(\U)+a)=a\ (a\in A)$
where $\yb_j$ is the coset of the indeterminate, $y_j$, in the $j$th Lindberg
extension $A_{\log a_1,\ldots,\log a_j}$ identified with its canonical image
in $A_{\log a_1,\ldots,\log a_n}$

\pf Uniqueness is obvious.

Let
$$\phi^0\mapto{A[x_1,\ldots,x_n]} { A_{\log a_1,\ldots,\log a_n}}
{\sum\nolimits_\nu a_\nu x^\nu}{\sum\nolimits_\nu a_\nu \yb_1^{\nu_1}\cdots\yb_n^{\nu_n} }.$$
Then $\phi^0$ is an evaluation homomorphism and for all $(a_\nu)_{\nu\in\No^n}\subseteq A$
with finitely-many non-zero $a_\nu$, we have
$$\eqalign{ \norm{ \sum\nolimits_\nu a_\nu \yb_1^{\nu_1}\cdots\yb_n^{\nu_n} }
&\le\sum\nolimits_\nu
\norm{ a_\nu \yb_1^{\nu_1}\cdots\yb_{n-1}
^{\nu_{n-1}} }t_{a_n}^{\nu_n}\le\cdots\cr
&\le\sum\nolimits_\nu\norm{ a_\nu}t_{a_1}^{\nu_1}\cdots t_{a_n}^{\nu_n}\cr
&=\norm {\sum\nolimits_\nu a_\nu x_1^{\nu_1}\cdots x_n^{\nu_n}}.\cr}$$
Thus $\phi^0$ extends to a contraction $\phi^1\map {B(\U)}{ A_{\log a_1,\ldots,\log a_n}}$.

We now assume that $n=2$ in the rest of this proof for notational simplicity.
The general case can be proved similarly. Observe
that there is a homeomorphism
$$\eta\mapto{\O(A_{\log \U})}{\O(A_{\log a_1,\log a_2})}{(\o,(\l_1,\l_2))}{((\o,\l_1),\l_2)}.$$
Let $\b=\sum_{j,k\in\No}b_{j,k}x_1^jx_2^k\in B(\U)$. Since $A_{\log a_1,\log a_2}$
is semisimple we have $\b\in\ker\phi^1$ if and only if for all $\ol\in\O(A_{\log \U})$,
$\eta\ol\left( \sum_{j,k\in\No}b_{j,k}\yb_1^j\yb_2^k\right)=0$. The left-hand side is
equal to
$$\eqalign{((\o,\l_1),\l_2)\left(\sum_{k\in\No}\left(\sum_{j\in\No}b_{j,k}\yb_1^j\right)\yb_2^k\right)
&= \sum_{k\in\No}(\o,\l_1)\left(\sum_{j\in\No}b_{j,k}\yb_1^j\right)\l_2^k\cr
&= \sum_{k\in\No}\sum_{j\in\No}\o(b_{j,k})\l_1^j\l_2^k\cr
&= (\o, (\l_1,\l_2))\left(\sum_{j,k\in\No}b_{j,k}x_1^jx_2^k\right).\cr}$$
So $\b\in\ker\phi^1$ if and only if for all $\ol\in h(J(\U))$ we have $\ol(\b)=0$.
That is, $\ker\phi^1=kh(J(\U))$. It is standard (see, for example, $\ref{Dal}$, p. 66)
that $khk=k$ and so $kh(J(\U))=J(\U)$.

So $\ker \phi^1= J(\U)$ and $\phi^1$ induces an isomorphism $\phi\map{\AlU}
{A_{\log a_1,\ldots,\log a_n}}$ of norm one. The map is a homeomorphism
by Banach's isomorphism theorem.\eop

\npar This lemma will be useful in allowing us to identify the uniform algebra
$\lUA$ in the following. We now define uniform-algebra extensions of $A$, a regular
uniform algebra on $\O$, in which every element of a set $\U\subseteq G(A)$
has a logarithm. The construction closely follows the treatment of
Cole's construction ($\ref{Col}$) which was given in Section 1.2.

When $\U$ is finite, we shall use the notation $\AlU$ for the extension
$B(\U)/J(\U)$
defined before Lemma 1 under a norm such that the canonical map $\th_\U\mapto A{\AlU}a{J(\U)+a}$
is an isometry. (Such a norm exists as we can regard the
algebra as an extension obtained by applying Lindberg's construction
finitely-many times.) We set $\AlU=A$ if $\U=\emptyset$.

Now let $\U\subseteq G(A)$ and set
$$\OlU:=\set
{\ol\in\O\times\C^\U\st \ab{\l_a}\le t_{a}\hbox{ and }\ e^{\l_a}=\ha(\o)\ \hbox{for all }a\in\U}
.$$
Then $\OlU$ is clearly compact. Let the canonical projections be
$$\eqalign{
&\pi\mapto\OlU\O\ol\o;\cr
&p_a\mapto\OlU\O\ol{\l_a}\qquad(a\in\U).\cr}$$
We set $\lUA$ to be the closed subalgebra of $C(\OlU)$ generated by $\ps(A)\cup
\set{p_a\st a\in\U}$. Evidently $\lUA$ is a uniform algebra on $\OlU$ and
can be described as a direct limit of a system of Lindberg extensions. We make this
statement precise in the next lemma. In order to do this it is useful to introduce
some more notation. For $\U\subseteq \V\in G(A)\fs$ with
$\U=\set{a_1,\ldots,a_m}$ and $\V=\U\udot\set{a_{m+1},\ldots,a_n}$,
we have
$\OlU=\Om{\AlU}$ by Lemma 1 and the preceeding comments so there is a canonical surjection
$$\pp\U\V\mapto{\O_{\log\V}}\OlU{(\o,(\l_j)_{j=1}^n)}{(\o,(\l_j)_{j=1}^m)}.$$
This induces
an isometric monomorphism
$$\pps\U\V\mapto{C(\OlU)}{C(\O_{\log\V})}f{f\circ\pp\U\V}.$$

\lem 5.3.4.2. Let $(A,\O)$ be a regular uniform algebra and $S\subseteq G(A)$.
Then
$$\eqalign{\OlS&=\il{ } \left( (\OlU),(\pp\U\V)_{\U\subseteq\V}\st\U,\V\in S\fs\right)
\qquad\hbox{ and}\cr
\lSA&=\dl{ } \left( (\lUA),(\pps\U\V)_{\U\subseteq\V}\st\U,\V\in S\fs\right).\cr}$$
(We have omitted the maps on the left-hand side for clarity; they will be given in
the course of the proof.)

\pf The justification is routine but we include the details for completeness and
to establish the notation.

Let the inverse limit of the compact, Hausdorff spaces be denoted by $K$. It
is standard that $K$ is non-empty and can be constructed as follows:
$$K=\set{ \k\in\times_{\U\in S\fs}\OlU\st\hbox{ for all }\U\subseteq\V\in S\fs
\ \pp\U\V(\k_\V)=\k_\U}.$$
The canonical maps associated with $K$ are the restrictions of the coordinate projections
to $K$
$$\tpp\U S\mapto K\OlU\k{\k_\U}\qquad\qquad(\U\in S\fs).$$
The canonical maps associated with the inverse limit in the statement are given by
$$\pp\U S\mapto \OlS\OlU\ol{(\o,(\l_a)_{a\in\U}})\qquad\qquad(\U\in S\fs).$$

Let
$$h\mapto\OlS K\ol{((\o,(\l_a)_{a\in\U}))_{\U\in S\fs}}.$$
Clearly $h$ is well-defined. The map is also continuous for if $((\o_t,\l^t))$ is
a net in $\OlS$ converging to $\ol$ then
$$\o_t\to\o\qquad\hbox{ and }\qquad\l^t_a\to\l_a\qquad(a\in S).$$
So for each $\U\in S\fs$ we have
$$(\o_t,(\l^t_a)_{a\in\U})\to (\o,(\l_a)_{a\in\U}).$$
By the definition of the product topology this implies that $h((\o_t,\l^t))\to
h(\ol)$. The map is clearly a bijection. Since $\OlS$ is compact and $K$ is Hausdorff,
$h$ is a homeomorphism. Finally for all $\uis$ the following diagram is commutative:
$$\matrix{\OlS &\mapright h &K\cr
\mapdown{\pp\U S}&\mapsw{\tpp\U K} &\cr
\OlU &&\cr}$$

Now we check that $\lSA$ is the direct limit of
the direct system above.
First note that this system is well-defined for if $\U\subseteq\V\in S\fs$ then
$\pps\U\V(\lUA)\subseteq A^{\log \V}$. To see this, it is clear from Lemma 1 that
$$\lUA=\overline{A_{\log \U}\htt},$$
as both algebras are equal to the closure in $C(\OlU)$ of the
algebra of polynomials in $\set{p_{a,\U}\st a\in\U}$ over $\pps\emptyset\U(A)$
where for all $a\in\U$, $p_{a,\U}$ is the coordinate map
$$p_{a,\U}\mapto\OlU\C\ol{\l_a}.$$
Let $f=\sum_\nu \pps\emptyset\U(b_\nu)p_{a_1,\U}^{\nu_1}\cdots p_{a_n,\U}^{\nu_n}
\in A_{\log \U}\htt$ where $(b_\nu)\subseteq A$, finitely many $b_\nu$
are non-zero, $\U=\set{a_1,\ldots,a_n}$,
$\V=\set{a_1,\ldots,a_{n\pri}}$, and $n\pri\ge n$. Then
$$\pps\U\V(f)=\sum_\nu \pps\emptyset\V(b_\nu)
\pps\U\V(p_{a_1})^{\nu_1}\cdots \pps\U\V(p_{a_n})^{\nu_n}.$$
It is easy to check that $p_{a,\V}=\pps\U\V(p_{a,\U})$. Therefore $\pps\U\V(f)\in
A_{\log\V}\htt$ and this implies
$$\pps\U\V(\lUA)\subseteq\overline{A_{\log\V}\htt}=A^{\log \V}\qquad\hbox{ as required.}$$

By results in Appendix 2, it is enough to show that $\lSA$
is generated  by the subalgebra $D=\bigcup_{\uis}\pps\U S(\lUA)$.

Again let $f=\sum_\nu \pps\emptyset\U(b_\nu)p_{a_1,\U}^{\nu_1}\cdots p_{a_n,\U}^{\nu_n}
\in A_{\log \U}\htt$ in the notation as above. It is easily checked that
$\pps\emptyset S=\ps$ and that for $a\in \U$ we have $\pps\U S(p_{a,\U})=p_a\in\lSA$.
So by its definition, $\lSA$ must be a subset of $\overline {D}$.
On the other hand we have shown that $D\subseteq\lSA$ so
the two uniform algebras are equal.\eop

\cor 5.3.4.3. Let $(A,\O)$ be a regular uniform algebra and $S\subseteq G(A)$. Then
\item{(i)} $\lSA$ is regular, and
\item{(ii)} if $G(A)$ is dense in $A$ then $G(\lSA)$ is dense in $\lSA$.

\pf (i) This follows from Proposition 5.3.3.3, Lemma 5.3.4.2 and the fact that the direct limit
of a direct system of regular uniform algebras is again regular.
(This is proved in essence in Proposition 2.2.1.2.)

(ii) Using Proposition 5.3.3.5 and induction, it follows from Lemma 5.3.4.1 that for
all $\U\in S\fs$, $A_{\log\U}$ has dense invertible group. The result now
follows from Lemma 5.3.4.2 and Proposition 4.2.2.5.\eop

\npar We now proceed as in Section 1.2.3 and repeat these extensions by
induction. As the details are so similar, we do not produce them here.

Let $(A,\O)$ be a regular uniform algebra and $\u>0$ be an ordinal number.
The direct system of uniform algebras
$$\left( (A_\t,\O_\t)_{\t\le\u},(\pps\s\t)_{\s\le\t\le\u}\right),$$
in which for all $\s\le\t\le\u$, $\Ps\s\t$ is a continuous surjection $\O_\t
\to\O_\s$ and $\pss\s\t\colon f\mapsto f\circ\Ps\s\t$,
is called a {\it system of logarithmic extensions} of $(A,\O)$ provided that
$$
(A_\t,\O_\t)=\cases{ (A,\O) &if $\t=0$,\cr
(A_\s^{\log S_\s}, \O_{\s,\log S_\s})
&for some $S_\s\subseteq G(A_\s)$, if $\t=\s+1$,\cr
\dl{} \left( (A_\s,\O_\s)_{\s<\t},(\pps\r\s)_{\r\le\s<\t}\right)&if $\t>0$
is a limit ordinal.\cr}$$

In the case where $\t$ is a limit ordinal, $A_\t$ is a uniform algebra on
$$\O_\t:=\il{}\left( (\O_\s)_{\s<\t},(\pp\r\s)_{\r\le\s<\t}\right),$$ and is the
closure of $\bigcup_{\s<\t}\pps\s\t(A_\s)$ (see Appendix 2). As in Corollary 3,
the new algebras generated in this way are all regular. Similarly
if $A$ has dense invertible group then so has $A_\t\;(\t\ge0)$.

The main application of these systems of extensions will be in Theorem 5.3.4.5. Before
stating this we give the following result.

\lem 5.3.4.4. Let $(A_\t,\O_\t)_{\t\le\u}$ be a system of logarithmic extension of
$A_0=C(\O_0)$. Then for all $\t\in[0,\u]$ $A_\t=C(\O_\t)$.

\pf As in Proposition 1.4.2.3 we have
that direct limits of trivial uniform algebras are trivial.
This lemma will therefore follow from the transfinite induction theorem and
Proposition 5.3.3.2, once
we have checked the case of extensions at a non-limit ordinal.

Let $A=C(\O)$ and $S\subseteq G(A)$. We must show $\lSA=C(\OlS)$. By Lemma 5.3.4.2 and
the result for direct limits again,
it is sufficient to prove that for all $\U\in S\fs$, $\lUA$ is self-adjoint.
So let $\U\in S\fs$. As noted before, $\lUA=\overline{\AlU\htt}$ and it is
clearly enough to prove that $\bar p_{a,\U}\in \AlU\htt$ for each $a\in\U$ where
$p_{a,\U}$ is as in Lemma 5.3.4.2.
A similar argument to the one used in Proposition 5.3.3.2 now gives the result.\eop

\thm 5.3.4.5. Let $A$ be a regular, semisimple, commutative, unital Banach algebra
with maximal ideal space $\O$.
Then there exists a regular uniform algebra, $A_\oo$, with maximal ideal space
$\O_\oo$ such that
\item{(i)} there is a continuous surjection $\pi\,\colon\,\O_\oo\to\O$,
\item{(ii)} $\pi$ induces a continuous embedding $\ps\colon A\hookrightarrow A_\oo;f\mapsto f\circ\pi$,
\item{(iii)} $H^1(\O_\oo,\Z)=0$, and
\item{(iv)} $\O_\oo$ is metrizable if $\O$ is.

If $A$ is a uniform algebra, then $\ps$ is isometric.

\pf If $A$ is not a uniform algebra then it embeds continuously in a regular
uniform algebra on $\O$ (see Appendix 1), so we may assume
to start with that $A$ is a regular uniform algebra.

Let $(A_\oo,\O_\oo)$ be the final algebra in a system of logarithmic extensions
of $A$ as above in which $S_n$ is taken to be a dense subset of
$G(A_n)$ for all $n<\oo$.

Item (iv) is as in the case of Cole extensions.
We can assume that $S_n$ is countable for each
$n<\oo$ if $\O_\oo$ is metrizable (equivalently, $A$ is separable); the
logarithmic extension $A_{n+1}$ will then be separable if $A_n$ is. The result now
follows from the easily-checked fact that an inverse limit of a sequence of metrizable, compact
spaces is metrizable.

It remains to check that $G(A_\oo)=e^{A_\oo}$.
Let $f\in G(A_\oo)$. Since $e^{A_\oo}$ is open (see Theorem 10.43
of $\ref{RudFA}$), there exists $\e>0$ such that $B_{A_\oo}(f,\e)
\subseteq e^{A_\oo}f$. Now $\oo$ is a limit ordinal so $\cup_{n<\oo}\pps n\oo(A_n)$
is dense in $A_\oo$. Hence there exists $n<\oo$ and $g_n\in A_n$ such that
$\norm{\pps n\oo(g_n)-f}<\e/2$. So $\pps n\oo(g_n)\in e^{A_\oo}f$.
In particular $g_n\in G(A_n)$ and so there exists $g\in S_n$ with
$\norm{g_n-g}=\norm{\pps n\oo(g_n)-\pps n\oo(g)}<\e/2$. Thus there exists
$h\in A_\oo$ such that $\pps n\oo(g)e^h=f$.

By construction there is $p\in A_{n+1}$ with $\pps n{{n+1}}(g)=e^p$. Since
$\Ps{n+1}{\oo}$ is continuous, $\pps {{n+1}}{\oo}(e^p)=\exp(\pps {{n+1}}{\oo}(p))$ so
we have $f=\exp(h+\pps {{n+1}}{\oo}(p))\in e^{A_\oo}$ as required.\eop

\npar So if $A$ is a non-trivial, regular uniform algebra with dense
invertible group (such as McKissick's example) then $A_\oo$ is a regular
uniform algebra with dense exponential group. We would like to know if
$A_\oo$ is trivial or not in these cases.

\def\na{{\bf na}}							
\def\mapse#1{\searrow${$\scriptstyle#1$}$}
\def\mapsw#1{ \swarrow${$\scriptstyle{#1}$}$  }
\def\dl#1{{\lim_\rightharpoondown}\mathstrut_{#1}}	
\def\il#1{{\lim_\leftharpoondown}\mathstrut_{#1}}	


\vskip 15pt
\noindent {\medtenbf 5.4. Conclusion}
\vskip 15pt

\noindent Sections 5.1 and 5.2 show that many sets of uniform algebras can not
be counterexamples to
the conjecture. But a counterexample could still exist.

Ramsay and Hoffman's example
($\ref{HofRam}$) of a non-trivial uniform algebra on $\b\N$, the Stone-\v Cech compactification
of $\N$, in which every non-negative element of $C(\b\N)$ is the modulus of some
element of the algebra seemed promising but does not appear to have
dense invertible group. Perhaps some other subalgebra of $C(\b\N)$ can
provide a counterexample. Even if one were found, this would leave open the question of
whether $\overline{e^A}=A$ implies that $A=C(\O)$ if $A$ is natural and $\O$ is metrizable.

It is extremely unsatisfactory that we do not know if systems of logarithmic extensions
preserve non-triviality. Numerical methods could conceivably be developed to test whether or not
the algebra constructed at the end of Section 5.3 is trivial or not.

None of the familiar characterisations of $C(X)$ (see for example $\ref{Bur}$ and $\ref{Sto}$)
seems to yield a simple reason why the conjecture should be true.

There are many more examples to check, and the conjecture leads to several other intriguing
problems such as whether or not $\overline{e^A}=A$ implies that $\S A=\O$.

\vfill\eject


\noindent {\bigtenbf Appendix 1}
\vskip 10pt
\noindent {\bigtenrm Gelfand Theory}
\vskip 25pt

\def\hp{\hat\phi}
\def\inc{\iota}
\def\is{\inc^*}
\def\iss{\inc^{**}}
\def\snc#1{\overline{\left(#1\right)\htt}}

\vskip 10pt

\noindent The purpose of this section is to clarify our notation and terminology associated
with the Gelfand theory of normed algebras.
A good account of Gelfand theory is given in Chapter 3 of $\ref{Pal}$.
We note some elementary results which will be used at
various points in thesis. As in the main body of the thesis, $A$ denotes a
commutative, unital, complex, normed algebra.

Recall that a {\it character} of an algebra, $A$, is a unital homomorphism $A\to\C$.
Let $\OA$ denote the space of continuous characters of $A$; when
$\O$ appears on its own it refers to $A$. Thus $\O\subseteq A^*$, the
continuous, linear maps $A\to\C$. 

Recall that the functions $(A^*\to\C\;;\;\l\mapsto\l(a))_{a\in A}$
induce a topology on $A^*$, called the weak *-topology.
As discussed in $\ref{Are}$, the space $\O$, with the relative weak *-topology, generalises
the notion of the character space of a Banach algebra. There is a bijection between
the continuous characters of $A$ and the closed, maximal ideals of $A$. If $A$ is
complete then characters are automatically continuous, of norm one, and the maximal ideals are closed.

It is easy to
check (see Lemma A.1.4) that $\O$ is homeomorphic to $\Om{\c A}$, the maximal ideal space
of the completion of $A$. In particular $\O$ is compact and Hausdorff. We briefly recall
some standard properties of the completion of a normed algebra.

A {\it completion} of $A$ is a pair
$(\tilde A, \inc)$ where $\Ac$ is a Banach algebra and $\inc$ is an isometric
monomorphism $A\to\Ac$ with dense image.\par
It is well-known (see for example $\ref{Bol}$, p. 35)
that completions exist and
are unique up to isometric isomorphism.
By this we mean that if $(B, j)$ is another completion then there is an
isometric isomorphism $\phi\map {\tilde A} B$ such that the following is commutative:
$$\matrix{ \tilde A&\mapright\phi&B\cr
	\mapup \inc&\mapne j&\cr
	A&&\cr}$$

Recall that the Gelfand transform of an element $a\in A$ is defined by
$$\ha\mapto\O\C\o{\o(a)}$$
and the map sending $a$ to $\ha$ is
a homomorphism, $\Gamma$, of $A$ into the algebra, $C(\O)$.
We denote the image of $\Gamma$ by $\hA$. We call $\hA$ {\it symmetric}
if $\hA$ contains the pointwise complex-conjugate of each of
its elements. Unless otherwise stated, we shall
regard $\hA$ as an algebra with the {\it supremum norm}
$$\norm\ha:=\sup_{\o\in\O}\ab{\ha(\o)}\qquad\qquad(\ha\in\hA).$$

\dfn A.1.1 ($\ref{Are}$). The normed algebra $A$ is called {\it topologically semisimple}
if $\Gamma$ is injective.

\npar 
As noted in $\ref{Are}$, if $A$ is a Banach algebra then topological semisimplicity is equivalent to
the usual notion of semisimplicity. We introduce some further standard notation for $a\in A$:
$$\eqalign{
\s(a)&:=\im\ha\hbox{ is the {\it spectrum} of }a\hbox{ and,}\cr
\ab a&:=\sup_{\o\in\O}\ab{\ha(\o)}\hbox{ is the {\it spectral radius} of }a.\cr}$$
It is a standard fact (see $\ref{Pal}$, p. 212),
known as the {\it spectral radius formula}, that, provided
$\ab a\le\norm a$ for each $a\in A$,
$$\ab a=\lim_{n\to+\infty}\norm{a^n}^{1/n}\qquad\qquad(a\in A).$$

We now state some results about the Gelfand
transform. The following five results are all surely well-known.
Recall that a continuous, unital homomorphism, $\th\map AB$, of normed algebras induces the
adjoint
$$\theta^*\mapto{\O(B)}{\OA}\o{\o\circ\theta}.$$

\lem A.1.2. Let

	$$\matrix{ A_1 &\mapright{\phi_1} &B_1 \cr
	\mapdown{\ths 1} & &\mapdown{\ths 2} \cr
	A_2 &\mapright{\phi_2} &B_2 \cr } $$

be a commutative diagram of continuous, unital homomorphisms. Then the following diagram
commutes:

$$\matrix{ \O\bigl(A_1\bigr) &\mapleft{\phi_1^*} &\O\bigl(B_1\bigr) \cr
	\mapup{\ths 1^*} & &\mapup{\ths 2^*} \cr
	\O\bigl(A_2\bigr) &\mapleft{\phi_2^*} &\O\bigl(B_2\bigr) \cr } $$
	\vskip 2pt
\pf We have $\phi_1^*\circ\ths 2^*=(\ths 2\circ\phi_1)^*=(\phi_2\circ\ths 1)^*
=\ths 1^*\circ\phi_2^*$.\eop

\lem A.1.3. Let $\phi\map A B$ be a unital, isometric isomorphism. Then there is an
isometric isomorphism $\hp\map\hA{\hat B}$ such that the following is commutative:

	$$\matrix{ A &\mapright{\phi} &B \cr
	\mapdown{} & &\mapdown{} \cr
	\hA &\mapright{\hp} &\hat B \cr } $$
where the vertical arrows are the Gelfand transforms.

\pf In this case it is known that $\f$ induces a homeomorphism
$\f^*\mapto{\O(B)}{\OA}\o{\o\circ\f}$
and an isometric isomorphism $\f^{**}\mapto{C\left(\OA\right)}{C\left(\O(B)\right)}f{f\circ\f^*}$.
Clearly $\phi^{**}(\hA)\subseteq\hat B$.
The map $\hp=\f^{**}\vert_{\hA}$ is surjective because for $b\in B$, let $a\in A$ be such that
$b=\f(a)$. We then have, for all $\o\in\O(B)$,
$$\hat{b}(\o)=\o(b)=\o(\f(a))=(\o\circ\f)(a)
=\ha(\o\circ\f)=\ha\bigl(\f^*(\o)\bigr)=\bigl(\f^{**}(\ha)\bigr)(\o).$$
So $\hat b=\f^{**}(\ha)$.\eop

\lem A.1.4. Let $B$ be a normed algebra and $\inc\map AB$ be an isometric monomorphism with dense image.
The adjoint maps induced between the spaces of closed, maximal ideals by
$\G\map{A}{\hA}$ and $\inc\map{A}{B}$
are homeomorphisms.

\pf
The first assertion is well-known but to fix notation we provide a proof. Remember that
the continuous map $\e\mapto{\OA}{\O(\hA)}{\o}{\eh\o}$ gives
all the evaluation characters on $\hA$. It is a formal exercise to check that
$\e\circ\G^*=\id{\O(\hA)}$ and $\G^*\circ\e=\id\OA$.

Since $\im \inc$ is dense in $B$, $\inc^*$ is injective. Let $\o\in\OA$ and
$\th\inv$ be the inverse map $\im \inc\to A$. Then $\o\circ \th$ extends to a
continuous homomorphism $\tilde\o\map B\C$ and $\inc^*(\tilde \o)=\o$.\eop

\lem A.1.5. Let $B$ be a commutative, unital normed algebra. Then, up to isometric
isomorphism, $\snc{\c B}=\overline{B\htt}$. (The closures are taken in $C(\Om B)$ with respect to
the supremum norms.)

\pf This is elementary; we provide a proof for completeness. By Lemma 4
the embedding $\inc\map B{\c B}$ induces a homeomorphism
$$\is\map{\Om{\c B} }{\Om B},$$
and so the adjoint map
$$\iss\mapto{ C(\Om{B}) }{ C(\Om{\c B}) }f{f\circ\is}$$
is an isometric isomorphism. We show that $\iss\left(\overline{B\htt} \right)=\snc {\c B}$.

Let $b\in B$ and $\o\in\Om{\c B}$. Then
$$\iss(\hat b)(\o)=\hat b(\is(\o))=\is(\o)(b)=\o(\inc(b))=\widehat{\inc(b)}(\o).$$
This shows that $\iss(B\htt)\subseteq \c B\htt$ and so $\iss\left(\overline{B\htt}\right)
\subseteq\snc{\c B}$.

Since $\iss$ is isometric, $\iss\left(\overline{B\htt}\right)$ is closed in $\snc{\c B}$
and it is sufficient to prove that $\iss(B\htt)$ is dense in  $\snc{\c B}$.

Let $c\in\left(\c B\right)\htt$. Then there exists $(b_n)\subseteq B$ such that $\widehat{\inc(b_n)}
=\iss(\widehat{b_n})
\to c\ (n\to+\infty)$.\eop

\noindent We end this appendix with a result which is used in Chapter 4. The standard unitisation,
$B_1$,
of a normed algebra, $B$, is defined in Appendix 3.

\prop A.1.6. Let $M$ be a closed, maximal ideal of a normed algebra $A$. Then the map
$$\th\mapto{M_1}A{(a,\l)}{a+\l 1_A}$$
is a topological isomorphism.

\pf The map is clearly a homomorphism. It follows from the definition of the norm on
$M_1$ that $\th$ is a contraction. Let $\o$ be the character associated with $M$ in the
usual way. Then $\o$ is continuous since $M$ is closed. Indeed $\norm\o=1$ since
$\o$ has an extension to a character on the completion of $A$.

For each $a\in A$, we have $a-\o(a)1_A\in M$ so $a=\th(a-\o(a)1_A,\o(a))$. This
shows that $\th$ is surjective. Now let $(a,\l)\in\ker \th$. Thus $a+\l 1_A=0$
and so $\o(a)+\l=0$. But $M=\ker \o$ so $\l=0$ and therefore $a=0$. Therefore
$(a,\l)=0$ and $\th$ is injective.

It remains to prove that $\th\inv$ is continuous. By a calculation above we have
that for all $a\in A$, $\th\inv(a)=(a-\o(a)1_A,\o(a))$. Thus
$$\eqalign{\norm{\th\inv(a)}&=\norm{a-\o(a)1_A}+\ab{\o(a)}\cr
&\le\norm a+\norm{\o(a)1_A}+\norm a\cr
&\le 3\norm a.\cr}$$
Therefore $\th\inv$ is continuous (and has norm at most 3).\eop

\def\snc#1{\overline{#1\htt}}	

\vfill\eject


\noindent {\bigtenbf Appendix 2}
\vskip 10pt
\noindent {\bigtenrm Direct and Inverse Limits}
\vskip 25pt

\def\inc{\iota}
\def\is{\inc^*}
\def\iss{\inc^{**}}
\def\ii{{i\in I}}

\vskip 10pt

\noindent We assume some fluency with direct and inverse systems. A good introduction
to the subject is in $\ref{Rot}$, while Chapter 1 of $\ref{Pal}$ gives useful
details of the special cases of normed algebras and normed spaces. In this section we
give details especially relevant to the the applications in this thesis.

Throughout this appendix, $I$ denotes
a directed set; this means that there exists a relation, $\le$, on $I$ such
that for all $i,j,k\in I$,
\item{(i)} $i\le i$;
\item{(ii)} $i\le j\le k\hbox{ implies that }i\le k$;
\item{(iii)} there exists $l\in I\hbox{ such that }l\ge i\hbox{ and }l\ge j$.

\dfn A.2.1. The family of normed algebras and homomorphisms
$$((A_i)_{i\in I},(\tss ij)_{i\le j})\eqno{(1)}$$
is a {\it directed system} if
for all $i,j\in I$,
$\tss ij\map{A_i}{A_j}$ is a homomorphism, called a {\it connecting homomorphism},
of normed algebras ($i\le j$) and
$$\eqalign{\hbox{for all }i,j,k\in I\qquad \tss ik&=\tss jk\circ\tss ij\qquad(i\le j\le k)\cr
\hbox{and}\qquad\tss ii&=\id{A_i}.\cr}$$

\npar
For our purposes
(except Lemma 3.2.3.2), we can further assume that each of the connecting homomorphisms is isometric.

\dfn A.2.2. The {\it normed direct limit} of (1) is a pair,
$$\dl{}((A_i)_{i\in I},(\tss ij)_{i\le j}):=(A,(\ths i)_{i\in I}),$$
in $\na$ consisting of a normed
algebra, $A$, and homomorphisms $\ths i\map {A_i}A$ such that
$$\hbox{for all }i,j\in I,\ i\le j\qquad \ths i=\ths j\circ \tss ij.$$
The pair is also required to satisfy the universal property that whenever
$B$ is a normed algebra and $(\f_i\map{A_i}B)_{i\in I}$ is a family of homomorphisms
such that for all $i,j\in I$ with $i\le j$ we have $\f_i=\f_j\circ\tss ij$ then
there exists a unique homomorphism $\f\map AB$ such that for all $i\in I$, $\f\circ\ths i=\f_i$.
This universal property makes the direct limit unique up to equivalence in the
appropriate category.

\npar
We now show how normed direct limits can be constructed (see $\ref{Pal}$, Chapter 1).
Assuming without loss of generality that the algebras are pairwise disjoint, then
$A$ may be taken as
$$\left(\,\udot_{i\in I}A_i\right)/\sim$$
where $\sim$ is the equivalence relation given by
$$a_i\sim a_j\quad\hbox{if and only if there exists }k\ge i,j\hbox{ such that }
\tss ik(a_i)=\tss jk(a_j),$$
for $a_i\in A_i$ and $a_j\in A_j$.
Note that if $j\ge i$, then, since $\tss ij$ is injective,
$$a_i\sim a_j\hbox{ if and only if }\ \tss ij(a_i)=a_j$$
for all $a_i\in A_i, a_j\in A_j$.
Let square brackets denote the equivalence classes of $A$.
Algebraic operations and the norm on $A$ are given by
$${\eqalign{[a_i]+\l[a_j]&:=[\tss ik(a_i)+\l\tss jk(a_j)];\cr
[a_i][a_j] &:=[\tss ik(a_i)\tss jk(a_j)];\cr
\norm{[a_i]} &:=\norm{a_i}.\cr}}
\qquad (i,j,k\in I\ (k\ge i,j),\
a_i\in A_i,a_j\in A_j,\l\in\C)$$
The canonical maps are given by $\ths i\mapto{A_i}Aa{[a]}$.
The norm is well-defined because all maps are isometric.

In the more general case where the connecting homomorphisms are contractions
then a construction for the direct limit, $A$, is as follows (see also $\ref{Pal}$). First
note that there is a seminorm, $p$, on the algebra $A_\sim:=\left(\,\udot_{i\in I}A_i\right)/\sim$
given before defined by
$$\eqalign{p([a_i])&:=\limsup_{j\ge i}\norm{\tss ij(a_i)}\qquad\qquad(i\in I,\ a_i\in A_i)\cr
&=\inf_{j\ge i}\sup_{k\ge j}\norm{\tss ik(a_i)}.\cr}$$
It is standard (see, for example, $\ref{RudFA}$, Section 1.43) that $p\inv(0)$
is an ideal in $A_\sim$ and that the quotient algebra $A_\sim/p\inv(0)$ has a norm
given by
$$\norm{p\inv(0)+b}:=p(b)\qquad\qquad(b\in A_\sim).$$
We define the normed direct limit of (1) to be $A=A_\sim/p\inv(0)$ in this case. The
canonical maps are compositions of the natural maps
$$\theta_i\mapto{A_i}A{a_i}{p\inv(0)+[a_i]}\qquad\qquad(i\in I).$$

In $\Ba$, the direct limit of (1) is defined to be the completion of the normed-algebra
direct limit.

We now draw the reader's attention to two special cases.

The first applies to $\na$ when all the connecting homomorphisms are given by inclusions.
The reader can easily check that the direct limit is identifiable with $\cup_{i\in I}A_i$
and inclusion maps.

The second case is more complicated, as it involves two systems of objects. It is
now convenient to recall the following definition.

\dfn A.2.3. An {\it inverse system} of topological spaces is a family
$$((X_i)_{i\in I},(\Ps ij)_{i\le j})\eqno{(2)}$$
of topological spaces and continuous maps $\Ps ij\map{X_j}{X_i}$ ($i\le j$)
such that the following hold:
$$\eqalign{\hbox{for all }i,j,k\in I\qquad \Ps ik&=\Ps ij\circ\Ps jk\qquad(i\le j\le k)\cr
\hbox{and}\qquad\Ps ii&=\id{X_i}.\cr}$$

An {\it inverse limit} for (2) is a pair $(X,(\pi_i)_\ii)$ consisting of a topological
space, $X$, and continuous maps, $\pi_i\map X{X_i}$ ($\ii$), such that whenever
$(Y,(\rho_i)_\ii)$ is a topological space with continuous maps $\rho_i\map Y{X_i}$ ($\ii$)
and for all $i,j\in I$ with $i\le j$, $\rho_i=\Ps ij\circ\rho_j$, then there exists a unique,
continuous map $\rho\map YX$ such that for all $\ii$, $\pi_i\circ\rho=\rho_i$. We write
$$\il{}((X_i)_{i\in I},(\Ps ij)_{i\le j}):=(X,(\pi_i)_{i\in I}),$$
but occasionally omit reference to the maps.

\npar It should now be clear to the reader what direct and inverse systems and limits should
be in other categories. Indeed, as noted in $\ref{Pal}$, we can
obtain constructions for direct limits of, for example, normed spaces
by omitting reference to multiplication in (1).

We can now describe the second situation.

Suppose that
$$((A_i)_{i\in I}, (\tss ij)_{i\le j})\eqno(3)$$
is a direct system of uniform algebras, $(A_i,X_i)$, in which $\tss ij$ is given by
$f\mapsto f\circ\Ps ij$ for some continuous
surjection $\Ps ij\map{X_j}{X_i}$ where for all $i,j,k\in I$, $\Ps ik=\Ps ij\circ\Ps jk$
$(i\le j\le k)$ and $\Ps ii=\id{X_i}$. Thus
$$((X_i)_{i\in I}, (\Ps ij)_{i\le j})$$
is an inverse system of compact, Hausdorff spaces and continuous surjections.
It is well-known (see, for example, $\ref{Lei}$, p. 209)
that the inverse limit, $(X,(\pi_i)_{i\in I})$, of (2) exists.
For us it is sufficient to know that $X$ is a compact, Hausdorff space homeomorphic
to the non-empty set
$$X=\set{\k\in\times_{j\in I}X_j\st\hbox{ for all }i\le j\quad \k_i=\Ps ij(\k_j)}$$
with the relative product-topology. The canonical maps $(\pi_i)_{i\in\ I}$ are the
restrictions to $X$ of the coordinate projections $\times_{j\in I}X_j\to X_i$.

It is standard (see, for example, $\ref{Col}$)
that the Banach-algebra direct limit, $A$, of a system of uniform
algebras like (3) is identifiable with the uniform algebra on $X$ generated by
$\cup_{i\in I}\ps_{i}(A_i)$, and this is how such limits are formed in the text.
The canonical maps, $A_i\to A$, are given by the restrictions to $A_i$ of the adjoint
maps $\ps_i\mapto {C(X_i)}{C(X)}f{f\circ\pi_i}$ for each $\ii$.

\vfill\eject


\noindent {\bigtenbf Appendix 3}
\vskip 10pt
\noindent {\bigtenrm Algebra}
\vskip 25pt

\def\res{{\rm res}\,}

\noindent In this reference section we collect certain standard concepts and notation from algebra.
We assume familiarity with groups, vector spaces, and homomorphisms.

The second part introduces some less commonly known algebraic notions. These are
used in Chapter 4.

\vskip 15pt
\noindent {\medtenbf A.3.1. Rings and Algebras}
\vskip 15pt
\noindent
Recall that a {\it ring} is a non-empty set, $R$, together with two
associative, binary operations, $+$ and $\cdot$ (which is often written by
juxtaposition), such that $(R,+)$ is an abelian group and
$$\eqalign{& a\cdot(b+c)=(a\cdot b)+(a\cdot c)\cr
&(a+b)\cdot c=(a\cdot c)+(b\cdot c)\cr}\qquad\qquad(a,b,c\in R).$$

\dfn A.3.1.1. We say that a ring, $R$, is {\it unital} if there exists $1\in R$ such that
for all $r\in R$, $1\cdot r=r\cdot 1=r$. If $R$ is unital, then such an element, $1$, is unique and it
is called the {\it identity} of $R$. The ring is {\it commutative} provided
$$a\cdot b=b\cdot a\qquad(a,b\in R).$$

\dfn A.3.1.2. An {\it algebra} over $\C$, is a ring, $A$, which is also a
vector space over $\C$ such that
$$\l (a\cdot b)=(\l a)\cdot b=a\cdot (\l b)\qquad(a,b\in A,\ \l\in\C).$$
The algebra is called {\it commutative} if it is commutative as a ring.

\dfn A.3.1.3. An element, $e$, of a ring is called an {\it idempotent} if $e\cdot e =e$.
An idempotent is called {\it trivial} if and only if it equals $0$ or $1$.

\npar Algebras can be defined over fields other than $\C$; such algebras do not
appear in the thesis except in passing.

%
%
%
\dfn A.3.1.4. The {\it standard unitisation} of a complex algebra, $A$, is the
algebra $A_1:=A\times\C$ in which algebraic operations are given as follows:
$$\eqalign{(a_1,\l_1)+(a_2,\l_2)&=(a_1+a_2,\l_1+\l_2)\cr
 (a_1,\l_1)(a_2,\l_2)&=(a_1a_2+\l_2a_1+\l_1a_2,\l_1\l_2)\cr
 \l_1(a_2,\l_2)&=(\l_1a_2,\l_1\l_2)\cr }\qquad\qquad(a_1,a_2\in A,\ \l_1,\l_2\in\C).$$

\npar It is standard that $A_1$ has identity $(0,1)$, and that if $A$ is normed
(see Section 1.1.1) then $A_1$ is a normed algebra under
$$\norm{(a,\l)}_{A_1}:=\norm {a}_A+\ab\l \qquad\qquad(a\in A,\ \l\in\C).$$
%

\vskip 15pt
\noindent {\medtenbf A.3.2. Resultants}
\vskip 15pt

\noindent
Some of our proofs rely on resultants. 
In this section, $A$
will denote a commutative, unital ring.

Let $\a(x)=a_0+\cdots+
a_{n-1}x^{n-1}+x^n$ be a monic polynomial over $A$ and
$\b(x)=b_0+b_1x
\cdots+b_{n-1}x^{n-1}\in A[x]$.
By definition (see for example $\ref{JacI}$, p. 325) the resultant of
$\a(x)$ and $\b(x)$ is
$${\rm res}(\a(x), \b(x)):=\left\vert\matrix{
\overbrace{\matrix{1 &a_{n-1} &\cdots &a_2\cr
			0 &1 &\cdots &\cdots\cr
			\cdots  &\cdots &\cdots &\cdots\cr
			0  &\cdots &\cdots &1\cr}}^{n-1}
&\overbrace{\matrix{    a_1 &a_0 &0 &\cdots &0\cr
			a_2 &a_1 &a_0 &\cdots &0\cr
			\cdots &\cdots &\cdots &\cdots &\cdots \cr
			a_{n-1} &a_{n-2} &\cdots &\cdots &a_0 \cr}}^{n}	\cr
\matrix{b_{n-1} &b_{n-2} &\cdots &b_1\cr
	0 &b_{n-1} &\cdots &\cdots\cr
	\cdots &\cdots &\cdots &\cdots \cr
	0 &0 &\cdots &0\cr}
&\matrix{b_0 &0 &0 &\cdots &0\cr
	b_1 &b_0 &0 &\cdots &0\cr
	\cdots &\cdots &\cdots &\cdots &\cdots\cr
	b_{n-1} &b_{n-2} &\cdots &\cdots &b_0\cr}\cr}\right\vert$$
\vskip 5pt
It can be shown (see $\ref{Bou}$, A.IV, Section 6.6.1) that $\prin+\b(x)$ is
invertible in the quotient algebra $A[x]/\prin$ if and only if
$\res(\a(x),\b(x))$ is invertible in $A$.

We see from the above that, writing $c$ for $b_0$,
$\res(\a(x),\b(x))=P(c)=
p_0+p_1c+\cdots+p_{n-1}c^{n-1}+c^n$ where $p_0,\ldots,p_{n-1}\in A$
are polynomials in $b_1,\ldots,b_{n-1},a_0,\ldots,a_{n-1}$ only with
coefficients in $\C$.

It is a standard fact ($\ref{GKZ}$, p. 398)
that $\res(\a(x),\b(x))$ is homogeneous of
degree $n$ in $b_0,\ldots,b_{n-1}$ and homogeneous of degree $n-1$ in $a_0,
\ldots,a_{n-1},1$. Recall that this means that each term of the
polynomial is (up to a scalar multiple) of the form
$$b_0^{s_0}\cdots b_{n-1}^{s_{n-1}}a_0^{t_0}\cdots a_{n-1}^{t_{n-1}}1^{t_n}$$
for some $s_0,\ldots,s_{n-1}\in\No$ and $t_0,\ldots, t_n\in\No$ such that $\sum_{k=0}^{n-1}s_k=n$
and $\sum_{k=0}^{n}t_k=n-1$ .





\vfill\eject


\noindent {\bigtenbf References}
\vskip 4mm

\item{\ref{Ale}} Alexander, H., and Wermer, J. (1997) {\sl Several Complex Variables and Banach Algebras.} 3rd ed. New York: Springer Verlag.

\item{\ref{Are}} Arens, R. and Hoffman, K. (1956) Algebraic Extension of Normed Algebras. {\sl Proc.\ Am.\ Math.\ Soc.}, 7, 203-210

\item{\ref{Bol}} Bollob\'as, B. (1999) {\sl Linear Analysis.} 2nd ed. Cambridge: Cambridge University Press.

\item{\ref{Bou}} Bourbaki, N. (1988) {\sl Algebra II, Chapters 4-7.} United States of America: Springer-Verlag.

\item{\ref{Bro}} Browder, A. (1969) {\sl Introduction to Function Algebras.} New York: W. A. Benjamin, Inc.

\item{\ref{BroTop}} Browder, A. (2000) Topology in the Complex Plane. {\sl Am.\ Math.\ Monthly}, 107, no. 5, 393-401

\item{\ref{Bur}} Burckel, R.~B. (1972) {\sl Characterizations of C(X) Among its Subalgebras.} Lecture Notes in Pure and Applied Mathematics, 6, New York: Marcel Dekker Inc.

\item{\ref{Cir}} \v Cirka, E.~M. (1966) Approximating Continuous Functions by Holomorphic Ones on Jordan Arcs in $\C^n$. {\sl Soviet Math. Dokl.}, 7, 336-338

\item{\ref{Col}} Cole, B.~J. (1968) {\sl One-Point Parts and the Peak-Point Conjecture}, Ph.D.~Thesis, Yale University.

\vskip 11pt
\item{\ref{Cor}} Corach, G., and Su\' arez, F.~D. (1988) Thin Spectra and Stable Range Conditions. {\sl J.\ Funct.\ Anal.}, 81, 432-442

\item{\ref{Cou}} Countryman, R.~S. (1967) On the Characterization of Compact Hausdorff $X$ for which $C(X)$ is Algebraically Closed. {\sl Pacif.\ J.\ Math.}, 20, no. 3, 433-448

\item{\ref{Dal}} Dales, H.~G. (2000) {\sl Banach Algebras and Automatic Continuity.} New York: Oxford University Press Inc.

\item{\ref{4DIS}} Dawson, T.~W. (2000) {\sl Algebraic Extensions of Normed Algebras}, M.Math.\ Dissertation, University of Nottingham, accessible from
$${\tt http://xxx.lanl.gov/abs/math.FA/0102131}$$

\item{\ref{DawS}} Dawson, T.~W. (2002) A Survey of Algebraic Extensions of Commutative, Unital Normed Algebras. {\sl Contemporary\ Math.} (to appear)

\item{\ref{DF}} Dawson, T.~W., and Feinstein, J.~F. (2002) On the Denseness of the Invertible Group in Banach Algebras. {\sl Proc.\ Am.\ Math.\ Soc.} (to appear).

\item{\ref{Fal}} Falc\'on Rodr\'\i guez, C. M. (1988) Sobre la densidad del grupo de los elementos invertibles de un \'algebra
uniforme. {\sl Revista Ciencias Matem\'aticas}, 9, no. 2, 11-17

\item{\ref{FeiThes}} Feinstein, J.~F. (1989) {\sl Derivations from Banach Function Algebras}, Ph.D.\ Thesis, University of Leeds

\item{\ref{FeiNTSR}} Feinstein, J.~F. (1992) A Non-Trivial, Strongly Regular Uniform Algebra. {\sl J.\ Lond.\ Math.\ Soc.}, 45, no. 2, 288-300

\item{\ref{GelSCF}} Gelfand, I.~M. (1960) On the Subrings of a Ring of Continuous Functions. {\sl \ Am.\ Math.\ Soc.\ Translations}, Series 2, 16

\item{\ref{GKZ}} Gelfand, I.~M., Kapranov, M.~M., and Zelevinsky, A.~V. (1994) {\sl Discriminants, Resultants, and Multidimensional Determinants.} United States of America: Birkh\" auser Boston.

\item{\ref{Gor}} Gorin, E.~A., and Lin, V.~J. (1969) Algebraic Equations with Continuous Coefficients and Some Problems of the Algebraic Theory of Braids. {\sl Math.\ USSR Sb.}, 7, no. 4, 569-596

\item{\ref{Gri}} Grigoryan, S.~A. (1984) Polynomial Extensions of Commutative Banach Algebras. {\sl Russian Math.\ Surveys}, 39, no. 1, 161-162

\item{\ref{Hal}} Halmos, P.~R. (1970) {\sl Naive Set Theory.} New York: Springer-Verlag.

\item{\ref{Hat}} Hatori, O., and Miura, T. (1999) On a Characterization of the Maximal Ideal Spaces of Commutative $C$*-Algebras in Which Every Element is the Square of Another. {\sl Proc.\ Am.\ Math.\ Soc.}, 128, no. 4, 1185-1189

\item{\ref{Heu}} Heuer, G.~A., and Lindberg, J.~A. (1963) Algebraic Extensions of Continuous Function Algebras. {\sl Proc.\ Am.\ Math.\ Soc.}, 14, 337-342

\item{\ref{HofRam}} Hoffman, K., and Ramsay, A. (1965) Algebras of Bounded Sequences. {\sl Pacif.\ J.\ Math.}, 15, no. 4, 1239-1248

\item{\ref{Hos}} Host-Madsen, A. (2000) Separable Algebras and Covering Spaces. {\sl Math.\ Scand.}, 87, 211-239

\item{\ref{JacI}} Jacobson, N. (1996) {\sl Basic Algebra I.} 2nd ed. New York: W.H. Freeman and Company.

\item{\ref{Kal}} Kallin, E. (1963) A Nonlocal Function Algebra. {\sl Proc.\ Nat.\ Acad.\ Sci.\ U.S.A.}, 49, 821-824

\item{\ref{Kar}} Karahanjan, M.~I. (1979) Some Algebraic Characterisations of the Algebra of All Continuous Functions on a Locally Connected Compactum. {\sl Math.\ USSR Sb.}, 35, 681-696

\item{\ref{Kel}} Kelley, J.~L. (1955) {\sl General Topology.} New-York: Van Nostrand.

\item{\ref{Lei}} Leibowitz, G.~M. (1970) {\sl Lectures on Complex Function Algebras.} United States of America: Scott, Foresman and Company

\item{\ref{LinAE}} Lindberg, J.~A. (1964) Algebraic Extensions of Commutative Banach Algebras. {\sl Pacif.\ J.\ Math.}, 14, 559-583

\item{\ref{LinFact}} Lindberg, J.~A. (1964) Factorization of Polynomials over Banach Algebras. {\sl Trans.\ Am.\ Math.\ Soc.}, 112, 356-368

\item{\ref{LinEANA}} Lindberg, J.~A. (1971) Extension of Algebra Norms and Applications. {\sl Studia Math.}, 40, 35-39

\item{\ref{LinIE}} Lindberg, J.~A. (1973) Integral Extensions of Commutative Banach Algebras. {\sl Can.\ J.\ Math.}, 25, 673-686

\item{\ref{McK}} McKissick, R. (1963) A Nontrivial Normal Sup Norm Algebra. {\sl Bull.\ Am.\ Math.\ Soc.}, 69, 391-395

\item{\ref{Mur}} Murphy, G.~J. (1990) {\sl C*-Algebras and Operator Theory.} Boston: Academic Press.

\item{\ref{Nar}} Narmaniya, V.~G. (1982) The Construction of Algebraically Closed Extensions of Commutative Banach Algebras. {\sl Trudy Tbiliss.\ Mat.\ Inst.\ Razmadze Akad.}, 69, 154-162

\item{\ref{Pal}} Palmer, T.~W. (1994) {\sl Banach Algebras and the General Theory of *-Algebras.} (Vol. 1) Cambridge: Cambridge University Press.

\item{\ref{Pea}} Pears, A.~R. (1975) {\sl Dimension Theory of General Spaces.} Cambridge: Cambridge University Press.

\item{\ref{Ped}} Pedersen, G.~K. (1995) {\sl Analysis Now.} (revised printing) New York: Springer-Verlag.

\item{\ref{Rao}} Rao, M., and Stetk\ae r, H. (1991) {\sl Complex Analysis, An Invitation.} Singapore: World Scientific Publishing Co. Pte. Ltd.

\item{\ref{RicMIS}} Rickart, C.~E. (1966) The Maximal Ideal Space of Functions Locally Approximable in a Function Algebra. {\sl Proc.\ Am.\ Math.\ Soc.}, 17, 1320-1327

\item{\ref{Rie}} Rieffel, M.~A. (1983) Dimension and Stable Rank in the K-Theory of C*-Algebras. {\sl Proc.\ Lond.\ Math.\ Soc.}, 46, no. 3, 301-333

\item{\ref{Rot}} Rotman, J.~J. (1979) {\sl An Introduction to Homological Algebra.} London: Academic Press Inc.

\item{\ref{RudRCA}} Rudin, W. (1987) {\sl Real and Complex Analysis.} 3rd ed. Singapore: McGraw-Hill Book Co.

\item{\ref{RudFA}} Rudin, W. (1991) {\sl Functional Analysis.} 2nd ed. Singapore: McGraw-Hill Book Co.

\item{\ref{Rud}} Rudin, W. (1997) {\sl The Way I Remember It.} United States of America: American Mathematical Society.

\item{\ref{Sid}} Sidney, S.~J. (1971) High-Order Non-Local Uniform Algebras. {\sl Proc.\ Lond.\ Math.\ Soc.}, 14, 735-752

\item{\ref{Sim}} Simmons, G.~F. (1963) {\sl Introduction to Topology and Modern Analysis.} Singapore: McGraw-Hill Book Co.

\item{\ref{StoCon}} Stolzenberg, G. (1963) Polynomially and Rationally Convex Sets. {\sl Acta Math.}, 109, 259-289

\item{\ref{StoMIS}} Stolzenberg, G. (1963) The Maximal Ideal Space of the Functions Locally in a Function Algebra. {\sl Proc.\ Am.\ Math.\ Soc.}, 14, 342-345

\item{\ref{Sto}} Stout, E.~L. (1973) {\sl The Theory of Uniform Algebras.} Tarrytown-on-Hudson, New York: Bogden and Quigley Inc.

\item{\ref{Tay}} Taylor, J.~L. (1975) Banach Algebras and Topology, In: Williamson, J.~H. (ed.)
{\sl Algebras in Analysis.} Norwich: Academic Press Inc. (London) Ltd. 118-186

\item{\ref{Zam}} Zame, W.~R. (1984) Covering Spaces and the Galois Theory of Commutative Banach Algebras. {\sl J.\ Funct.\ Anal.}, 27, 151-171

\vfill\eject

\end